\newif\iffinal
\finaltrue	

\documentclass[11pt,bezier,amstex]{article}  

\iffinal
\else\usepackage[notref,notcite]{showkeys}
\else\usepackage[notref]{showkeys}

\fi

\usepackage{color}              
\usepackage{epsfig}
\usepackage{graphicx}
\usepackage[]{amsmath}
\usepackage{amssymb}
\usepackage{epsfig}
\usepackage{hyperref}
\usepackage{amsfonts}
\usepackage{pstricks,pst-node,pst-tree}
\definecolor{MyDarkBlue}{rgb}{0,0.08,0.50}
\definecolor{BrickRed}{rgb}{0.65,0.08,0}
\usepackage[textsize=small]{todonotes}       

\newcommand{\BF}{{\bf {BF}}}
\newcommand{\BP}{{\bf {BP}}}

\newcommand{\DD}{\mathcal{D}}
\newcommand{\CC}{\mathcal{C}}
\newcommand{\JJ}{\mathcal{J}}
\newcommand{\ti}{\tilde}
\newcommand{\thtpar}{\theta}
\newcommand{\betapar}{\vartheta}

\newcommand{\clc}{\mathcal{C}}
\newcommand{\clh}{\mathcal{H}}
\newcommand{\clt}{\mathcal{T}}
\newcommand{\clu}{\mathcal{U}}
\newcommand{\cle}{\mathcal{C}}
\newcommand{\clv}{\mathcal{V}}
\newcommand{\clg}{\mathcal{G}}

\newcommand{\KK}{\mathcal{K}}
\newcommand{\VV}{\mathcal{V}}
\newcommand{\NNN}{\mathbb{N}}
\newcommand{\RRR}{\mathbb{R}}
\newcommand{\cals}{\mathcal{S}}
\newcommand{\cla}{\mathcal{A}}
\newcommand{\squig}{\leftrightarrow}
\newcommand{\gampar}{\gamma}
\newcommand{\IA}{{\bf{IA}}}
\newcommand{\RG}{{\bf{RG}}}
\newcommand{\XX}{\mathcal{X}}
\newcommand{\bfx}{{\bf{x}}}
\newcommand{\bfy}{{\bf{y}}}
\newcommand{\bfG}{{\bf{G}}}

\newcommand{\be}{\begin{equation}}
\newcommand{\ee}{\end{equation}}
\newcommand{\beq}{\begin{eqnarray*}}
\newcommand{\eeq}{\end{eqnarray*}}

\newcommand{\beqn}{\begin{eqnarray}}
\newcommand{\eeqn}{\end{eqnarray}}

\newcommand{\ba}{\begin{aligned}}
\newcommand{\ea}{\end{aligned}}
\newcommand{\bes}{\begin{equation*}}
\newcommand{\ees}{\end{equation*}}

\newtheorem{Lemma}{Lemma}[section]
\newtheorem{Proposition}[Lemma]{Proposition}

\newtheorem{Theorem}[Lemma]{Theorem}

\newtheorem{Corollary}[Lemma]{Corollary}

\newcommand{\calC}{\mathcal{C}}
\newcommand{\EE}{\mathcal{E}}
\newcommand{\bfe}{{\bf e}}
\newcommand{\prob}{\mathbb{P}}
\newcommand{\PP}{\mathcal{P}}
\newcommand{\E}{\mathbb{E}}

\newcommand{\Acal}{\mathcal{A}}

\newcommand{\FF}{\mathcal{F}}
\newcommand{\GG}{\mathcal{G}}

\newcommand{\TT}{\mathcal{T}}
\newcommand{\WW}{\mathcal{W}}
\newcommand{\eps}{\varepsilon}

\newcommand{\qed}{\ \ \rule{1ex}{1ex}}

\newcommand{\set}[1]{\left\{#1\right\}}

\newcommand{\Rbold}{{\mathbb{R}}}

\newcommand{\expec}{\mathbb{E}}

\newcommand{\bfd}{{\bf d}}
\newcommand{\bfzero}{{\bf 0}}

\newcommand{\bfC}{\boldsymbol{C}}
\newcommand{\bfX}{\boldsymbol{X}}
\newcommand{\ldown}{l^2_{\downarrow}}
\newcommand{\calS}{\mathcal{S}}
\newcommand{\bars}{\bar{s}}

\newcommand{\sss}{\scriptscriptstyle}

\newcommand{\barx}{\bar{x}}
\def\1{{\mathchoice {1\mskip-4mu\mathrm l}      
{1\mskip-4mu\mathrm l}
{1\mskip-4.5mu\mathrm l} {1\mskip-5mu\mathrm l}}}

\newcommand{\com}{\mathcal{COM}}
\newcommand {\convd}{\stackrel{d}{\longrightarrow}}
\newcommand {\convp}{\stackrel{\sss {\mathbb P}}{\longrightarrow}}

\numberwithin{equation}{section}

\author{Shankar Bhamidi\thanks{{\tt bhamidi@email.unc.edu}} \and Amarjit Budhiraja\thanks{{\tt budhiraj@email.unc.edu}} \and Xuan Wang\thanks{{\tt wangxuan@email.unc.edu}} \\
Department of Statistics and Operations Research,\\
304 Hanes Hall, \\
University of North Carolina, \\
Chapel Hill, NC.
}

\title{Bohman-Frieze processes at criticality and emergence of the giant component}

\begin{document}
\maketitle

\begin{abstract}
	The evolution of the usual Erd\H{o}s-R\'{e}nyi random graph model on $n$ vertices can be described as follows: At time $0$ start with the empty graph, with $n$ vertices and no edges. Now at each time $k$, choose 2 vertices uniformly at random and attach an edge between these two vertices. Let $\bfG_n(k)$ be the graph obtained at step $k$.   Following  \cite{er-1} and \cite{janson1994birth}, the work of Aldous in \cite{aldous1997brownian} shows that for any fixed $t\in \Rbold$,  when $k(n) = n/2+ n^{2/3} t/2$, the sizes of the components in $\bfG_n(k(n))$ scale like $n^{2/3}$ and, as $n \to \infty$, 
	the scaled component  size vector converges to  the standard multiplicative coalescent at time $t$.  Furthermore, when $k(n) = sn/2$, for any $s< 1$, the size of the largest component is $O(\log{n})$ while for  $s> 1$, the size of the largest component scales like $f(s)n$, where $f(s)> 0$.  Thus $t_c = 1$ is the critical parameter for this model at which time a giant component emerges.
	
	The last decade has seen variants of this process introduced, loosely under the name {\it Achlioptas processes}, to understand the effect of simple changes in the edge formation  scheme on the emergence of the giant component. Stimulated by a question of Achlioptas, one of the simplest and most popular of such models is the Bohman-Frieze (BF) model wherein at each stage $k$,  2 edges $e_1(k)=(v_1,v_2)$ and $e_2(k) = (v_3, v_4)$ are chosen uniformly at random. If at this time $v_1, v_2$ are both isolated then this edge is added, otherwise $e_2$ is added. Then \cite{bohman2001avoiding} (and further analysis in \cite{spencer2007birth}) show that  once again there is a critical parameter, which is larger than $1$, above and below which the asymptotic behavior is as in the Erd\H{o}s-R\'{e}nyi setting, 
	in particular this simple modification in the attachment scheme delays the time of emergence of the giant component.
	 While an intense study for this and related models seems to suggest that at criticality, this model should be in the same universality class as  the original Erd\H{o}s-R\'{e}nyi process, a precise mathematical treatment
	 of the dynamics in the critical window has to date escaped analysis. In this work we study the component structure of the BF model in the critical window
	and show that  at criticality the sizes of components properly rescaled and re-centered converge to the same limits as for the original Erd\H{o}s-R\'{e}nyi process. Our proofs rely on a series of careful approximations of the original random graph model leading eventually to the treatment of certain near critical
	 multitype branching processes with general type spaces through functional and stochastic analytic tools. 
	
	
\end{abstract}

\noindent
{\bf Key words:} critical random graphs, multiplicative coalescent, entrance boundary, giant component, branching processes, inhomogeneous random graphs. 

\noindent
{\bf MSC2000 subject classification.} 60C05, 05C80, 90B15.



\section{Introduction}
\label{sec-int}

Random graph models of various systems in the real world have witnessed a tremendous growth in the last decade. The availability of a large amount of empirical data on real world networks and systems such as road and rail networks, bio-chemical networks, social networks and data transmission networks such as the internet has stimulated an inter-disciplinary effort to formulate models to understand such systems ranging from (largely) static systems such as road and rail networks, to more dynamic systems such as social networks, information networks  and models of coagulation and aggregation in physics and colloidal chemistry. 

The classical Erd\H{o}s-R\'{e}nyi(ER) random graph can be thought of as a network evolving in time via the following prescription: Start with the empty graph at time zero. Then at each discrete time point, choose an edge uniformly at random and place it in the network (irrespective of whether it was present before or not). This simple model exhibits a tremendous amount of complexity, see \cite{bollobas-rg-book,janson-luczak-bk}. As one adds more and more edges, eventually the system transitions from a ``sub-critical'' regime wherein the largest component is of order $\log{n}$, to the super-critical regime wherein there exists a unique giant component of order $n$ and the next largest component is $O(\log{n})$. Understanding the emergence of this giant component and the properties of near critical Erd\H{o}s-R\'{e}nyi random graph have stimulated an enormous amount of work, see e.g \cite{nachmias2007component,bollobas2010asymptotic,ding2010diameters} and the references therein.

 This model has been modified in various ways the last few years, to understand the effect of choice in the emergence of the giant component. For example, in \cite{achlioptas2009explosive} and \cite{d2010local} simulation based studies are carried out for an attachment scheme where at each stage one chooses two edges uniformly at random,  and then one uses the edge which delays the emergence of a giant component (for example by choosing the edge that minimizes the product of the two adjoining components). A number of fascinating conjectures about the ``explosive emergence'' of the giant component are stated in these articles.

Such modifications of the ER random graph model have generated a tremendous amount of interest in many different communities, ranging from combinatorics to statistical physics. Although a number of such models have been analyzed in the sub-critical and super-critical regime, see for example the marvelous and comprehensive \cite{spencer2007birth}, understanding how the giant component emerges, even for very simple modifications to the ER setting, for example for the  Bohman-Frieze model (which we describe below) has escaped a mathematical understanding. The aim of this paper is to prove the first rigorous results about such models at criticality and obtain precise asymptotics for the merging dynamics in the critical window that lead to the emergence of the giant component.

The Bohman-Frieze (BF) model is as follows: Let ${\bf 0}_n$ be the empty graph on the vertex set $[n] \stackrel{\sss def}{=} \{1,2, ... , n\}$ with no edges. The random graph process $\{\bfG_n^{\sss BF}(k)\}_{k \in \mathbb{N}}\equiv \{\bfG_n(k)\}_{k \in \mathbb{N}}$ evolves as follows:
\begin{itemize}
	\item At time $k=0$ let $\bfG_n(0) = {\bf0}_n$,
	\item The process evolves in discrete steps. For $k \geq 0$ let the state of the graph at time $k$ be $\bfG_n(k)$. Then at time $k+1$, the graph $\bfG_n(k+1)$ is constructed as follows: Choose two edges $e_1=
	(v_1, v_2)$ and $e_2 = (v_3, v_4)$ uniformly at random amongst all the possible ${n\choose 2}$ edges. If $v_1, v_2$ are isolated vertices then let $\bfG_n(k+1)$ be the graph $\bfG_n(k)$ with edge $e_1$ added else let $\bfG_n(k+1)$ be the graph $\bfG_n(k)$ with the edge $e_2$ added.
\end{itemize}
Note that the standard Erd\H{o}s-R\'{e}nyi (ER) process can be thought of as a variant of the above process wherein we add edge $e_1$ at each stage (irrespective of whether the end points are isolated vertices or not). We shall use $\bfG_n^{\sss ER}(k)$ to denote the usual Erd\H{o}s-R\'{e}nyi process. 

For the standard Erd\H{o}s-R\'{e}nyi process classical results (see e.g. \cite{er-1}, \cite{er-2}, \cite{bollobas-rg-book}) tell us that by scaling time by $n$ and looking at the process $\bfG_n^{\sss ER}(\lfloor nt/2 \rfloor) $, a phase transition occurs at $t_c(er)=1$, namely
\\(a) {\bf Subcritical regime:} For fixed $t<1 $, the largest component in $\bfG_n^{\sss ER}(\lfloor nt/2 \rfloor)$ is $O(\log{n})$.
\\(b) {\bf Supercritical regime:} For fixed $t> 1$, the largest component has size $\sim f(t)n$ for some positive function $f$ (namely there is a giant component) while the second largest component is of size $O(\log{n})$.
\\(c) {\bf Critical regime: } For $t=1$, the first and the second largest components both are $\Theta(n^{2/3})$.

For the BF model, the original work in \cite{bohman2001avoiding} shows that there exists $t_0 > 1 $ such that the size of the largest component in $\bfG_n^{\sss BF}( \lfloor nt_0/2 \rfloor)$ is $o_p(n)$. Thus this model ``avoids'' the giant component for more time than the original Erd\H{o}s-R\'{e}nyi model. The comprehensive paper \cite{spencer2007birth} showed for  the BF model and for a number of other related processes,  that there exists a critical value $t_c$,  such that conclusions (a) and (b) hold for these models as well with $1$ replaced
by $t_c$. For the BF model, the supercritical regime, namely the regime $t> t_c\equiv t_c(bf)$ was studied further in \cite{janson2010phase}. Numerical calculations from this paper show that $t_c(bf)\approx 1.17$.

Understanding what happens in the critical regime for these models is of tremendous interest for many different fields, as it turns out that much of the ``action'' happens at criticality. For the Erd\H{o}s-R\'{e}nyi processes, Aldous's result in \cite{aldous1997brownian} shows  that   if we denote the component size vector(components arranged in decreasing order), in  the graph $\bfG_n^{\sss ER}(\lfloor nt/2 \rfloor)$, by 
\[\bfC_n^{\sss ER}(t)=(\calC_n^{\sss(i)}(t): i\geq 1),\]
then, for any fixed $\lambda\in \Rbold$, as $n \to \infty$, 
\be \bar \bfC_n^{\sss ER}(\lambda) \stackrel{\sss def}{=} n^{-2/3}\bfC_n^{\sss ER}\left( 1+\frac{\lambda}{n^{1/3}}\right)\convd \bfX(\lambda),\label{ins1629n}\ee
 where the process $(\bfX(\lambda): -\infty< \lambda < \infty)$ is a Markov process called the {\it eternal standard multiplicative coalescent} on the space 
\begin{equation}
	\ldown = \{(x_1,x_2,\ldots): x_1\geq x_2\geq \cdots \geq 0, \sum_i x_i^2< \infty\}.
	\label{eqn:ldown}
\end{equation}
The space $\ldown$ is endowed with the same topology as inherited from $l^2$, and the above $\convd$ denotes weak convergence in this space.
  In fact the paper \cite{aldous1997brownian} shows weak convergence of $\bar \bfC_n^{\sss ER}$ to $\bfX$
in $\DD((-\infty, \infty), \ldown)$ (the space of RCLL functions from $(-\infty, \infty)$ to $\ldown$ endowed with the usual Skorohod topology).
We defer a full definition of the above Markov process to Section \ref{sec:results}.



Despite the tremendous amount of interest in the critical regime of such models and the significant progress in analyzing these models above and below criticality (see e.g. \cite{spencer2007birth} for a general treatment of such models and \cite{janson2010phase} for an analysis of the BF model), analysis of the component sizes at criticality and an understanding of the dynamic behavior of these components and how they merge through the critical scaling window for the giant component to emerge, for settings beyond the classical Erd\H{o}s-R\'{e}nyi model has remained an extremely challenging problem. The most commonly used technique in getting refined results about the sizes of the largest component at criticality,  which is to show that the breadth-first exploration processes of components converge to a certain `inhomogeneous Brownian motion' does not seem to extend to such settings and one needs to develop a rather different machinery.   In this work we take a significant step towards the understanding of the  behavior of such models by
treating the asymptotics in the critical regime for the Bohman-Frieze process. Although this model makes a simple modification to the attachment scheme for the Erd\H{o}s-R\'{e}nyi  setting, the analysis
requires many new constructions and mathematical ideas. Our main result, Theorem \ref{theo:main}, shows that at criticality the sizes of components properly rescaled and re-centered converge, as  for the original Erd\H{o}s-R\'{e}nyi process, to the eternal standard multiplicative coalescent.   The proof proceeds through a sequence of approximations that reduce the study to that for a certain inhomogeneous random graph model which allows us to bring to bear techniques from the theory of multitype branching processes on general type spaces and integral operators on Hilbert spaces.  These culminate in an estimate, obtained in Proposition \ref{prop:com-size}, on the size of the largest component through the ``subcritical" window.  Using this estimate and other careful estimates on quadratic variations of certain martingales associated with the Bohman-Frieze process we then obtain (see Proposition \ref{prop:main}) precise asymptotics on the sum of squares and cubes of component sizes just below the critical window.  The latter result is a key ingredient to the proof as it allows us to show that just before the critical window, the component sizes satisfy the regularity conditions required for convergence to the standard multiplicative coalescent, and thus allow us to couple the BF process with a standard multiplicative coalescent through the critical window to prove the main result. 
Although not pursued here,  the techniques developed in this work for establishing
the asymptotic behavior in Proposition \ref{prop:main} are potentially of use for the analysis of the critical regime for a very general class of Achlioptas processes (called bounded size rules).

{\bf Organization of the paper:} Due to the length of the paper, let us now briefly guide the reader through rest of this work. We begin in Section \ref{sec:cont-time-bf}
 with the construction of the continuous time version of the BF-process and then state the main result (Theorem \ref{theo:main}) in Section \ref{sec:results}. 
 Next, in Section  \ref{sec:disc} we give a wide ranging discussion of related work and the relevance of our results. We shall then start with the proof   of the main result by giving an intuitive sketch in Section \ref{sec:proof-idea} wherein we also provide details on the organization of the various steps in the proof carried out in Sections 5-7.  Finally Section \ref{sec:proof-main} combines all these ingredients to complete the proof.


\section{The Bohman-Frieze process}
\label{sec:cont-time-bf}
Denote the vertex set by $[n]=\{1,2,\ldots, n\}$ and the edge set by $\EE_n=\{\{v_1,v_2\}: v_1\neq v_2\in [n]\}$. To simplify notation we shall suppress $n$ in the notation unless required.  Denote by $\BF(t)=\BF_n(t)$, $t \in [0,\infty)$,  the continuous time  Bohman-Frieze random graph process, constructed as follows:\\
Let $\EE^2=\EE \times \EE$ be the set of all ordered pairs of edges. For every ordered pair of edges $\bfe = (e_1,e_2) \in \EE^2$ let $\PP_\bfe$ be a Poisson process on $[0,\infty)$ with rate $2/n^3$, and let these processes be independent as $\bfe$ ranges over $\EE^2$. We order the points generated by all the ${n \choose 2} \times  {n \choose 2}$ Poisson processes by their natural order as $0 < t_1<t_2<...$. Then we can define the BF-process iteratively as follows:\\
(a) When $t \in [0, t_1) $, $\BF(t)={\bf 0}_n$, the empty graph with $n$ vertices; \\
(b) Consider $t \in [t_k, t_{k+1})$, $k \in \mathbb{N}^+$, where $t_k$ is a point in $\PP_\bfe$ and $\bfe = (e_1,e_2) =(\{v_1,v_2\}, \{v_3,v_4\})$. 
 If $v_1$, $v_2$ are both singletons (i.e. not connected to any other vertex) in $\BF(t_k-)$, then $\BF(t)=\BF(t_k-) \cup \{e_1\}$, else let $\BF(t)=\BF(t_k-) \cup \{e_2\}$. \\
Note that multiple edges are allowed between two given vertices, however this has no significance in our analysis which
is primarily concerned with  component sizes.\\

Consider the same construction but with the modification that we always add $e_1$ to the graph and disregard the second edge $e_2$. Note that the total rate of adding new edges is:
$${n\choose 2}\times {n \choose 2} \frac{2}{n^3} \approx \frac{n}{2}.$$
Then this random graph process is just a continuous time version of the standard Erd\H{o}s-R\'{e}nyi process 
and the convergence in (\ref{ins1629n}) continues to hold, where $\bfC_n^{\sss ER}(t)$ now represents the component size vector at time $t$ for this continuous
time ER process and once more, $t_c=1$ is the critical parameter for the model.

 As proved in \cite{spencer2007birth}, the Bohman-Frieze model also displays a phase transition and the critical time $t_c \approx 1.1763$.  We now summarize
some results from the latter paper that characterize this critical parameter in terms of the behavior of certain differential equations.


The following notations and definitions mostly follow \cite{janson2010phase}. Let $\CC_n^{\sss (i)}(t)$ denote the size of the $i^{th}$ largest component in $\BF_n(t)$, and $\bfC_n(t)=(\CC_n^{\sss (i)}(t) : i \ge 1)$ the component size vector. For convenience, we define $\CC_n^{\sss (i)}(t)=0$ whenever $t < 0$. 

For fixed time $t \ge 0$, let $X_n(t)$ denote the number of singletons at this time and $\bar{x}(t)=X_n(t)/n$ denote the density of singletons. For simplicity,
we have suppressed the dependence on $n$ in the notation. For $ k = 2,3$, let
\be \calS_k(t)=\sum_{i \ge 1} (\CC_n^{\sss (i)}(t))^k \label{ins2141}\ee
and let $\bar{s}_k(t)=\calS_k(t)/n$. Then from \cite{spencer2007birth}, there exist deterministic functions $x(t), s_2(t), s_3(t)$ such that for each fixed $t\geq 0$:
$$\bar{x}(t)\convp x(t),\qquad \bar{s}_k(t) \convp s_k(t) \qquad \mbox{for } k=2,3,$$
as $n\to\infty$. The limiting function $x(t)$ is continuous and differentiable for all $t\in \Rbold_+$. For $k\geq 2$, there exists $1<t_c< \infty$ such that
 $s_k(t)$ is finite, continuous and differentiable for $0\le t< t_c$, and $s_k(t) = \infty$ for $t\geq  t_c$. Furthermore, $x, s_2, s_3$ solve the following differential equations.
 \begin{align}
	x'(t)&=-x^2(t)-(1-x^2(t))x(t) \qquad &\mbox{ for } t \in [0,\infty,) \qquad x(0)=1
	\label{eqn:xss-def}\\
	s_2'(t)&=x^2(t)+(1-x^2(t))s_2^2(t)  \qquad &\mbox{ for } t \in [0,t_c), \qquad s_2(0)=1\label{ins-s2}\\
	s_3'(t)&=3x^2(t)+3(1-x^2(t))s_2(t)s_3(t)  \qquad &\mbox{ for } t \in [0,t_c), \qquad s_3(0)=1.\label{ins-s3}
\end{align}
This constant $t_c = t_c(bf)$ is the critical time such that whp, for $t< t_c$, the size of the largest component in $\BF(t)$ is $O(\log{n})$, while for $t> t_c$ there exists a giant component of size $\Theta(n)$ in $\BF(t)$.  Furthermore from \cite{janson2010phase} (Theorem 3.2) there exist constants \begin{align}
	\alpha &= (1-\barx^2(t_c))^{-1} \approx 1.063 \label{eqn:alpha-def}\\
	\beta &=\lim_{t\uparrow t_c} \frac{s_3(t)}{[s_2(t)]^3}\approx .764 \label{eqn:beta-def}
\end{align} 
such that as $t\uparrow t_c$
\begin{align}
	s_2(t) &\sim \frac{\alpha}{t_c-t} \label{eqn:s2-scaling-crit} \\
	s_3(t) &\sim \beta (s_2(t))^3 \sim \beta \frac{\alpha^3}{(t_c-t)^3}. \label{eqn:s3-scaling-crit}
\end{align}



\subsection{Main results}
\label{sec:results}
Our goal in this work is to establish a limit theorem of the form in (\ref{ins1629n}) for the BF process.

We begin with a precise definition of the standard eternal multiplicative coalescent process $\bfX$ introduced in Section \ref{sec-int}. 
For $x \in l^2$, let $\mbox{ord}(x) \in \ldown$ be a reordering of the vector $x$ in the decreasing order. For $x\in \ldown$,  $1\le i< j<\infty$, let 
\[x^{ij}=\mbox{ord}(x + \langle x,e_j\rangle (e_i-e_j)),\] where $\{e_i\}$ is the canonical basis in $l^2$. Namely, $x^{ij}$ is the vector obtained by merging
the $i$-th and $j$-th  `components' in $x$ and reordering the vector.  Aldous(\cite{aldous1997brownian}) showed that there is a Feller Markov process
with sample paths in $D([0,\infty):\ldown)$ with infinitesimal generator
\be \mathcal{A}_{\mbox{\tiny{MC}}} f(x) = \sum_{i<j}x_ix_j (f(x^{ij})-f(x)), \; x \in \ldown, \; f: \ldown \to \RRR.\label{ins1723n}\ee
This Markov process describes a coalescence dynamics where at any time instant the $i$-th and $j$-th clusters merge at rate equal to the product of the
sizes of the two clusters.  One special choice of initial distribution for this Markov process is particularly relevant for the study of asymptotics
of random graph models.  We now describe this distribution. Let $\{W(t)\}_{t\ge 0}$ be a standard Brownian motion, and for a fixed $\lambda \in \Rbold$, define
\[W_\lambda(t) = W(t)+\lambda t-\frac{t^2}{2},\; t \ge 0.\]
Let $\bar W_{\lambda}$ denote the reflected version of $W_{\lambda}$, i.e., 
\begin{equation}
	\bar{W}_\lambda(t) = W_\lambda(t) - \min_{0\leq s\leq t} W_\lambda(s), \; t \ge 0.
	\label{eqn:inh-ref-bm}
\end{equation}
 Define an excursion of $\bar W_\lambda$ as an interval $(l,u) \subset [0,+\infty)$ such that $\bar W_\lambda(l)=\bar W_\lambda(u)=0$ and $\bar W_\lambda(t)>0$ for all $t \in (l,u)$. Define $u-l$ as the size of the excursion. Order the sizes of excursions of $\bar W_\lambda$ as 
\[\xi_1(\lambda)> \xi_2(\lambda)> \xi_3(\lambda)> \cdots\] and write $\Xi(\lambda) = (\xi_i(\lambda):i\geq 1).$
Then $\Xi(0)$ defines a $\ldown$ valued random variable. Denote by $\{\bfX(\lambda)\}_{\lambda \ge 0}$ the $\ldown$ valued Markov process with initial
distribution as the probability law of $\Xi(0)$ and infinitesimal generator as in  (\ref{ins1723n}). Then \cite{aldous1997brownian} shows that
for each $\lambda \in (0,\infty)$, $\bfX(\lambda)$ has the same law as $\Xi(\lambda)$.  In fact, \cite{aldous1997brownian} shows that the process $\bfX$ can be extended
for $\lambda \in (-\infty, \infty)$, namely there is a process with sample paths in $\DD((-\infty,\infty):\ldown)$, denoted once more as $\bfX$, such that for every
$\lambda \in (-\infty, \infty)$, $\bfX(\lambda) =_d \Xi(\lambda)$ and $\{\bfX(\lambda+t)\}_{t\ge 0}$ is a Markov process with generator $\mathcal{A}_{\mbox{\tiny{MC}}}$.  This description
uniquely characterizes a probability measure on $\DD((-\infty,\infty):\ldown)$ representing the probability law of $\bfX$. The process $\bfX$ is called
the {\em eternal standard multiplicative coalescent}.

We are now in a position to state our main result.
For the rest of this work $t_c=t_c(bf)$ will denote the critical point for the continuous time Bohman-Frieze process as defined above.  Recall the constants $\alpha$ and $\beta$  defined in \eqref{eqn:alpha-def} and \eqref{eqn:beta-def}.
\begin{Theorem}
\label{theo:main}
For $\lambda \in \Rbold$,  let
\begin{equation}
	\label{eqn:res-com-def}
	{\bar \bfC}_n^{\sss BF}(\lambda) = \left(\frac{\beta^{1/3}}{n^{2/3}}\calC^{\sss (i)}_n \left(t_c+ \beta^{2/3} \alpha\frac{\lambda}{n^{1/3}} \right): i\geq 1\right) 
\end{equation}
be the rescaled component sizes of the Bohman-Frieze process at time $t_c+\alpha\beta^{2/3} \frac{\lambda}{n^{1/3}} $. Then
\[\bar {\bfC}_n^{\sss BF} \stackrel{d}{\longrightarrow} \bfX\]
as $n\to\infty$, where $\bfX$ is the eternal standard multiplicative coalescent, and  $\convd$ denotes weak convergence in the space $\DD((-\infty, \infty):\ldown)$. In particular for each fixed $\lambda \in \RRR$, the rescaled component sizes in the critical window satisfy ${\bar \bfC}_n^{\sss BF}(\lambda)\convd \Xi(\lambda)$ where $\Xi(\lambda)$ are the sizes of the excursions for the reflected inhomogeneous Brownian motion in \eqref{eqn:inh-ref-bm}. 
\end{Theorem}

\section{Discussion}
\label{sec:disc}

Let us now give a wide ranging discussion of the implications of the above result and proof techniques.
\renewcommand{\labelenumi}{(\alph{enumi})}
\begin{enumerate}
	\item {\bf Dynamic network models:} The formulation and study of dynamic networks models, such as the Bohman-Frieze process wherein one has both randomness and choice,  are relatively new.
	 Motivated by a question by Achlioptas, that whether  simple modifications of the Erd\H{o}s-R\'{e}nyi random graph model could allow one to avoid the giant component for a larger amount of time, such models    now fall loosely under the term {\it Achlioptas processes} (see e.g. \cite{bohman2001avoiding,bohman2004avoidance,spencer2007birth,krivelevich2010hamiltonicity}).  The models considered in the above papers use {\it bounded size} rules, wherein there is a fixed $K$ such that when making a choice between the two randomly selected edges, all component sizes greater than some $K$ are treated identically (e.g. $K=1$ for the model considered in this study). Much recent interest has been generated by recent models formulated in \cite{achlioptas2009explosive}, wherein interesting examples of {\it unbounded size} rules have been analyzed. Simulations of such models seem to suggest a new phenomenon called ``explosive percolation'', wherein the phase transition  from a subcritical  to a supercritical regime with a giant component appears much more ``abruptly'', in that the largest component seems to increase from a size smaller than $\sqrt{n}$ to a size larger than $n/2$ in a very small scaling window. Examples of such rules include the product rule wherein one chooses two edges at random and connects the edge that \emph{minimizes} the product of the component sizes on the 2 end points of the edge. See \cite{d2010local,chen2010explosive} for additional models of this type. Recently in \cite{riordan2011achlioptas} it was shown that for a number of such models, the phase transition is actually ``continuous", so for such models one expects similar behavior as what one sees in the Erd\H{o}s-R\'{e}nyi random graph model.  However,  understanding what happens at criticality and the behavior of the scaling window is a challenging math program. The techniques in this paper have the potential to be extended to the analysis of a number of such models, which we attempt to do in work in progress. 
	\item {\bf Recent results on the Bohman-Frieze model:} 
	After this work was submitted for publication, we came across the interesting preprint \cite{bf-spencer-perkins-kang} that was announced around the same time. The latter paper studies the Bohman-Frieze process in the cases when $t=t_c-\eps$ and $t=t_c+\eps$, for fixed $\eps>0$. In \cite{bf-spencer-perkins-kang} the largest and second largest components are studied, bounds on the sizes of these components,  as well as the number of surplus edges amongst all components in these regimes are derived. The techniques  used in \cite{bf-spencer-perkins-kang} for understanding the size of the largest component in the subcritical regime are very different from those used in the current paper. In particular, our work requires an understanding of the fine scale asymptotics of the entire vector of component sizes at criticality (properly rescaled). For identifying the scaling window and how the giant component emerges via the merger of small components through the scaling window, we need more refined results at the critical value, in particular we need to study the behavior as $\eps = \eps(n)\to 0$ reasonably quickly -- see Proposition \ref{prop:main} where precise asymptotic results for   $\eps= 1/n^{\gamma}$, $\gamma \in (1/6,1/5)$ are obtained. Conjecture 1 in  \cite{bf-spencer-perkins-kang} states that one expects an upper bound of $\epsilon^{-2}\log n$ on the
	largest component at the time instant $(t_c-\epsilon)$.  This should be compared with the $(\log n)^4/\eps^2$ upper bound for times $t_c-\epsilon$ established in Proposition \ref{prop:com-size} of the current work, not just for a fixed $\eps$ but for $\eps=\eps(n)\to 0$. In fact, the proposition establishes an estimate that
	is {\em uniform} over the time interval $(0, t_c - n^{-\gamma})$, $\gamma \in (0, 1/5)$. This proposition is at the heart of our analysis and its proof requires significant work and  new ideas and techniques for general random graph models with immigrating vertices and near critical multi-type branching processes with general state spaces, see Sections \ref{sec:related-models} and \ref{sec:largest-com}. 
	\item {\bf Multiplicative coalescent:} As shown in the current work and several other papers, see e.g \cite{aldous2000random,bhamidi-hofstad-van,bhamidi2009novel,nachmias-peres},
	multiplicative coalescent arises as the limit object in a large number of random graph models at criticality.  Thus asymptotics of large random graphs has intimate connections with 
	the general theory of coalescent (and fragmentation) processes that is currently a very active area in probability, see e.g. \cite{aldous1999deterministic,bertoin-frag-book,pitman-book}. 
	\item {\bf Starting from an arbitrary configuration:} 
	 In this work we have only considered the case where we start with the empty configuration $\bfzero_n$. One can imagine starting from a different configuration and then attempt to analyze the emergence of the giant component. Along the lines of \cite{aldous1998entrance},  under the assumption that starting configuration satisfies suitable regularity conditions,
	it should be possible to study the asymptotic behavior of this model in terms of the entrance boundary of the standard multiplicative coalescent. This will
	be pursued in future work.
	\item {\bf Proof techniques:} The study of critical random graphs have generated an enormous amount of interest in the probabilistic combinatorics community and a large number of techniques have been developed for the fine scale asymptotics of such models, ranging from generating function and counting arguments (see e.g. \cite{janson1994birth}, \cite{bollobas-rg-book}); branching processes exploring local neighborhoods of the graph (see e.g \cite{janson-luczak-bb},and the references therein,\cite{molloy-reed-crit}); multi-type branching processes (see e.g. \cite{bollobas-riordan-janson} and the references therein); differential equation based techniques (see e.g. \cite{spencer2007birth}, \cite{bohman2001avoiding} and for a comprehensive survey \cite{wormald1999differential}); and the the breadth first walk approach coupled with the martingale central limit theorem (see e.g.\cite{karp1990transitive},\cite{martin-lof} for one of the first few studies using this technique; \cite{aldous1997brownian} where this was used to explore the fine scale structure of the classical Erd\H{o}s-R\'{e}nyi random graph at criticality and connections with the multiplicative coalescent; see \cite{nachmias-peres}, \cite{bhamidi-hofstad-van} for further results using this technique). Our methods are inspired by
	 \cite{aldous2000random} which uses estimates on the size of largest component in the subcritical window to analyze the asymptotic behavior of the 
the sum of squares and cubes of component sizes near criticality.  Estimates on the latter in turn allow the verification of the sufficient conditions given in  \cite{aldous1997brownian}
for convergence to the standard multiplicative coalescent. Although our general outline of the proof is similar to  \cite{aldous2000random} as described above, it turns out that for general Achlioptas processes, completing this program is a rather challenging problem and for the BF model treated in the current work we need to develop quite a bit of machinery, which is done in Sections \ref{sec:related-models} and \ref{sec:largest-com}, in order to obtain the required estimates.
	\item {\bf Other structural properties of components:} In this study we focus on the sizes of the  components in the critical scaling window. Developing probabilistic methodology for understanding the actual structure of the large components is also of great interest, see for example \cite{braunstein2003optimal} for one of the first studies at a non-rigorous level, exploring the connection between the structure of these components at criticality and the internal structure of the minimal spanning tree in various random graph models (the ``strong disorder'' regime in statistical physics).  In particular, it was conjectured in \cite{braunstein2003optimal} that the diameter of the minimal spanning trees in various random graph models scales like $n^{1/3}$. Rigorous studies have now been carried out for the Erd\H{o}s-R\'{e}nyi random graph model and fascinating connections have been discovered between the critical random graphs and the famous continuum random tree of Aldous (see \cite{addario2009critical}, \cite{addario2009continuum}). Proving such structural convergence of the entire components in dynamic network settings such as Achlioptas process at criticality,  to random fractals such as continuum random trees,  would be of great interest in a number of different fields.  
\end{enumerate}

\section{Proof idea}
\label{sec:proof-idea}
Let us now give an idea of the proof.  We begin by showing in Proposition \ref{prop:main} below that, just before the critical window, the configuration of the components satisfies some important regularity properties. 
This proposition will be used in Section \ref{sec:proof-main}  in order to apply a result of  \cite{aldous1997brownian} that gives sufficient conditions for convergence to the
multiplicative coalescent.
\begin{Proposition}
\label{prop:main}
Let $\gamma \in (1/6,1/5)$ and define $t_n = t_c- n^{-\gamma}$. Then we have
\begin{align}
\frac{n^2 \calS_3(t_n)}{\calS_2^3(t_n)} &\convp \beta \label{eqn:s3-s2}\\
\frac{n^{4/3}}{\calS_2(t_n)} - \frac{n^{-\gamma+1/3}}{\alpha} &\convp 0 \label{eqn:s2-tn-n-alpha}\\
\frac{n^{2/3}\calC_n^{\sss (1)}(t_n)}{\calS_2(t_n)} &\convp 0. \label{eqn:max-s2}
\end{align}

\end{Proposition}
Now note that $t_n$ can be written as
\[t_n = t_c+ \beta^{2/3}\alpha\frac{\lambda_n}{n^{1/3}}\]
where
\[\lambda_n = -\frac{n^{-\gamma+1/3}}{\alpha\beta^{2/3}} \to -\infty\]
as $n\to\infty$. The above proposition implies that the configuration of rescaled component sizes, for large $n$ at time ``$-\infty$'', satisfy the regularity conditions for the standard multiplicative coalescent (see Proposition 4 in \cite{aldous1997brownian}). 

Once the above has been proved, the second step is to show that through the critical window, the component sizes merge as in the multiplicative coalescent, at rate proportional to the product of the rescaled component sizes. This together with arguments similar to \cite{aldous2000random} will complete  the proof of the main result. \\

Let us now outline the framework of the proof:

\begin{itemize}
	\item The bound on the largest component $\CC_n^{\sss (1)}(t)$ when $t \uparrow t_c$ (Proposition \ref{prop:com-size}) plays a crucial role in proving  the statements in Proposition \ref{prop:main}. In order to achieve this, we introduce a series of related models from Section \ref{sec:model-equiv} through Section \ref{sec:model-irg}.
	\item Section \ref{sec:largest-com} uses these models to prove asymptotically tight bounds on the size of the largest component through the  subcritical window. The main goal of this Section is to prove Proposition \ref{prop:com-size}. 
	\item Section \ref{sec:analysis-s2s3} uses the bounds on the largest component from Proposition \ref{prop:com-size} to analyze the sum of squares and cubes of component sizes near the critical window. As one can imagine when working so close to the critical window, one needs rather careful estimates.
	Our arguments are based on a precise analysis of quadratic variation processes for certain martingales associated with the BF model. This section completes the proof of Proposition \ref{prop:main}. 
	\item Finally, in Section \ref{sec:proof-main}  we use Proposition \ref{prop:main} and a coupling with the standard multiplicative coalescent, in a manner similar to
	\cite{aldous2000random}, to prove the main result.    
\end{itemize}
Without further ado let us now start with the proofs.

\section{An estimate on the largest component }
\label{sec:related-models}

The following estimate on  the largest component is the key ingredient in our analysis.    Recall that $t_c$ denotes the critical time for the BF process.

\begin{Proposition}
\label{prop:com-size}
Let $\gamma \in (0, 1/5)$ and let  $I_n(t) \equiv \CC_n^{\sss (1)}(t)$ be the largest component of $\BF_n(t)$.  Then, for some $B\equiv B(\gamma) \in (0, \infty)$,
$$\prob \{ I_n(t) \le m(n,t), \forall t < t_c-n^{-\gamma}\} \to 1, \mbox{ when } n \to \infty$$
where
\be \label{ins1711} m(n,t)= B \frac{(\log n)^4}{(t_c-t)^2} .\ee

\end{Proposition}

The proof of Proposition \ref{prop:com-size} will be  completed in Section \ref{sec:proof-prop-reduced}.
In the current section we will give constructions of some auxiliary random graph processes that are key to our analysis. Although not pursued here, we believe that analogous constructions will be key ingredients
in treatment of more general random graph models as well.
The section is organized as follows.
 \begin{itemize}
 \item In Section \ref{sec:notation-con} we  give the basic notation and mathematical conventions used in this paper.\\
 \item In Section \ref{sec:model-equiv} we will carry out a preliminary analysis of the BF process and identify three deterministic maps $a_0, b_0, c_0$ from $[0, \infty)$ to
 $[0,1]$  that play a fundamental role in our analysis.\\
 \item Guided by these deterministic maps, in Section \ref{sec:model-rgiva} we will define a random graph process with immigrating vertices and attachments (RGIVA) which is simpler to analyze than, and is suitably `close' to, the Bohman-Frieze process.  A precise estimate on the approximation error introduced through this model
 is obtained in Section \ref{sec:proof-prop-reduced}.\\
 \item In Section \ref{sec:model-irg} we will introduce an inhomogeneous random graph (IRG) model associated with a given RGIVA model such that
 the two have identical component volumes at all times. This allows for certain functional analytic techniques to be used in estimating the maximal component size.
 We will also make an additional approximation to the IRG model which will facilitate the analysis. \\
 \item In Section \ref{sec:summary-model} we  summarize   connections between  the various models introduced above.\\
 \end{itemize}

\subsection{Notation }
\label{sec:notation-con}

\subsubsection{ Graphs  and random graphs}
A graph $\bfG=\{\VV, \EE\}$ consists of a vertex set $\VV$ and an edge set $\EE$, where $\VV$ is a subset of some type space $\XX$ and $\EE$ is a subset of all possible edges $\{ \{v_1,v_2\}:v_1 \ne v_2 \in \VV\}$.  An example of a type space is $[n]=\{1,2,...,n\}$. Frequently we will assume $\XX$ to have additional structure, for example to be a measure space $(\XX,\TT,\mu)$. When $\VV$ is a finite set, we write $|\VV|$ for its cardinality.

 $\bfG$ is called \textbf{null graph} if $\VV=\emptyset$, and we write $\bfG=\emptyset$. $\bfG$ is called an \textbf{ empty graph} if $|\VV|=n$ and $\EE=\emptyset$, and we write $\bfG = {\bf0}_n$.

Given two graphs, $\bfG_i=\{\VV_i,\EE_i\}$ for $i=1,2$, $\bfG_1$ is said to be a \textbf{subgraph} of $\bfG_2$ if and only if $\VV_1 \subset \VV_2$ and
$\EE_1 \subset \EE_2$ and we denote this as $\bfG_1 \le \bfG_2$ (or equivalently $\bfG_2 \ge \bfG_1$). We write $\bfG_1 =\bfG_2$ if $\bfG_1 \le \bfG_2$ and $\bfG_1 \ge \bfG_2$. 

A connected component $\CC=\{\VV_0,\EE_0\}$ of a graph $\bfG=\{\VV, \EE\}$ is a subgraph which is connected (i.e. there is a path between any two vertices
in $\CC$).  The number of vertices in $\CC$ will be called the size of the component and frequently we will denote the size and the component by the
same symbol.


Let $\GG$ be the set of all  possible graphs $(\VV, \EE)$ on a given type space $\XX$.  When $\VV$ is countable, we will
consider $\GG$ to be endowed with the discrete topology and the corresponding Borel sigma field and refer to a random element of $\GG$ as a random graph. All random graphs in this work are given on a fixed probability space $(\Omega, \mathcal{F}, \mathbb{P})$ which will usually be
suppressed in our proofs. 

\subsubsection{Probability and  analysis}
All the unspecified limits are taken as $n \to +\infty$.
Given a sequence of events $\{E_n\}_{n\ge 1}$, we say $E_n$ (or $E$) occurs with high probability (whp) if $\prob\{E_n\} \to 1$. 
For functions $f, g: \mathbb{N} \to \mathbb{R}$, we write
 $g=O(f)$ if for some $C \in (0, \infty)$, $\limsup g(n)/f(n) < C$ and $g=\Theta(f)$ if $g=O(f)$ and $f=O(g)$.  Given two sequences of random variables $\{\xi_n\}$ and $\{\zeta_n\}$, we say $\xi_n=O(\zeta_n)$ whp if there is a $C \in (0, \infty) $ such that $\xi_n < C \zeta_n$ whp, and write $\xi_n=\Theta(\zeta_n)$ whp if there exist  $0 < C_1 \le C_2 < \infty$ such that $C_1 \zeta_n<\xi_n<C_2 \zeta_n$ whp. Occasionally, when clear from the context, we suppress `whp' in the statements. 

 We also use the following little $o$ notation: For a sequence of real numbers $g(n)$, we write $g=o(f)$ if $\limsup|g(n)/f(n)|=0$. For a sequence of random variables $\xi_n$, we write ``$\xi_n=o_p(f)$'' if $\xi_n/f(n)$ converges to $0$ in probability.

For a real measurable function $\psi$ on a measure space $(\XX,\TT,\mu)$, the norms $\|\psi\|_2$ and $\|\psi\|_\infty$ are defined in the usual way. We use $\convp$ and $\convd$ to denote the convergence in probability and in distribution respectively.

We use $=_d$ to denote the equality of random elements in distribution. Suppose that $(S, \cals)$ is a measurable space and we are given a partial
ordering on $S$. Given two $S$ valued random variables $\xi_1, \xi_2$, we say a pair of $S$ valued random variables $\xi_1^*, \xi_2^*$ given on
a common probability space define a coupling of $(\xi_1, \xi_2)$ if $\xi_i =_d\xi_i^*$, $i=1,2$.
We say 
the $S$ valued random variable $\xi_1$ \textbf{ stochastically dominates} $\xi_2$, and write $\xi_1 \ge_d \xi_2$ if there exists a coupling between the two random variables, say $\xi_1^*$ and $\xi_2^*$, such that $\xi_1^* \ge \xi_2^*$ a.s. 

 For two sequences of $S$ valued random elements $\xi_n$ and $\tilde \xi_n$, we say ``$\xi_n \le_d \tilde \xi_n$ whp.'' if there exist a coupling between $\xi_n$ and $\tilde \xi_n$ for each $n$ (denote as $\xi_n^*$ and $\tilde \xi_n^*$) such that $\xi_n^* \le \tilde \xi_n^*$ whp.

Two examples of $S$ that  are relevant to this work are $\DD([0,T]: \mathbb{R})$ and $\DD([0,T]: \GG)$ with the natural associated partial ordering.

\subsubsection{Other conventions}

We always use $n, m, k, i, j$ to denote non-negative integers unless specified otherwise. We use $s, t, T$ to denote the time parameter for continuous
time (stochastic) processes. The scaling parameter is denoted by $n$. Throughout this work $T = 2t_c$ which is a convenient upper bound for the time parameters of interest.


We use $d_1, d_2, ...$ for constants whose specific value are not important. Some of them may appear several times and the values might not be the same. We use $C_1, C_2, ...$ for constants that appear in the statement of theorems.

\subsection{A preliminary analysis of Bohman-Frieze process}
\label{sec:model-equiv}
Recall that $\BF_n(t)$ denotes the BF process at time $t$ and note that $\BF_n$ defines a stochastic process with sample paths in
$\DD([0, T]:\clg)$. Also recall that
$\CC_n^{\sss (i)}(t)$ denotes the size of the $i^{th}$ largest component in $\BF_n(t)$,  $\bfC_n(t)=(\CC_n^{\sss (i)}(t) : i \ge 1)$ is the vector of component sizes  and
$X_n(t)$ denotes the number of singletons in  $\BF_n(t)$.  We let $\FF_t \equiv \FF^n_t = \sigma \{ \BF_n(s), s \le t\}$ and refer to it as the natural filtration for the BF process.

At any fixed time $t>0$, let $\com(t)$ denote the collection of all non-singleton components 
\[\com(t)= \set{\CC_n^{\sss (i)}(t): |\CC_n^{\sss (i)}(t)|\geq 2 }.\] 
Recall that $\bar{x}(t)=X_n(t)/n$.
 We will now do an informal calculation of the rate at which an  edge $e=\{v_1,v_2\}$ is added to the graph $\BF(t)$.  There are  three different ways an edge can be added: (i) both $v_1$ and $v_2$ are singletons, (ii) only one of them is a singleton, (iii) neither of them is a singleton. \\
 
\textbf{Analysis of the three types of events:}   \\

(i) Both $v_1$ and $v_2$ are singletons.  We will refer to such a component that is formed by connecting  two singletons as a \textbf{doubleton}.  This will happen at rate
\begin{equation}
\frac{2}{n^3} \left[{X_n(t) \choose 2} {n \choose 2} +\left({n \choose 2}-{X_n(t) \choose 2}\right){X_n(t) \choose 2}\right] \stackrel{\sss def}{=} n \cdot a_n^*(\bar{x}(t)). \label{eqn:ins609}
\end{equation}
The first product in the squared brackets is the count of all possible $\bfe=(e_1,e_2) \in \EE^2$ such that $e_1$ joins up two singletons and thus will be added to the graph, while the second product is the count of all $\bfe=(e_1,e_2) \in \EE^2$ such that the first edge $e_1$ does {\bf not} connect two singletons while $e_2$ connects two singletons and will be added. \\
Define $a_0: [0,1] \to [0,1]$ as 
\be
a_0(y) = 2\left(\frac{y^2}{2} \cdot \frac{1}{2}+\left(\frac{1}{2}-\frac{y^2}{2}\right) \frac{y^2}{2}\right)=\frac{1}{2}
	(y^2+(1-y^2)y^2).
	\label{eqn:a-def}
\ee
It is easy to check that
\be
\label{eqn:a-astar}
a_n^*(\bar{x}(t)) = a_0(\bar{x}(t))+r_a(t), \mbox{ where, } \sup_{t}|r_a(t)| \le 5/n.
	\ee
	Recall that  $x(t)$ is the solution of the differential equation  \eqref{eqn:xss-def}. 	
To simplify notation we will write $a_n^*(\barx (t))=a^*(\barx)=a^*(t)=a^*$ and $a_0(t)=a_0(x(t))$ exchangeably. 
Similar conventions will be followed for the functions $c_n^*, c_0$ and $b_n^*, b_0$ that will be introduced below. 
We shall later show that  $\sup_{t\le T}|\barx_n(t) - x(t)| \to 0$ in probability (see Lemma \ref{lemma:error-diff-eqn} of this paper,  also see \cite{spencer2007birth}). This in particular implies that $\sup_{t\le T}|a^*_n(t) - a_0(t)| \to 0$ 
in probability. 
 \\

(ii) Only one of them is a singleton: This will happen if and only if $e_1$ does not connect two singletons while $e_2$ connects a singleton and a non-singleton, thus at the rate
\begin{equation}
\frac{2}{n^3} \left({n \choose 2}-{X_n(t) \choose 2}\right) (n-X_n(t)) X_n(t). \label{eqn:ins633}
\end{equation}
We are also interested in the rate that a given non-singleton vertex (say, $v_0$) is connected to any singleton, which is
\begin{equation}
\frac{2}{n^3} \left({n \choose 2}-{X_n(t) \choose 2}\right) X_n(t) \stackrel{\sss def}{=} c^*_n(\barx(t)).
\end{equation}
Thus at time $t$  a singleton will be added to $\com(t)$ during the small time interval $(t, t+dt]$, by attaching to a given vertex $v_0 \in \com(t)$, with the rate $c^*(t)$.\\Define $c_0: [0,1] \to [0,1]$ as 
\begin{equation}
	c_0(y)=(1-y^2)y, \; y  \in [0,1].
	\label{eqn:c-def}
\end{equation}
Then
\begin{equation}
	c^*(\barx(t))=c_0(\barx(t))+r_c(t) \qquad \mbox{ and } \sup_{t}|r_c(t)| \le 2/n.
	\label{eqn:c-def1}
\end{equation}

(iii) Neither of them is a singleton:  This will happen at the rate
\begin{equation}
\frac{2}{n^3} \left({n \choose 2}-{X_n(t) \choose 2}\right) {n-X_n(t) \choose 2}.
\end{equation}
Also, the event that two fixed non-singleton vertices are connected has the rate
\begin{equation}
\frac{2}{n^3} \left({n \choose 2}-{X_n(t) \choose 2}\right) \stackrel{\sss def}{=} \frac{1}{n} b^*_n(\barx(t)). \label{ins1516}
\end{equation}
Let $b_0: [0,1] \to [0, 1]$ be defined as
\begin{equation}
	b_0(y)=1-y^2 , \; y \in [0,1].
	\label{eqn:b-def}
\end{equation}
Then 
\begin{equation}
	b^*(\barx(t))=b_0(\barx(t))+r_b(t) \qquad \mbox{ and } \sup_{t}|r_b(t)| \le  2/n.
	\label{eqn:b-def1}
\end{equation}



Note that for the study of the largest component one may restrict attention to the subgraph $\com(t)$.  The evolution of this subgraph is described in terms of  
stochastic processes $a^*(\bar x(t)), b^*(\bar x(t))$ and $c^*(\bar x(t))$.
In the next subsection, we will introduce a random graph process that is ``close" to $\com(t)$ but easier to analyze. Intuitively, we replace $a^*(t), b^*(t), c^*(t)$ with  deterministic functions $a(t),b(t),c(t)$ which are close to $a_0(t),b_0(t),c_0(t)$ (and  thus, from Lemma \ref{lemma:error-diff-eqn}, whp close to $a^*(\bar x(t)), b^*(\bar x(t)), c^*(\bar x(t))$)
and construct a random graph with similar dynamics as $\com(t)$. \\

\subsection{The random graph process with immigrating vertices and attachment}
\label{sec:model-rgiva}
In this subsection, we introduce a random graph process with immigrating vertices and attachment (RGIVA). This construction is inspired by \cite{aldous2000random} where
a random graph with immigrating vertices (RGIV) is constructed --   we generalize this construction by including attachments.
RGIVA process will be governed by three continuous maps $a,b,c$ from $[0, T] \to [0, 1]$ (referred to as {\bf rate functions}) and the graph at time $t$ will be denoted by
$\IA_n(t)=\IA_n(a,b,c)_t$.   When $(a,b,c)$ is sufficiently close to $(a_0, b_0, c_0)$ , the RGIVA model well approximates the BF model in a sense that will be made
precise in Section \ref{sec:proof-prop-reduced}.\\

\textbf{The RGIVA process} $\IA_n(t)=\IA_n(a,b,c)_t$. Given the rate functions $a,b,c$, define $\IA_n(t)$ as follows:\\
(a) $\IA_n(0)=\emptyset$, the null graph;\\
(b) For $t \in [0,T)$, conditioned on $\IA_n(t)$, during the small time interval $(t,t+dt]$,
\begin{itemize}
\item (immigration) a doubleton (consisting of two vertices and a joining edge) will be born at rate $n \cdot a(t)$,
\item (attachment) for any given vertex $v_0$ in $\IA_n(t)$, a new vertex will be created and connected to $v_0$ at rate $c(t)$,
\item (edge) for any given pair of vertices $v_1,v_2$ in $\IA_n(t)$, an edge will be added between them at rate $\frac{1}{n} \cdot b(t)$.
\end{itemize}
The events listed above occur independently of each other.\\

In the special case where  $a(t)\equiv b(t) \equiv 1$, $c(t) \equiv 0$, and doubletons are replaced by singletons, the above model reduces to  the  RGIV model of \cite{aldous2000random}. 
We note that the above construction  closely follows our analysis of three types of events in Section \ref{sec:model-equiv}, replacing stochastic processes
$a^*(\bar x_n(t)), b^*(\bar x_n(t)), c^*(\bar x_n(t))$ with deterministic maps $a(t), b(t), c(t)$.

The following lemma establishes a connection between the Bohman-Frieze process and the RGIVA process.  Recall the partial order on the
space
$\DD([0,T]: \bfG)$.
\begin{Lemma}
\label{lemma:couple-bfia} Let $(a_L,b_L,c_L)$ and $(a_U,b_U,c_U)$ be rate functions. Further, let $U \equiv U_n$ be the event that $\{ a^*(t)\le a_U(t), b^*(t)\le b_U(t), c^*(t)\le c_U(t) \mbox{ for all } t \in [0,T]\}$ and 
$L \equiv L_n$ be the event that $\{ a^*(t) \ge a_L(t), b^*(t) \ge b_L(t), c^*(t) \ge c_L(t) \mbox{ for all } t \in [0,T]\}$.
Define for $t \in [0,T]$
\[
\com_n^U(t) = \left\{ 
  	\begin{array}{l l}
    	\emptyset & \quad \text{on $U^C$}\\
    	\com_n(t) & \quad \text{on $U$}\\
  \end{array} \right. ;\;\;  \com_n^L(t) = \left\{ 
  	\begin{array}{l l}
    	\IA_n(a_L,b_L,c_L)_T & \quad \text{on $L^C$}\\
    	\com_n(t) & \quad \text{on $L$}\\
  \end{array} \right.\]
Then\\
(i)Upper bound: $\com_n^U \le_d \IA_n^U \equiv \IA_n(a_U,b_U,c_U)$.\\
(ii)Lower bound: $\com_n^L \ge_d \IA_n^L \equiv \IA_n(a_L,b_L,c_L)$.
\end{Lemma}
\textbf{Proof:} We only argue the upper bound.  The lower bound is proved similarly.  
Construct $\IA_n^U(t)$ iteratively on $[0, T]$ as described in the definition, and construct $\com_n^U(t)$ simultaneously by rejecting the proposed change on the graph with probabilities $(1-a^*/a_U)^+$, $(1-b^*/b_U)^+$ and $(1-c^*/c_U)^+$ according to the three types of the events. Let $\tau= \inf\{0\le t \le T:  a^*(t)>a_U(t) \mbox{ or } b^*(t)>b_U(t) \mbox{ or } c^*(t)>c_U(t) \}$ and set $\com_n^U(t)$ to be the null graph whenever $t \ge \tau$. This construction defines a coupling of 
$\IA_n^U$ and $\com_n^U$ such that $\com_n^U \le \IA_n^U$ a.s.   The result follows. \qed\\


\subsection{An inhomogeneous random graph with a weight function}
\label{sec:model-irg}
In this section we introduce a inhomogeneous random graph (IRG) associated with $\IA_n(a,b,c)$ for given rate functions $a,b,c$.
For a general treatment of  IRG models we refer the reader to  \cite{bollobas-riordan-janson},  which our presentation largely follows. We  generalize the setting
of \cite{bollobas-riordan-janson} somewhat  by including a weight function and considering the volume of a component instead of the number of vertices of a component. We begin with 
a description and some basic definitions for  a general IRG model.\\

A \textbf{type space} is a measure space $(\XX, \TT, \mu)$ where $\XX$ is a complete separable metric space (i.e. a Polish space), $\TT$ is the Borel $\sigma$-field and $\mu$ is a finite measure. 

A \textbf{kernel} on the type space $(\XX, \TT, \mu)$ is a measurable function $\kappa :\XX\times \XX \to [0, \infty)$. The kernel $\kappa$ is said to be symmetric if $\kappa(\bfx, \bfy)=\kappa(\bfy,\bfx)$ for all $\bfx,\bfy \in \XX$.
 We will also use $x, y$ instead of $\bfx, \bfy$ for  elements in $\XX$ when there is no confusion between an $x \in \XX$ and the function $x(t)$  defined in \eqref{eqn:xss-def}.

A \textbf{weight function} $\phi$ is a measurable, non-negative function on $(\XX, \TT, \mu)$. 

A \textbf{basic structure} is a triplet $\{ (\XX, \TT, \mu), \kappa, \phi\}$, which consists of a type space, a kernel and a weight function.\\

\textbf{The IRG model:} Given a type space $(\XX, \TT, \mu)$, symmetric kernels $\{\kappa_n\}_{n\ge 1}$, and a weight function $\phi$, a random graph $\RG_n(\kappa_n)$ ( $\equiv \RG_n(\kappa_n, \mu) \equiv \RG_n(\kappa_n, \mu, \phi)$), for any integer $n>0$,
is constructed as follows:\\
(a) The vertex set $\VV$ are the points of a Poisson point process on $(\XX,\TT)$ with intensity $n\cdot \mu$.\\
(b) Given $\VV$, for any two vertices $x,y \in \VV$, place an edge between them with probability $\left(\frac{1}{n} \cdot \kappa_n(x,y) \right)\wedge 1$.\\

One can similarly define an IRG model associated with a basic structure $\{ (\XX, \TT, \mu), \kappa, \phi\}$, where $\kappa$ is a symmetric kernel, by letting $\kappa_n=\kappa$ for all $n$ in the above definition.

The weight function $\phi$ is used in defining the volume of a connected component in the above construction of a random graph. Given a component of $\RG_n(\kappa, \mu, \phi)$ whose vertex set is $\VV_0$, define $\sum_{x\in \VV_0}\phi(x)$ as the \textbf{volume} of the component. 

One can associate $\kappa$ with an integral opertor $\KK : L^2(\mu) \to L^2(\mu)$ defined as
\begin{equation}
	\KK f (x)=\int_{\XX} \kappa(x,y)f(y)\mu(dy)
	\label{eqn:calk-def}
\end{equation}
Denote by $\rho = \rho(\kappa)$ the operator norm of $\KK$.  Then $\rho=\rho(\kappa)=\|\KK\|=\sup_{\|f\|_2=1}\|\KK f\|_2$. \\

Given rate functions $a,b,c$, there is a natural basic structure and the corresponding IRG model associated with $\IA_n(a,b,c)$, which we now describe. \\

Fix $t \in [0,T]$. 
Then the following two stage construction describes an equivalent (in law) procedure for obtaining $\IA_n(a,b,c)_t$ :\\

\textbf{Stage I}:  Recall that transitions in $\IA_n(a,b,c)$ are caused by three types of events: {\em immigration}, {\em attachment} (to an existing vertex) and {\em edge formation} (between existing
vertices). Consider the random graph obtained by including all the immigration and attachment events until time $t$ but ignoring the edge formation events.  We call the components resulting from this construction as clusters. Note that each cluster consists of exactly one doubleton (which starts the formation of the cluster) and possibly other vertices obtained through later attachments. Note that doubletons immigrate  at rate $a(s)$ and supposing that a  doubleton is born at time $s$,  the size of the cluster at time $s \le u \le t$ denoted by $w(u)$ evolves according to a integer-valued time-inhomogeneous  jump Markov process starting at $w(s)=2$ and infinitesimal generator $\cla(u)$
given as  
\be \label{eqn:ws-def}\cla(u) f(r) = c(u) r \cdot \left(f(r+1) - f(r)\right), \; f: \NNN \to \RRR, s \le u \le t.\ee
We set $w(u) = 0$ for $0 \le u < s$ and  denote this cluster which starts at instant $s$ by $(s,w)$.

\textbf{Stage II:} Given a realization of  the random graph of Stage I, we add edges to the graph. Each pair of vertices will be connected during $(s,s+ds]$ with rate $\frac{1}{n}b(s)$. Thus the number of edges between two clusters $\bfx=(s, w), \bfy=(r, \tilde w)$ at time instant $t$ is a Poisson random variable with mean 
$\frac{1}{n} \int_0^t w(u)\tilde w(u)b(u)du$. Consequently, 
\begin{align}
	\prob\{ \bfx \mbox{ and } \bfy \mbox{ is connected } | \mbox{ Stage I} \} &=1-\exp \{ -\frac{1}{n} \int_0^t w(u)\tilde w(u)b(u)du \}\\
	&\le \frac{1}{n} \int_0^t w(u)\tilde w(u)b(u)du.
	\label{eqn:approx-ker}
\end{align}
\\

It is easy to see that the graph resulting from this two stage construction has the same distribution as $\IA_n(a,b,c)_t$.  

We now introduce an IRG model associated with the above construction in which each cluster is treated as a single point in a suitable type space and the size of the cluster
is recorded using an appropriate weight function.  Let $\XX = [0,T] \times \WW$, where $\WW=\DD([0,T]: \mathbb{N})$ is the Skorohod $D$-space with the usual
Skorohod topology.  Denote by $\TT$ the Borel sigma field on $[0,T] \times \WW$. For future use, we will refer to this particular choice of
type space $(\XX, \TT)$  as the {\em cluster space}. For a fixed  time $t\geq 0$, consider a weight function defined as 
\be
\label{eqn:phi-def}\phi_t(\bfx) = w(t), \;\; \bfx = (s,w) \in [0,T] \times \WW.
\ee
Then this weight function associates with each `cluster' $\bfx$ its size at time $t$.  We now describe the finite measure $\mu$ that governs the intensity
of the Poisson point process $\PP_t(a,b,c)$ of clusters (regarded as points in $\XX$).  Denote by $\nu_s$ the unique probability measure on the space $\WW$ under which, a.s.,  $w(u) = 0$ for all $u < s$,
$w(s) = 2$ and $w(u), u \in [s, T]$ has the probability law of the time inhomogeneous Markov process with generator $\{\cla(u), s \le u \le T\}$ defined in \eqref{eqn:ws-def}.
Let $\mu$ be a finite measure on $\XX$ defined as $\mu(ds dw) = \nu_s(dw) a(s) ds$, namely, for a non-negative real measurable function $f$ on $\XX$ 
$$\int_{\XX} f(\bfx) d\mu(\bfx)  =  \int_0^T a(s)\left (\int_{\WW} f(s,w) d\nu_s(w) \right) ds.$$
We also define for each $t \in [0, T]$, a finite measure $\mu_t$ on $\XX$ by the relation $\mu_t(A) = \mu(A \cap ([0,t] \times \WW))$.  Then for $f$ as above,
\be \label{eqn:mu-def} \int_{\XX} f(\bfx) d\mu_t(\bfx)  =  \int_0^t a(s)\left (\int_{\WW} f(s,w) d\nu_s(w) \right) ds.\ee
The measure $\mu_t$ will be  the intensity of the Poisson point process on $\XX$ which will be used in our construction of the IRG model associated with $\IA_n(a,b,c)_t$.
Now we describe the kernel that will govern the edge formation amongst the points. Define
\be
	\kappa_{n,t}(\bfx, \bfy) =\kappa_{n,t}((s,w),(r,\tilde w))=n\left(1-\exp \{ -\frac{1}{n} \int_0^t w(u)\tilde w(u)b(u)du\}\right).
	\label{eqn:kappan-def}\ee
We will also use the following modification of the kernel $\kappa_{n,t}$.
\be	
	\kappa_t(\bfx, \bfy) =\kappa_t((s,w),(r,\tilde w))=\int_0^t w(u)\tilde w(u)b(u)du.
	\label{eqn:kappa-def}
\ee
With the above definitions we can now define IRG models $\RG_n(\kappa_{n,t}, \mu_t, \phi_t)$ and $\RG_n(\kappa_{t}, \mu_t, \phi_t)$ associated with  the type space
$ (\XX, \TT, \mu_t)$.

  Denote the size of the largest component [resp. the component containing the first immigrating doubleton]  in $\IA_n(a,b,c)_t$ by $\clc^{(1)}(a,b,c)_t$ [resp. $\clc^{(0)}(a,b,c)_t$].  Also, denote the volume of the largest component [resp. the component containing the first cluster]
in $\RG_n(\kappa_t, \mu_t, \phi_t)$ by $\clc^{(1)}(\kappa_t, \mu_t, \phi_t)$ [resp. $\clc^{(0)}(\kappa_t, \mu_t, \phi_t)$].  Define
$\clc^{(1)}(\kappa_{n,t}, \mu_t, \phi_t)$, $\clc^{(0)}(\kappa_{n,t}, \mu_t, \phi_t)$ in a similar fashion. The following is an immediate consequence of the above construction.
\begin{Lemma}
	\label{lemma:rgiva-irg-def}
	 We have \[(\clc^{(1)}(a,b,c)_t, \clc^{(0)}(a,b,c)_t) =_d (\clc^{(1)}(\kappa_{n,t}, \mu_t, \phi_t), \clc^{(0)}(\kappa_{n,t}, \mu_t, \phi_t))\]
	and
	$$\clc^{(1)}(\kappa_{n,t}, \mu_t, \phi_t) \le_d \clc^{(1)}(\kappa_{t}, \mu_t, \phi_t), \; \clc^{(0)}(\kappa_{n,t}, \mu_t, \phi_t) \le_d \clc^{(0)}(\kappa_{t}, \mu_t, \phi_t).$$
		%
\end{Lemma}
For future use we will write $\RG_n(\kappa_t, \mu_t, \phi_t) \equiv \RG_{n,t}(a,b,c)$.  

\subsection{A summary of the models}
\label{sec:summary-model}
As noted earlier, the key step in the proof of Proposition \ref{prop:main} is a good estimate on the size of the largest component in the Bohman-Frieze process $\BF_n(t)$
as in Proposition \ref{prop:com-size}.  For this  we have introduced a series of approximating models. We summarize the relationship between these models below.\\
\begin{itemize}
	\item We can decompose the Bohman-Frieze process as $\BF_n= \com_n\cup X_n$, namely the non-singleton components and singleton components at any time $t$. 
	\item We shall show that $\com_n \approx \IA_n(a_0,b_0,c_0)$, where $a_0, b_0, c_0$ are defined in \eqref{eqn:a-def}, \eqref{eqn:b-def}, \eqref{eqn:c-def}. More precisely we shall show that as $n\to\infty$, for any fixed $\delta>0$, we have, whp.
$$\IA_n((a_0-\delta)^+,(b_0-\delta)^+,(c_0-\delta)^+) \le_d \com_n \le_d \IA_n((a_0+\delta)\wedge 1,(b_0+\delta)\wedge 1,(c_0+\delta)\wedge 1). $$
This is a consequence of Lemma \ref{lemma:couple-bfia}.
	\item Given rate functions $(a,b,c)$, for all $t \in [0, T]$,	 
	$$\clc^{(i)}(a,b,c)_t =_d \clc^{(i)}(\kappa_{n,t}, \mu_t, \phi_t) \le_d \clc^{(i)}(\kappa_t, \mu_t, \phi_t), \; i = 0,1.$$  
	Here $\kappa_{n,t}, \kappa_t, \mu_t, \phi_t$ and $a,b,c$ are related through \eqref{eqn:kappan-def}, \eqref{eqn:kappa-def}, \eqref{eqn:mu-def} (see also \eqref{eqn:ws-def}), \eqref{eqn:phi-def}, respectively.
\end{itemize}

\section{Analysis of the largest component at sub-criticality: Proof of Proposition \ref{prop:com-size}} 
\label{sec:largest-com}

This section proves Proposition \ref{prop:com-size}.  The section is organized as follows:
\begin{itemize}
\item In Section \ref{sec:largest-first} we reduce the problem to proving Proposition \ref{prop:reduced}.  We give the proof of Proposition \ref{prop:com-size} using this
result.  Rest of Section \ref{sec:largest-com} is devoted to the proof of Proposition \ref{prop:reduced}.\\
\item In preparation for this proof, in Section \ref{sec:approx-error} we present some key lemmas that allow us to estimate the errors between various models summarized in Section \ref{sec:summary-model}.  Proofs of some lemmas (Lemmas \ref{theo:irg-first-com}, \ref{theo:error-cont-norm} and \ref{theo:error-cont-int}) are deferred to  later sections.  \\
\item Using these lemmas, in Section \ref{sec:proof-prop-reduced} we prove the key proposition, Proposition \ref{prop:reduced}. The rest of Section \ref{sec:largest-com}
proves the supporting Lemmas \ref{theo:irg-first-com}, \ref{theo:error-cont-norm} and \ref{theo:error-cont-int}.\\
\item In Section \ref{sec:model-branching} we introduce a branching process related to the IRG model, and prove Lemma \ref{theo:irg-first-com}.  A key step in the proof is
Lemma \ref{theo:branching} whose proof is left to Section \ref{sec:proof-branching}.\\
\item Section \ref{sec:analysis-ker} analyzes the kernel $\kappa_t$ associated with $\RG_{n,t}(a,b,c)$ and proves Lemma \ref{theo:error-cont-norm} .\\
\item Finally, in Section \ref{sec:measure-theoretic} we give the proof of Lemma \ref{theo:error-cont-int}.
\end{itemize}

%
%
%




\subsection{ From the largest component to the first component}
\label{sec:largest-first}
In this section we will reduce the problem of proving the estimate on the largest component in Proposition \ref{prop:com-size} to an estimate on the 
{\em first} component as in Proposition \ref{prop:reduced}. This  reduction, although somewhat different, is inspired by a similar idea used in \cite{aldous2000random}. \\
Recall that $\CC_n^{\sss (1)}(t) \equiv I_n(t)$ denotes the largest component in $\BF_n(t)$.
Let $\CC_n^s(t)$, $ 0 \le s \le t$, denote the component whose first doubleton is born at time $s$ in $\BF_n(t)$. In particular 
$\CC_n^s(t) = \emptyset$ if there is no doubleton born at time $s$. Without loss of generality, we assume that the first doubleton is born at time $0$.   Then $\CC_n^0(t)$ denotes the component of the \emph{first doubleton} at time $t$ of the BF process. 
The following lemma estimates the size of the largest component $I_n(t)$ in terms of the  size of the first component.
\begin{Lemma} \label{lem2105} For any $n \in \NNN$, $t_0 \in [0,T]$ and deterministic function $\alpha: [0,T] \to   [0, \infty)$
$$\prob\{ I_n(t) > \alpha(t), \mbox{ for some } t < t_0\} \le n T\prob\{ \CC_n^0(t) > \alpha(t), \mbox{ for some } t < t_0 \}.$$
\end{Lemma}
\textbf{Proof:} 
Let $\{\BF_n^{(i)}(t), t \ge 0\}_{i \in \NNN_0}$ be an i.i.d. family of $\{\BF_n(t), t \ge 0\}$ processes on the same vertex set $[n]$.
Let $N$ be a rate $n$ Poisson process independent of the above collection.  Denote by $\{\tau_i, i \in \NNN\}$ the jump times of the Poisson process.  Set $\tau_0=0$.
Denote the first component of $\BF_n^{(i)}$ at time $t$ by $\JJ_n^{(i)}(t)$.  Consider the random graph
$$\bfG_n^t = \cup_{i\in \NNN_0: \tau_i \le t} \JJ_n^{(i)}(t)$$
and let $I_n^{\bfG}(t)$ denote the size of the largest component in $\bfG_n^t$.  Then since $a_n^*(t) \le 1$ for all $t$, $I_n \le_d I_n^{\bfG}$.
Thus
\beq
\prob\{ I_n(t) > \alpha(t), \mbox{ for some } t < t_0\} &\le & \prob\{ I_n^{\bfG}(t) > \alpha(t), \mbox{ for some } t < t_0\}\\
&= & \sum_{k \in \NNN_0} \prob\{ I_n^{\bfG}(t) > \alpha(t), \mbox{ for some } t < t_0, N(T) = k\}\\
& \le &   \sum_{k \in \NNN_0} \prob\{ \JJ_n^{(i)}(t) > \alpha(t), \mbox{ for some } t < t_0, \mbox{ for some } i \le k\} \prob\{N(T) = k\}\\
& \le &   \sum_{k \in \NNN_0} k \prob\{ \CC_n^0(t) > \alpha(t), \mbox{ for some } t < t_0\} \prob\{N(T) = k\}.
\eeq
The result follows. \qed \\

Next, in the following lemma, we reduce an estimate on the probability of the event $\{ \CC_n^0(t) > \alpha(t), \mbox{ for some } t < t_0 \}$ to an estimate on 
 $\sup_{t \in [0, t_0]}\alpha(t)\prob\{ \CC_n^0(t) > \alpha(t)\}$.
\begin{Lemma}
\label{lem2103}
There exists an $N_0 \in \NNN$ such that for all $n \ge N_0$,  $t_0 \in [0,T]$ and continuous $\alpha: [0,T] \to   [0, \infty)$
\begin{equation}
\prob \{ \CC_n^0(t) >2 \alpha(t), \mbox{ for some } 0< t \le t_0 \}  \le 16nT^2 \sup_{0 \le s \le t_0} \left \{\alpha(s)\prob \{ \CC_n^0(s) > \alpha(s)\} \right \}.\end{equation}
\end{Lemma}
\textbf{Proof:}
Fix $N_0 \in \NNN$ such that for all $n \ge N_0$, $\sup_{s \in [0, T]}\{a_n^*(s) \vee b_n^*(s)\} \le 2$.  Consider now $n \ge N_0$.
Define  $\tau = \inf \{ t >0: \CC_n^0(t) > 2\alpha(t) \}$.  Then
\be
\label{ins1614}
\prob\{\CC_n^0(t) > 2\alpha(t) \mbox{ for some } t \in [0, t_0]\} = \prob\{\tau \le t_0\}.
\ee
Denote by $\CC^0_n \squig_t \CC^s_n$ the event that components $\CC^0_n$ and $\CC^s_n$ merge at time $t$.  By convention this event is taken to
be an empty set if no doubleton is born at time instant $s$. 
Then
$$
\{\tau = t\} = \{\CC^0_n(t-) < 2 \alpha(t)\} \cap \{\CC^0_n(t-) + \CC_n^s(t-) \ge 2 \alpha(t); \CC^0_n \squig_t \CC^s_n, \mbox{ for some }  s < t \}.$$
Next note that
\begin{itemize}
\item  Since $a_n^*(s) \le 2$, the rate at which doubletons are born can be bounded by $2n$.
\item Given a doubleton was born at instant $s$, the event $\{\CC^0_n \squig_u \CC^s_n, \mbox{ for some }  u \in (t, t+dt]\}$ occurs, conditionally on $\FF_t$, with
probability $\frac{1}{n} \CC^0_n(t)\CC^s_n(t)b^*_n(t) dt$. This probability, using the fact that $b^*_n(s) \le 2$ and $\CC_n^s(t) \le n$, on the event 
$\{\CC^0_n(t) < 2 \alpha(t)\}$
is bounded by $4\alpha(t) dt$. 
\item  $\prob \{\CC^0_n(t) + \CC_n^s(t) \ge 2 \alpha(t)\}$ is bounded by $2\prob \{\CC^0_n(t) \ge \alpha(t)\}$.
\end{itemize}
Using these observations we have the following estimate
\beq
\prob\{\tau \le t_0\} 
&\le& 
\E\int_{[0, t_0]} {\bf 1}_{\{\CC^0_n(t) < 2 \alpha(t)\}}\left[  \int_{[0,t]} na^*_n(s) \cdot (\frac{1}{n}\clc^0_n(t)\clc^s_n(t)) \cdot (b^*_n(t)) ds\right] dt\\
&\le&\int_{[0, t_0]} \left[  \int_{[0,t]} 2n \cdot 2\prob(\CC^0_n(t) \ge \alpha(t)) \cdot (4\alpha(t)) ds\right] dt\\
&\le& \int_{[0, t_0]} (2nt) \cdot 2\prob(\CC^0_n(t) \ge \alpha(t)) \cdot (4\alpha(t)) dt \\
& \le & 16nT^2 \sup_{t \in [0, t_0]}\left\{\alpha(t)\prob\{ \CC_n^0(t) > \alpha(t)\}\right\}.
\eeq
Result follows on combining this estimate with \eqref{ins1614}.
\qed\\

The following proposition will be proved in Section \ref{sec:proof-prop-reduced}.
\begin{Proposition}
\label{prop:reduced}
Given $\eta \in (0, \infty)$ and $\gamma \in (0, 1/5)$, there exist $B,C, N_1 \in (0, \infty)$ such that for all $n \ge N_1$
\begin{equation}
\prob \set{ \CC_n^0(t) \ge  m(n,t)/2 } \le Cn^{-\eta} \mbox{ for all }  0<t<t_c-n^{-\gamma},\label{ins2046n}
\end{equation}
where $m(n,t)$ is as defined in \eqref{ins1711}.
\end{Proposition}
\textbf{Remark:}  Intuitively, one has that in the subcritical regime, i.e. when $t< t_c$,  $\prob\{ \CC_n^0(t) > m  \} < d_1 e^{-d_2 m}$ for some constants $d_1,d_2$. This suggests a bound as in (\ref{ins2046n}) for each fixed $t < t_c$.  However, the constants $d_1$ and $d_2$ depend on $t$, and in fact one expects
that, $d_2(t) \to 0$ when $t \uparrow t_c$. On the other hand, in order to prove the above proposition one requires estimates that are uniform for all $t < t_c -n^{-\gamma}$ as $n\to\infty$. This analysis is substantially more delicate as will be seen in subsequent sections. \\

We now prove Proposition  \ref{prop:com-size} using the above results.\\

\textbf{Proof of Proposition \ref{prop:com-size}:} Fix $\gamma \in (0, 1/5)$ and fix $\eta > 2+ 2\gamma$.
Let $B,C, N_1$ be as determined in Proposition \ref{prop:reduced} for this choice of $\eta, \gamma$ and let $m(n,t)$ be as defined in \eqref{ins1711}.  
Without loss of generality we can assume that $N_1 \ge N_0$ where $N_0$ is as in Lemma \ref{lem2103}.  Then applying Lemmas \ref{lem2105} and \ref{lem2103}
with $t_0 = t_c - n^{-\gamma}$ and $\alpha(t)=m(n,t)$, we have 
\beq
\prob (I_n(t) \ge m(n,t) , \mbox{ for some } 0< t < t_c-n^{-\gamma} \} &\le & nT \prob (\CC^0_n(t) \ge m(n,t) , \mbox{ for some } 0< t < t_c-n^{-\gamma} \}\\
&\le & 16 n^2T^3 \sup_{s \in [0, t_c-n^{-\gamma}]}\left \{m(n,s) \prob\{\CC^0_n(s) \ge  m(n.s)/2\}\right \}\\
&\le & 16CBn^{2-\eta + 2\gamma}T^3 (\log n)^4.
\eeq
Since $\eta > 2+ 2\gamma$, the above probability converges to $0$ as $n \to \infty$.  The result follows.
 \qed \\

\subsection{Some Preparatory Results}
\label{sec:approx-error}

This section collects some results that are helpful in estimating the errors between various models described in Section \ref{sec:summary-model}.

The first lemma estimates the error between $\barx_n(t)\equiv \barx(t) = X_n(t)/n$ and its deterministic limit $x(t)$ defined in \eqref{eqn:xss-def}.

\begin{Lemma}
\label{lemma:error-diff-eqn}
For any $T> 0$,  there exists a  $C(T) \in (0, \infty)$ such that, for all $\gamma_1\in [0, 1/2)$,
\[\prob \{ \sup_{0\leq t\leq T}|\barx_n(t) - x(t)| > \frac{1}{n^{\gamma_1}} \} \leq \exp \{ -C(T) n^{1-2\gamma_1}\}.\]
\end{Lemma}
\textbf{Proof:} 
Recall that  $[n] = \set{1,2,\ldots, n}$. Let $E_n = n^{-1} [n]$ and let $E=[0,1]$. Recall the three types of events described in Section \ref{sec:model-equiv}
that lead to edge formation in the BF model.  Of these only events of type (i) and (ii) lead to a change in the number of singletons.  For the
events of type (i), i.e.
in the case when a doubleton is created, $\barx$ decreases by $2/n$. Two key functions (see \eqref{eqn:ins609}) for this case are
\begin{align*}
	f^*_{-2}(y) &= a^*_n(y)\\
	f_{-2}(y) &=a_0(y)= \frac{1}{2} \left(y^2+ (1-y^2)y\right).
\end{align*}
For the events of type (ii), i.e. in the case when a singleton  attaches to a non-singleton component, $\barx$ decreases by $1/n$. Two key functions (see \eqref{eqn:ins633}) for this case are
\begin{align*}
	f^*_{-1}(y) &= (1-y)c_n^*(y)\\
	f_{-1}(y)&= (1-y)c_0(y) =y(1-y^2)(1-y).
\end{align*}

Note that $0\le f^*_{l}(\barx)\le 1$ for $l=-1,-2$, and  that $\barx(t)$ is a Markov process on the state space $E_n$ for which at time $t$ we have the transitions $\barx(t)\leadsto \barx(t)-1/n$ at rate $n f_{-1}^*(\barx(t))$ and $\barx(t)\leadsto\barx(t) -2/n$ at rate $n f_{-2}^*(\barx(t))$.  Furthermore
\begin{equation}
	|f_{-1}^*(y) - f_{-1}(y)|\leq \frac{2}{n} \qquad |f_{-2}^*(y) - f_{-2}(y)| \leq \frac{5}{n}, \; \mbox { for all } y \in [0,1]. \label{eqn:f1n-approx}
\end{equation}
Let $Y_{-1}(\cdot), Y_{-2}(\cdot)$ be independent rate one Poisson processes. Then the process $\barx(t)$ started with $\barx(0)=1$ can be 
constructed (see eg. \cite{kurtz-density}, \cite{ethier-kurtz}) as the unique solution of the stochastic equation 
\begin{equation}
	\label{eqn:stoc-xn-integ}
	\barx(t) = 1 -  \frac{1}{n} Y_{-1}\left(n\int_0^t f_{-1}^*(\barx(s)) ds\right)-\frac{2}{n} Y_{-2}\left(n\int_0^t f_{-2}^*(\barx(s)) ds\right).
\end{equation}

By Equation \eqref{eqn:xss-def}, the limiting function $x(\cdot)$ is the unique solution of the integral equation 
\begin{equation}
	x(t) = 1- \int_0^t f_{-1}(x(s)) ds - \int_0^t 2 f_{-2}(x(s)) ds. \label{eqn:xt-integ}
\end{equation}
Also note that $\forall y,z\in E$
\begin{equation}
	|(f_{-1}(y)+2f_{-2}(y)) - (f_{-1}(z)+2f_{-2}(z)) |\leq 6|y-z|.
	\label{eqn:f-sum-lip}
\end{equation}

Using \eqref{eqn:xt-integ} and \eqref{eqn:stoc-xn-integ} we get 
\begin{align*}
	|\barx(t) - x(t)| \leq A_1^n(t)+ A_2^n(t)+ A_3^n(t)
\end{align*}
where 
\begin{align*} A_1^n(t) = \left|\sum_{l=-1,-2}l \left[\frac{1}{n} Y_{l}\left(n\int_0^t f_{l}^*(\barx(s)) ds\right)- \int_0^t f_{l}^*(\barx(s))ds\right]\right| \leq 4 \sup_{l=-1,-2}\sup_{t< T} \left|\frac{Y_{l}(n t)}{n} -t\right|.
	\end{align*}
and by \eqref{eqn:f1n-approx}
\begin{align*}
	A_2^n(t)  =\left| \int_0^t \sum_{l=-1,-2} l \left[f_{l}^*(\barx(s)) - f_{l}(\barx(s)) \right] ds \right| \leq \frac{7}{n}T.
\end{align*} 
and finally by \eqref{eqn:f-sum-lip} 
\begin{align*}
	A_3^n(t) & = \left|\int_0^s \sum_{l=-1,-2}l \left[f_{l}(\barx(s)) - f_{l}(x(s))\right] ds\right|\\ 
	  &\leq 6\int_0^t |\barx(s)-x(s)| ds .
\end{align*}
Combining these estimates we get 
\[|\barx(t)-x(t)| \leq \left(\frac{7}{n}+ 4 \sup_{l=-1,-2}\sup_{t\le T} \left|\frac{Y_{l}(n t)}{n} -t\right|\right) + 6\int_0^t |\barx(s)-x(s)| ds . \]
This implies, by Gronwall's lemma (see e.g. \cite{ethier-kurtz}, p498) 
\[\sup_{s\leq T}|\barx(s)-x(s)| \leq \left(\frac{7}{n}+ 4 \sup_{l=-1,-2}\sup_{t\le T} \left|\frac{Y_{l}(n t)}{n} -t\right|\right) e^{6T}. \]
Proof is completed using standard  large deviations estimates for Poisson processes. 
\qed\\

In the next lemma we note some basic properties of the integral operator associated with a kernel $\kappa$ on a finite measure space.
\begin{Lemma}
\label{lemma:prop-kernel}
Let $\kappa$, $\kappa'$ be  kernels  given on a finite measure space $(\XX,\TT,\mu)$.  Assume that $\kappa, \kappa' \in L^2(\mu\times \mu)$.
Denote the  associated integral operators by $\KK$ and $\KK'$ (see \eqref{eqn:calk-def}) and there norms by $\rho(\kappa), \rho(\kappa')$ respectively. Then\\
(i) $\KK$ is a compact operator. In particular $\rho(\kappa) =\|\KK\| \le \|\kappa\|_2 = \left(\int_{\XX \times \XX} \kappa^2(x,y) \mu(dx)\mu(dy) \right )^{1/2} <\infty$.\\
(ii) If $\kappa \le \kappa'$, then $\rho(\kappa) \le \rho(\kappa')$.\\
(iii) $\rho(\kappa+\kappa') \le \rho(\kappa)+\rho(\kappa')$ and $\rho(t \kappa)=t \rho(\kappa)$ for $t \ge 0$.\\
(iv) $|\rho(\kappa)-\rho(\kappa')| \le \rho(|\kappa-\kappa'|)$.\\
(v) $\rho(\kappa) \le \|\kappa\|_{\infty} \mu(\XX).$\\
\end{Lemma}
\textbf{Proof:} (i) is a standard result, see Theorem VI.23 of \cite{reed-simon-fna}.\\
(ii) For any nonnegative $f$ in $L^2(\mu)$, $\KK f (x) \le \KK' f (x)$ pointwise. Thus for such $f$, $\|\KK f\|_2\le \|\KK' f\|_2$.  Result follows on observing that the suprema of 
$\|\KK f\|_2,  \|\KK' f\|_2$
over  $\{f \in L^2: \|f\|_2=1\}$ is the same as the suprema over $\{f \in L^2: \|f\|_2=1, f \ge 0 \}$ .\\
(iii) This follows immediately from the facts that $\|(\KK+\KK')f\|_2 \le \|\KK f\|_2 + \|\KK' f\|_2$ and $\KK(tf) =t \KK f$.\\
(iv) Note that $\kappa \le \kappa' + |\kappa - \kappa'|$ and $\kappa' \le \kappa + |\kappa - \kappa'|$. Result follows on combining this observation
with  (ii) and (iii).\\
(v) This follows immediately from (i) and the fact that $\|\kappa\|_2 \le \|\kappa\|_\infty \mu(\XX)$.\\
  \qed \\

We now present some auxiliary estimates for the IRG model from Section \ref{sec:model-irg}.  The following lemma will be proved in Section \ref{sec:model-branching}.  Recall the definition of a {\em basic structure} from Section \ref{sec:model-irg}.
\begin{Lemma}
\label{theo:irg-first-com}
Let $\{ (\XX, \TT, \mu), \kappa, \phi\}$ be a  basic structure, where $\kappa$ is symmetric.  Suppose that $\mu$ is non-atomic and $\rho(\kappa)=\|\KK\|<1$.
 For fixed $x_0 \in \XX$, denote by $\CC_n^{\sss RG}(x_0)$ the volume of the component of $\RG_n(\kappa)$ that contains $x_0$. 
Define $\CC_n^{\sss RG}(x_0) = 0$ if $x_0$ is not a vertex in $\RG_n(\kappa)$.
Then  for all  $m\in \NNN$\\
\begin{equation}
\label{ins1531}
\prob \{ \CC_n^{\sss RG}(x_0) > m\} < 2 \exp \{-C_1\Delta^2 m\}
\end{equation}
where 
\be \label{ins2134} \Delta=1-\rho(\kappa), \;\; C_1=\frac{1}{8\|\phi\|_{\infty}(1 + 3\|\kappa\|_{\infty}\mu(\XX))}.\ee
\end{Lemma}
The above result will be useful for estimating the size of a given component in $\RG_n(a,b,c)$.  One difficulty in directly using this result
is that the kernel $\kappa_t$ and the weight function $\phi_t$ defined in \eqref{eqn:kappa-def} and \eqref{eqn:phi-def} are not bounded.  We will
overcome this by using a truncation argument.  In order to control the error caused by truncation, the following two results will be useful. For rest of this subsection the type space $(\XX, \TT)$ will be taken to be the {\em cluster space} introduced above \eqref{eqn:phi-def}.

\begin{Lemma}
\label{theo:error-wst}
Given rate functions $(a,b,c)$ and $t \in [0, T]$, let $\mu_t$ be the finite measure on $(\XX, \TT)$ defined as in \eqref{eqn:mu-def}.
Let $\PP_n$  be a Poisson point process on $(\XX, \TT)$ with intensity $n \cdot \mu_t$. Define 
$$ Y_n \stackrel{\sss def}{=} \sup_{(s,w) \in \PP_n} w(t).$$
Then for every $A \in (0, \infty)$
\begin{equation*}
\prob\{Y_n > A\} <  2T \cdot  n(1-e^{-T})^{A/2}.
\end{equation*}
\end{Lemma}
\textbf{Proof:} Let $N$ be the number of points in $\PP_n$, then $N$ is  Poisson with mean $\int_0^t na(s)ds \le nT$. 
Let $\{Z_2^{(i)}\}_{i \ge 1}$ be  independent copies of $Z_2$ (also independent of $N$), where $Z_2$ is a pure jump Markov process on $\NNN$ with initial
condition $Z_2(0) = 2$ and infinitesimal generator $\cla_0$ defined as
$$\cla_0f(k) = k(f(k+1) - f(k)), k \in \NNN, \;\; f: \NNN \to \RRR. $$
Thus $Z_2$ is just a Yule process started with two individuals at time zero. Note that 
$$Y_n \le \sup_{(s,w) \in \PP_n} w(T) \le_d \sup_{1 \le i \le N} Z_2^{(i)}(T),$$
where the first inequality holds a.s and the second inequality uses the fact that $c \le 1$.
Standard facts about the Yule process (see e.g.\cite{norris-mc}) imply that $Z_2^{(i)}(T)$ is distributed as sum of two independent $\mbox{Geom}\{e^{-T}\}$. 
Thus
\begin{align*}
\prob\{Y_n > A\} 
&\le \mathbb{E}(N) \cdot \prob \{ Z_2(T) > A \}\\
&\le nT  \cdot 2 (1-e^{-T})^{A/2}.
\end{align*}
This completes the proof of the lemma. \qed\\

The following corollary follows on taking $A = C\log n$ in the above lemma.

\begin{Corollary}\label{cor1210}
Let $Y_n$ be as in the above lemma and fix $\eta \in (0, \infty)$.  Then there exist $C_1(\eta), C_2(\eta) \in (0, \infty)$ such that
for any    rate functions $(a,b,c)$
$$\prob\{Y_n > C_1(\eta) \log n\} < C_2(\eta) n^{-\eta}, \;\; \mbox{ for all } n \in \NNN .$$
\end{Corollary}

From Section \ref{sec:model-equiv}, recall the definitions of the functions $a_0, b_0, c_0$ associated with the BF model.   The following lemma will allow us to argue that $\RG_n(a_0,b_0,c_0)$ is well approximated by $\RG_n(a,b,c)$ if the rate functions $(a,b,c)$ are sufficiently
close to $(a_0, b_0, c_0)$.  Let $((\XX, \TT, \mu_t), \kappa_t, \phi_t)$ be the basic structure associated with rate functions $(a,b,c)$.  Let $\KK_t$ be the integral
operator defined by \eqref{eqn:calk-def}, replacing $(\mu, \kappa)$ there by $(\mu_t, \kappa_t)$.  Let $\rho_t = \rho(\kappa_t)$.  In order to emphasize the dependance
on rate functions $(a,b,c)$, we will sometimes write $\rho_t=\rho_t(a,b,c)$.  Similar notation will be used for $\kappa_t, \mu_t, \phi_t$ and $\KK_t$.
\begin{Lemma}
\label{theo:error-cont-norm}
Fix rate functions $(a,b,c)$.  Suppose that $\inf_{s \in [0, T]} a(s) > 0$ and for some $\theta \in (0, \infty)$,
$ c(s)\ge \theta s$, for all $s \in [0, T]$.  
Given $\delta > 0$ and $t \in [0, T]$, let  
$$\rho_{+,t} = \rho_t((a+\delta)\wedge 1, (b+\delta)\wedge 1, (c+\delta)\wedge 1),\;
\rho_{-,t} = \rho_t((a-\delta)^+, (b-\delta)^+, (c-\delta)^+).$$
Then there exists $C_2 \in (0, \infty)$ and $\delta_0 \in (0, 1)$ such that for all $\delta \le \delta_0$ and $t \in [0, T]$
$$
\max \{ |\rho_t- \rho_{+,t}|,  |\rho_t-\rho_{-,t}|\}  \le C_2 (-\log \delta)^3 \delta^{1/2}.$$
\end{Lemma}
The proof of the above lemma is quite technical and deferred to Section \ref{sec:analysis-ker}.   \\

The next lemma gives some basic properties of $\rho_t(a_0, b_0, c_0)$. Recall that $t_c$ denotes the critical time for the emergence of the giant component in the BF model. 
 
\begin{Lemma}
\label{theo:error-cont-int}
Let $\rho(t)=\rho_t(a_0, b_0, c_0)$.  Then:\\
(i) $\rho(t)$ is strictly increasing in $t \in [0, T]$;\\
(ii) $ \rho(t_c)=1$;\\
(iii) $\lim_{s \to 0^+} (\rho(t_c)-\rho(t_c-s))/s = \rho'_- (t_c) >0$. 
\end{Lemma}
The proof of the lemma is given in Section \ref{sec:measure-theoretic}. \\

\subsection{Proof of Proposition \ref{prop:reduced}}
\label{sec:proof-prop-reduced}
This section is devoted to the proof of Proposition \ref{prop:reduced}.
Fix $\eta \in (0, \infty)$ and $\gamma \in (0, 1/5)$.\\

\textbf{Step 1: from $\BF_n$ to $\IA_{n,\delta}$ }\\

Let $\gamma_1 = 2/5$ and define $E_n = \{ \sup_{0 \le t \le T} |\barx_n(t)-x(t)| \le n^{-\gamma_1}\}$.

From Lemma \ref{lemma:error-diff-eqn}, 
\begin{equation}
\prob \{ E_n^c \} \leq \exp \{ -C(T) n^{1-2\gamma_1}\} =  \exp \{ -C(T) n^{1/5}\}.
\label{eqn:ineq-step1}
\end{equation}
From \eqref{eqn:a-astar} and recalling that the Lipschitz norm of $a_0$ is bounded by 2 (see \eqref{eqn:a-def}), we have that on $E_n$
$$|a^*(t) - a_0(t)| \le 5n^{-1} + 2 n^{-\gamma_1}, \mbox{ for all } t \in [0, T].$$
Similar bounds can be shown to hold for $b^*$ and $c^*$.  Thus we can find $n_1 \in \NNN$ and $d_1 \in (0, \infty)$ such that, for $n \ge n_1$, on $E_n$
$$a_n^*(t)\le a_0(t)+\delta_n, b_n^*(t)\le b_0(t)+\delta_n, c_n^*(t)\le c_0(t)+\delta_n, \mbox{ for all } t \in [0, T],$$
where $\delta_n = d_1n^{-\gamma_1}$.  Since $a_n^*, b_n^*, c_n^*$ are all bounded by $1$, setting
$(a_0(t)+\delta_n)\wedge 1 = a_{n,\delta}$ and similarly defining $b_{n,\delta}, c_{n,\delta}$,
we in fact have that
$$a_n^*(t) \le a_{n,\delta}(t), b_n^*(t)\le b_{n,\delta}(t), c_n^*(t) \le c_{n,\delta}(t), \mbox{ for all } t \in [0, T].$$
Let $\CC_{n,\delta}^{\sss IA}(t)$ denote the size of the first component in 
$\IA_n(a_{n,\delta}, b_{n,\delta}, c_{n,\delta})_t$.  From Lemma 
 \ref{lemma:couple-bfia}, we have for any $m \in \NNN$
\begin{equation}
\prob\{ \CC_n^{0}(t) > m, E_n\} \le \prob\{ \CC_{n,\delta}^{\sss IA}(t) > m, E_n \} \le \prob\{ \CC_{n,\delta}^{\sss IA}(t) > m\}.
\label{eqn:prob-step1}
\end{equation}

\textbf{Step 2: from  $\IA_{n,\delta}$ to $\RG_{n,\delta,A}$}\\
For $t \in [0, T]$, and rate functions $a_{n,\delta}, b_{n,\delta}, c_{n,\delta}$, consider the IRG model
$\RG_n(\kappa_{t,\delta},\mu_{t,\delta}, \phi_t)$, where 
$\kappa_{t, \delta} = \kappa_t(a_{n,\delta}, b_{n,\delta}, c_{n,\delta})$ and  $\mu_{t,\delta}$ is the measure for the IRG model corresponding to these rate functions as defined in \eqref{eqn:mu-def}.  Let $A_n = C_1(\eta) \log n$, where $C_1(\eta)$ is as in Corollary \ref{cor1210}.
Consider the following truncation of the kernel $\kappa_{t,\delta}$ and weight function $\phi_t(s,w) = w(t)$:
$$\kappa_{t,\delta,A}(\bfx,\bfy)=\kappa_{t,\delta} (\bfx, \bfy) {\bf1}_{\{w(T) \le A_n\} }{\bf1}_{\{\tilde w(T) \le A_n\} }, \; \bfx = (s,w), \bfy = (r, \tilde w)$$
and
$$\phi_{t,A}(s,w)=\phi_t(s,w) {\bf1}_{\{w(T) \le A_n\} }.$$
Then $\|\phi_{t,A}\|_\infty \le A_n$, $\|\kappa_{t,\delta,A}\|_\infty \le T A_n^2$ .\\

Recall the Poisson point process $\PP_t(a,b,c)$ associated with rate functions $(a,b,c)$, introduced below \eqref{eqn:phi-def} and write
$\PP_{t,\delta} = \PP_t(a_{n,\delta}, b_{n,\delta}, c_{n,\delta})$.  Let $Y_{n,\delta} = \sup_{(s,w) \in \PP_{t,\delta}}w(T)$.  From Corollary \ref{cor1210}
\begin{equation}
\prob\{Y_{n,\delta} > A_n\} <C_2(\eta) n^{-\eta}.
\label{eqn:ineq-step2}
\end{equation}
Let $\CC_{n,\delta}^{\sss RG}(t) = \CC_{n,t}^{\sss RG}(a_{n,\delta}, b_{n,\delta}, c_{n,\delta})$  be the volume of the `first' component in $\RG_{n,t}(a_{n,\delta}, b_{n,\delta}, c_{n,\delta})\equiv \RG_{n}(\kappa_{t,\delta}, \mu_{t,\delta}, \phi_{t})$.  Then
from Lemma \ref{lemma:rgiva-irg-def}
\be
\prob\{\CC_{n,\delta}^{\sss IA}(t) > m\} \le \prob\{ \CC_{n,\delta}^{\sss RG}(t) >m\}.\label{ins1632} \ee
Letting $\CC_{n,\delta, A}^{\sss RG}(t)$ denote the volume of the first component in $\RG_{n}(\kappa_{t,\delta, A}, \mu_{t,\delta}, \phi_{t,A})$, namely the random graph formed using the truncated kernel. Then 
\beqn
 \prob\{ \CC_{n,\delta}^{\sss RG}(t) >m\}
&\le& \prob\{Y_{n,\delta} > A_n\}+ \prob\{ \CC_{n,\delta}^{\sss RG}(t) >m, Y_{n,\delta} \le A_n\}\nonumber \\
&=& \prob\{Y_{n,\delta} > A_n\}+ \prob\{ \CC_{n,\delta,A}^{\sss RG}(t) >m, Y_{n,\delta} \le A_n\} \label{eqn:trunc-no-effect}\nonumber\\
&\le& \prob\{Y_{n,\delta} > A_n\}+ \prob\{ \CC_{n,\delta,A}^{\sss RG}(t) >m\} \label{eqn:prob-step2}
\eeqn

\textbf{Step 3: Estimating $\CC_{n,\delta,A}^{\sss RG}$}\\
We will apply Lemma \ref{theo:irg-first-com}, replacing $\{(\XX, \TT, \mu), \kappa, \phi\}$ there by
$\{(\XX, \TT, \mu_{t,\delta}), \kappa_{t,\delta, A}, \phi_{t,A}\}$, where $t \in (0, t_c - n^{-\gamma})$.
From \eqref{ins1531} we have
\be 
\label{ins1532}
\prob\{ \CC_{n,\delta,A}^{\sss RG}(t) >m\} \le 2 \exp \{-C_1\Delta^2m\},
\ee
where $$C_1 = \frac{1}{8\|\phi_{t,A}\|_{\infty}(1+3\|\kappa_{t,\delta, A}\|_{\infty}\mu_{t,\delta}(\XX))}, $$
and $\Delta = 1 - \rho(\kappa_{t,\delta, A})$.  We now estimate $\rho(\kappa_{t,\delta, A})$.
Since  $ \kappa_{t,\delta, A} \le \kappa_{t,\delta}$, by (ii) of Lemma \ref{lemma:prop-kernel}, we have $\rho(\kappa_{t,\delta, A}) \le \rho(\kappa_{t,\delta})$.
Note that rate functions $(a_0, b_0, c_0)$ satisfy conditions of Lemma \ref{theo:error-cont-norm}.  Thus, recalling that
$\delta_n = d_1 n^{-2/5}$, we have from this result, that for some $d_2 \in (0, \infty)$,
$\rho(\kappa_{t,\delta})< \rho(\kappa_t)+ d_2(\log n)^3 n^{-1/5}$, for all $t \le T$.  Here $\kappa_t = \kappa_t(a_0, b_0, c_0)$.

Next, by Lemma \ref{theo:error-cont-int}, there exists  $d_3 \in (0, \infty)$ such that 
$\rho(\kappa_t)<1-d_3(t_c-t)$ for all $t \in [0,t_c)$. Combining these estimates, we have for  $t < t_c-n^{-\gamma}$,
\begin{align*}
\rho(\kappa_{t,\delta, A}) 
&< 1-d_3(t_c-t) +d_2(\log n)^3 n^{-1/5}.
\end{align*}
Recalling that $\gamma \in (0, 1/5)$ we have that, for some  $n_2 \in (n_1, \infty)$ and $d_4 \in (0, \infty)$,  
$$ \rho(\kappa_{t,\delta, A}) \le  1-d_4 (t_c-t), \mbox { for all }  t \in (0, t_c-n^{-\gamma}) \mbox{ and } n \ge n_2.$$
Using this estimate in \eqref{ins1532} and recalling that  $\|\phi_{t,A}\|_\infty \le A_n$, $\|\kappa_{t,\delta,A}\|_\infty \le T A_n^2$ , we have that for some $d_5 \in (0, \infty)$
\be
\label{ins1625}
\prob\{ \CC_{n,\delta,A}^{\sss RG}(t) >m\} \le 2 \exp \{ -\frac{d_5}{(\log n)^3} (t_c-t)^2 m\}, \mbox{ for all } m \in \NNN, t \in (0, t_c-n^{-\gamma}) \mbox{ and } n \ge n_2.\ee

\textbf{Step 4: Collecting estimates:}\\
Combining \eqref{eqn:ineq-step1}, \eqref{eqn:prob-step1}, \eqref{ins1632}, \eqref{eqn:ineq-step2}, \eqref{eqn:prob-step2} and \eqref{ins1625}, we have 
\beqn
\prob\{ \CC_n^{0}(t) > m\} &\le& \prob\{ E_n^c\} +\prob\{Y_{n,\delta} > A_n\} +\prob\{ \CC_{n,\delta,A}^{\sss RG}(t) >m \}\nonumber \\
&\le &
 e^{-C(T) n^{1/5}} + C_2(\eta) n^{-\eta} +2 \exp\{-d_5 \frac{(t_c-t)^2}{(\log n)^3} m\}.
\eeqn
Finally, result follows on replacing $m$ in the above display with
$\frac{\eta (\log  n)^4}{d_5 (t_c-t)^2}$. \qed \\

The following lemma will be used in the proof of Lemma \ref{theo:error-cont-int}.  We will use notation and arguments similar to that in the
proof of Proposition \ref{prop:reduced} above.

\begin{Lemma}
\label{lemma:logn-bound}
Let $(a,b,c)$ be rate functions. Fix $t \in [0, T]$.  Let $I_n^{\sss IA}(t)$ denote the largest component in $\IA_n(a,b,c)_t$.
Suppose that $\rho_t(a,b,c)< 1$. Then for some $C_0 \in (0, \infty)$
$$\prob\{I_n^{\sss IA}(t) > C_0 (\log n)^4\} \to 0 \mbox{ when } n \to \infty.$$
\end{Lemma}
\textbf{Proof: } 
 Let $\CC_n^{\sss IA}(t)$ be the first component of $\IA_n(a,b,c)_t$. Then an elementary argument (cf. proof of Lemma \ref{lem2105}) shows
 that for $m > 0$
$$\prob\{I_n^{\sss IA}(t) > m \} \le T n \prob\{ \CC_n^{\sss IA}(t) > m\}.$$
By an argument as in \eqref{eqn:prob-step2}, we have 
$$
 \prob\{\CC_n^{\sss IA}(t) > m\} 
 \le \prob\{ \CC_n^{\sss RG}(t) > m\}\\
 \le \prob\{ Y_n > A_n\}+\prob\{ \CC_{n,A}^{\sss RG}(t) > m\},
$$
where $\CC_n^{\sss RG}$, $Y_n$ and $\CC_{n,A}^{\sss RG}$ correspond to $\CC_{n,\delta}^{\sss RG}$, $Y_{n,\delta}$ and $\CC_{n,\delta,A}^{\sss RG}$
introduced above in the proof of  Proposition \ref{prop:reduced}, with $(a_{n,\delta}, b_{n,\delta}, c_{n,\delta})$ replaced with $(a,b,c)$.
From Corollary \ref{cor1210} we can find $d_1 \in (0, \infty)$ such that
$\prob(Y_n \ge d_1 \log n) = O(n^{-2})$.   Let $A_n = d_1 \log n$. 
 Then, recalling that $\rho_t(a,b,c) < 1$, we gave by Lemma \ref{theo:irg-first-com} that, for some $d_2 \in (0, \infty)$,
$$ \prob\{ \CC_{n,A}^{\sss RG}(t) > m\} < 2 \exp \{ - d_2m/ (\log  n)^3\}.$$
Taking $m= \frac{2}{d_3} (\log n)^4$, we have $ \prob\{ \CC_{n,A}^{\sss RG}(t) > m\}=O(n^{-2})$.  Combining the above
estimates we have   $\prob\{I_n^{\sss IA}(t) > \frac{2}{d_3} (\log n)^4 \}=O(n^{-1})$. The result follows.\qed \\

\subsection{Proof of Lemma \ref{theo:irg-first-com}: A branching process construction}
\label{sec:model-branching}
The key idea in the proof of Lemma \ref{theo:irg-first-com} is the coupling of the breadth first exploration of components in the IRG model with  a certain 
continuous type branching process.   
This coupling will reduce the problem of establishing the estimate in  Lemma \ref{theo:irg-first-com} to a similar bound
 on the total volume of the branching process (Lemma \ref{theo:branching}). We refer the reader to  \cite{bollobas-riordan-janson} where a similar coupling in a setting where the type space $\XX$ is finite
 using a finite-type branching process is constructed.
 In this subsection we will give the proof of Lemma  \ref{theo:irg-first-com} using
 Lemma \ref{theo:branching}. Proof of the latter result is given in Section \ref{sec:proof-branching}.\\

Throughout this section we will fix a basic structure   $\{(\XX,\TT,\mu), \kappa, \phi\}$, where $\kappa$ is a symmetric kernel, and a $x_0 \in \XX$. 
  Let $\RG_n(\kappa)$ be the IRG constructed using this structure as in  Section \ref{sec:model-irg}.
We now describe a branching process associated with the above basic structure.  The process starts in the $0$-th generation with   a single vertex of type $x_0 \in \XX$
and in the $k$-th generation,  a vertex $x$   will have offspring, independently of the remaining $k$-th generation vertices,
 according to a Poisson point process on $\XX$ with intensity $\kappa(x,y)\mu(dy)$ to form the $(k+1)^{th}$ generation. We denote this branching process as $\BP(x_0)$. \\

Denote by $\{\xi_i^{(k)}\}_{i=1}^{N_k} \subset \XX$  the $k^{th}$ generation of the branching process. Define the volume of the $k$-th generation as $G_k=\sum_{i=1}^{N_k} \phi(\xi_i^{(k)})$. The \textbf{total volume} of $\BP(x_0)$ is defined as $G=G(x_0)=\sum_{k=0}^{\infty} G_k$.\\

The following lemma, proved at the end of the section, shows that $\CC_n^{\sss RG}(x_0)$ is stochastically dominated by $G(x_0)$.
\begin{Lemma}
\label{lemma:irg-branching}
For all $m_1 \in \NNN$,
$$\prob\{ \CC_n^{\sss RG}(x_0) > m_1\} \le \prob\{ G(x_0) > m_1 \}.$$
\end{Lemma}
Next lemma, proved in Section \ref{sec:proof-branching} shows that the estimate in Lemma \ref{theo:irg-first-com} holds with
$\CC_n^{\sss RG}(x_0)$ replaced by $G(x_0)$.
\begin{Lemma}
\label{theo:branching}
  Suppose that $\rho(\kappa)=\|\KK\|<1$.
  Then  for all  $m\in \NNN$\\
\begin{equation}
\label{ins1531b}
\prob \{ G > m\} < 2 \exp \{-C_1\Delta^2 m\}
\end{equation}
where $\Delta$ and $C_1$ are as in \eqref{ins2134}.
\end{Lemma}

Using the above lemmas we can now  complete the proof of Lemma \ref{theo:irg-first-com}.\\

\noindent \textbf{Proof of Lemma \ref{theo:irg-first-com}: } Proof is  immediate from Lemmas \ref{theo:branching} and \ref{lemma:irg-branching} . \qed \\

We conclude this section with the proof of Lemma \ref{lemma:irg-branching}.\\

\noindent {\bf Proof of Lemma \ref{lemma:irg-branching}:}  Without loss of generality assume that $\CC_n^{\sss RG}(x_0)\neq 0$.  We now
explore the component $\CC_n^{\sss RG}(x_0)$ in the standard breadth first manner.   

Define the sequence of \textbf{unexplored sets} $\{U_m\}_{m \ge 0}$ and the set of \textbf{removed vertices} $\{R_m\}_{m \ge 0}$ iteratively as follows: 
Let
$R_0=\emptyset, U_0=\{x_0\}$ and $y_1=x_0$. 
Suppose we have defined $R_j, U_j$, $j = 0, 1, \cdots, m-1$ and $U_{m-1}= \{y_m, y_{m+1}, \cdots y_{t_m}\}$.  Then set
\begin{align*}
R_m &=R_{m-1} \cup \{y_m\}\\
U_m &=U_{m-1} \cup E_m \setminus \{y_m\}
\end{align*}
where \[E_m=\{ x \in \XX : x \text{ is a neighbor of $y_m$ in $\RG_n(\kappa)$ (i.e. $\{x,y_m\}$ is an edge) and } x \notin R_{m-1} \cup U_{m-1} \}.\] 
If $U_{m-1}=\emptyset$ we set $U_j=E_j = \emptyset$ and $R_j=R_{m-1}$ for all $j \ge m$.
Thus $U_m$ are the vertices at step $m$ that have been revealed by the exploration but whose neighbors have not been explored yet. 
Note that the number of vertices in $R_{m-1}\cup U_{m-1}$ equals $t_m$.  Label the vertices in $E_m$ as $y_{t_m+1}, y_{t_m+2}, \ldots y_{t_m+|E_m|}$.
With this labeling we have a well defined specification of the sequence $\{R_j, U_j, E_{j+1}\}_{j \in \NNN_0}$.
%
%
%
%
Note that $\CC_n^{\sss RG}(x_0)= m_0$ if and only if $U_{m_0-1}\neq \emptyset, U_{m_0}=\emptyset$ and $|R_{m_0}|=m_0$. 

We will now argue that for every $m \in \NNN$, conditioned on $\{U_{m-1}, R_{m-1}\}$, $E_m$ is a Poisson point process on the space $\XX$
with intensity 
$$\Lambda_m^*(dx)=\beta_m(x) (\kappa(y_m, x)\wedge n) \mu(dx),$$
where $\beta_m: \XX \to [0, 1]$ is given as  $\beta_1 \equiv 1$ and, for $m >1$,
$$\beta_m(x)=\Pi_{y \in R_{m-1}}\left[1- \left(\frac{\kappa(y,x)}{n}\wedge 1\right)\right], \; x \in \XX .$$
Consider first the case $m=1$.
Denote the Poisson point process on $(\XX, \TT)$ used in the construction of $\RG_n(\kappa,\mu)$ by $N_n(\kappa,\mu)$.
From the complete independence property of Poisson point processes and the non-atomic assumption on $\mu$, conditioned on the existence of a vertex $x_0$
in $N_n(\kappa,\mu)$, $N_n(\kappa,\mu)\setminus \{x_0\}$ is once again a Poisson point process with intensity $n \cdot \mu(dx)$ on $\XX$. Also, conditioned on
$N_n(\kappa,\mu)$, a given type $x$ vertex in $N_n(\kappa,\mu)$ would be connected to $x_0$ with probability $(\kappa(x_0,x)/n)\wedge 1 $.  Thus the neighbors of $x_0$, namely $E_1$, define a Poisson point process with intensity $(\kappa(x_0,x)\wedge n)\mu(dx)$. This proves the above statement on $E_m$ with $m=1$.

Consider now
$m > 1$. Since $\mu$ is non-atomic and $U_{m-1}\cup R_{m-1}$ consists of only finitely many elements, it follows that conditioned on  vertices in $R_{m-1} \cup U_{m-1}$ belonging to $N_n(\kappa,\mu)$, $N_n(\kappa,\mu)\setminus (R_{m-1} \cup U_{m-1})$ 
is once again a Poisson point process on $\XX$  with intensity  $n \cdot \mu(dx)$. Note that a vertex $x \in N_n(\kappa,\mu)\setminus (R_{m-1} \cup U_{m-1})$ is in $E_m$ 
if and only if $x$ is a neighbor of $y_m$ and $x$ is not a neighbor of any vertex in $R_{m-1}$. So conditioned on  $\{R_{m-1}, U_{m-1}\}$,
the probability that $x$ is in $E_m$ equals 
\[(\kappa(y_m, x)/n \wedge 1) \cdot \Pi_{y \in R_{m-1}} [1- (\kappa(y,x)/n) \wedge 1].\]
From this and the fact that the edges in $\RG_n(\kappa, \mu)$ are placed in a mutually independent fashion, it follows that the points in $E_m$, conditioned on $\{R_{m-1}, U_{m-1}\}$, describe a Poisson point process with intensity
\begin{align*}
& n \mu(dx) \cdot (\kappa(y_m, x)/n \wedge 1) \cdot \Pi_{y \in R_{m-1}} [1- (\kappa(y,x)/n) \wedge 1]\\
&= \Pi_{y \in R_{m-1}} [1- (\kappa(y,x)/n) \wedge 1] \cdot (\kappa(y_m,x) \wedge n)\mu(dx)\\
&=\Lambda_m^*(dx).
\end{align*}
 Thus conditioned on $\{R_{m-1}, U_{m-1}\}$, $E_m$ is a Poisson point process with the claimed intensity.

Next note that one can carry out an analogous breadth first exploration of $\BP(x_0)$.  Denoting the corresponding vertex sets once more by
$\{R_j, U_j, E_{j+1}\}_{j \in \NNN_0}$ we see that conditioned on  $\{R_{m-1}, U_{m-1}\}$, $E_m$ is a Poisson point process with intensity $\kappa(y_m,x)\mu(dx)$.

 As $ 0 \le \beta_m(x) \le 1$ and $\kappa \wedge n \le \kappa$, we can now construct a coupling between $\BP(x_0)$ and $\CC_n^{RG}(x_0)$ by first
constructing $\BP(x_0)$ and then by iteratively rejecting each offspring of type $x$ in $E_m$ (and all of its descendants) with probability
 $$1-\frac{\beta_m(x)(\kappa(y_m,x)\wedge n)}{\kappa(y_m,x)}.$$
The lemma is now immediate. \qed \\






\subsection{Proof of Lemma \ref{theo:branching}}
\label{sec:proof-branching}
Assume throughout this subsection, without loss of generality, that $\max\{\|\phi\|_{\infty}, \|\kappa\|_{\infty}, \mu(\XX)\} < \infty$. Recall that $\kappa$ is a symmetric kernel.
Define, for $k \in \NNN$,  the kernels $\kappa^{(k)}$ recursively as follows. $\kappa^{(1)}=\kappa$ and for all $k \ge 1$
$$\kappa^{(k+1)}(x,y)=\int_{\XX} \kappa^{(k)}(x,u)\kappa(u,y) \mu(du).$$
Recall that $\{\xi_i^{(k)}\}_{i=1}^{N_k}$ denotes the $k$-th generation of $\BP(x_0)$ and note that it describes a
 Poisson point process with intensity $\kappa^{(k)}(x_0,y) \mu(dy)$.  This observation allows us to compute exponential moments of the form in the lemmas below.
\begin{Lemma}
\label{theo:bp-onestep}
Let $g: \XX \to \RRR_+$ be a bounded measurable map. Fix $\delta>0$ and let $0<\epsilon < \log(1+\delta)/\|g\|_{\infty}$.
Then
$$\E \exp\{\epsilon \sum_{i=1}^{N_1} g(\xi_i^{(1)})\} \le \exp\{ \epsilon (1+\delta) (\KK g)(x_0)\}.$$
\end{Lemma}
\textbf{Proof:}  Fix $\delta, \epsilon$ as in the statement of the lemma.
By standard formulas for  Poisson point processes
\begin{align*}
\E \exp\{\epsilon \sum_{i=1}^{N_1} g(\xi_i^{(1)})\} 
&= \exp \{ \int_\XX \kappa(x_0,u) (e^{\epsilon g(u)}-1) \mu(du)\}\\
&\le \exp \{ \int_\XX \kappa(x_0,u) (1+\delta) \epsilon g(u) \mu(du) \}\\
&= \exp\{ \epsilon (1+\delta) (\KK g)(x_0)\},
\end{align*}
where the middle inequality follows on noting that  $e^{\epsilon g(u)}-1 \le (1+\delta) \epsilon g(u)$, whenever $\epsilon g(u) < \log(1+\delta)$.\\
\qed\\
Using the above lemma and a recursive argument, we obtain the following result. Recall that $G_k = \sum_{i=1}^{N_k} \phi(\xi_i^{(k)})$ denoted the volume of generation $k$ where volume is measured using the function $\phi$. 
\begin{Lemma}
\label{theo:bp-keybound}
Fix $k \in \NNN$ and $\delta > 0$. Given a weight function $\phi$, define $\phi_0=\phi + \sum_{i=1}^k (1+\delta)^i \KK^i \phi$.  Then for all $\epsilon \in (0, \frac{\log(1+\delta)}{\|\phi_0\|_{\infty}})$
\be \label{ins2253}\E \exp \{ \epsilon \sum_{i=0}^k G_i  \} \le \exp \{ \epsilon [ \phi(x_0) + \sum_{i=1}^k (1+\delta)^i \KK^i \phi (x_0) ] \} = \exp \{\epsilon \phi_0(x_0)\}.\ee
\end{Lemma}

\textbf{Proof:} Define $\{\phi_i\}_{i=0}^k$  using a backward recursion, as follows. Let $\phi_k=\phi$. For $0 \le i <k $
$$\phi_{i}=\phi + (1+\delta) \KK \phi_{i+1}.$$
Let $\FF_l = \sigma \{ \{\xi_i^{(k)}\}_{i=1}^{N_k}, k = 1, \cdots l\}$.
We will show recursively, as $l$ goes from $k$ to $0$, that 
\begin{equation}
\E [ \exp \{ \epsilon \sum_{i=l}^k G_i \} | \FF_l] \le \exp \{ \epsilon \sum_{i=1}^{N_l} \phi_l (\xi_i^{(l)}) \} \}.
\label{eqn:bp-induction}
\end{equation}
The lemma is then immediate on setting   $l=0$ in the above equation.\\

When $l=k$, \eqref{eqn:bp-induction} is in fact an equality, and so \eqref{eqn:bp-induction} holds trivially for $k$. Suppose now that
\eqref{eqn:bp-induction} is true for $l+1$, for some $l \in \{0, 1, \cdots, k-1\}$.  Then
\begin{align*}
\E [ \exp \{ \epsilon \sum_{i=l}^k G_i \} | \FF_l] 
&=  \exp \{ \epsilon G_l\} \E[ \E [ \exp \{ \epsilon \sum_{i=l+1}^k G_i \} | \FF_{l+1}] | \FF_l]\\
&\le \exp \{ \epsilon G_l\} \E[ \exp \{ \epsilon \sum_{i=1}^{N_{l+1}} \phi_{l+1} (\xi_i^{(l+1)}) \}| \FF_l] \hspace{.5in}\\
&\le \exp \{ \epsilon G_l\} \exp \{ \epsilon (1+\delta) \sum_{i=1}^{N_l} \KK \phi_{l+1}(\xi_i^{(l)})\} \hspace{.5in} \\
&= \exp \{ \epsilon \sum_{i=1}^{N_l} \phi (\xi_i^{(l)})\} \exp \{ \epsilon (1+\delta) \sum_{i=1}^{N_l} \KK \phi_{l+1}(\xi_i^{(l)})\}\\
&= \exp \{ \epsilon \sum_{i=1}^{N_l} [\phi (\xi_i^{(l)}) + (1+\delta)\KK \phi_{l+1}(\xi_i^{(l)})]\} \\
&= \exp \{ \epsilon \sum_{i=1}^{N_l} \phi_l (\xi_i^{(l)}) \}.
\end{align*}
For the first inequality above we have used the fact that by assumption \eqref{eqn:bp-induction} holds for $l+1$ and for the
second inequality we have applied Lemma \ref{theo:bp-onestep} along with the observation that $\epsilon \|\phi_l\|_{\infty} < \log(1+\delta)$ holds for all $l=1,2,...,k$,
since for all $l$,  $\phi_l \le \phi_0$ and $\epsilon \|\phi_0\|_{\infty} < \log(1+\delta)$.

This completes the recursion and the result follows. \qed\\

To emphasize that $\phi_0$ in the above lemma depends on $\delta$ and $k$, write $\phi_0 = \phi^{(k)}_{\delta}$.  Note that $\phi^{(k)}_{\delta}$ is increasing in $k$.
Let $\phi^*_{\delta} = \lim_{k \to \infty} \phi^{(k)}_{\delta}$.  The following corollary follows on sending $k \to \infty$ in \eqref{ins2253}.

\begin{Corollary}
\label{cor2310}
Fix $\delta > 0$ and $ \epsilon \in (0, \log(1+\delta)/\|\phi^*_{\delta}\|_{\infty})$.  Then
\begin{equation}
\E\{ \exp{\epsilon G}\} \le \exp\{ \epsilon \phi^*_{\delta}(x_0)\}.
\label{eqn:bp-mct}
\end{equation}
\end{Corollary}

\begin{Lemma}
\label{theo:bp-contraction}
For $n \in \NNN$ and $x \in \XX$
$$\KK^n \phi (x) \le \rho^{n-1} \|f_x\|_2 \|\phi\|_2, \mbox{ where } f_x(\cdot)=\kappa(x,\cdot) \mbox{ and } \rho = \rho(\kappa).$$
\end{Lemma}
\textbf{Proof:} Note that
\beq
\KK^n \phi (x) 	&=&\int_\XX \kappa^{(n)}(x,u)\phi(u) \mu(du)
				\le \| \int_\XX f_x(u) \kappa^{(n-1)}(u,\cdot) \mu(du)\|_2 \|\phi\|_2\\
				&=&\| \KK^{n-1} f_x\|_2 \|\phi\|_2
				\le \rho^{n-1} \|f_x\|_2 \|\phi\|_2.
				\eeq

\qed \\

Now we can finish the proof of Lemma \ref{theo:branching}.\\
\textbf{Proof of Lemma \ref{theo:branching}:} 
Observing that $\|\phi\|_2 \le \|\phi\|_{\infty} \mu(\XX)^{1/2} $ and $\|f_x\|_2 < \|\kappa\|_{\infty}\mu(\XX)^{1/2}$, we have
for $\delta \in (0, \infty)$ such that $(1+\delta)\rho<1$, and $x \in \XX$
\begin{align*}
\phi^*_{\delta}(x) 
&= \phi(x)+\sum_{i=1}^{\infty} (1+\delta)^i \KK^i \phi (x)\\
&\le \|\phi\|_{\infty} + \|f_x\|_2 \|\phi\|_2 (\sum_{i=1}^{\infty} (1+\delta)^i \rho^{i-1}) \hspace{.5in} \\
&\le \|\phi\|_{\infty} + \|\kappa\|_{\infty}\|\phi\|_{\infty}\mu(\XX)\frac{(1+\delta)}{1-(1+\delta)\rho},\end{align*}
where the first inequality above follows from Lemma \ref{theo:bp-contraction}.
Setting $\delta = \frac{\Delta}{2}$, we see
$$(1+\delta)\rho=(1+\Delta/2)(1-\Delta)<1-\Delta/2.$$
Using this and that $\Delta < 1$, we have
$$
\phi^*_{\delta}(x)\le \|\phi\|_{\infty} \left(1+ \frac{3\|\kappa\|_{\infty}\mu(\XX)}{\Delta} \right)\equiv d_1.$$
Let $\epsilon = \log(1+\delta)/(2d_1)$.  Clearly $\epsilon \in (0, \log(1+\delta)/\|\phi^*_{\delta}\|_{\infty})$.
Using Corollary \ref{cor2310} we now have that
\begin{align*}
\prob\{ G > m \} 
&\le \exp\{-\epsilon m\} \exp\{\epsilon \phi^*_{\delta}(x_0)\} \\
&\le  \exp\{-\epsilon m\} \exp \{\frac{\log(1+\delta)}{2}\}\\
&\le2 \exp \{-\frac{\log(1+\delta)}{2d_1}m\}.
\end{align*}
Finally, noting that $\log(1+\delta) \ge \frac{\delta}{2}$, we have
$$\frac{\log(1+\delta)}{2d_1} \ge \frac{\Delta^2}{8\|\phi\|_{\infty}(1 + 3\|\kappa\|_{\infty}\mu(\XX))}.$$
The result follows. \qed \\

\subsection{Proof of Lemma \ref{theo:error-cont-norm}}
\label{sec:analysis-ker}

We begin with a general result for integral operators on general type spaces.
\begin{Lemma}
\label{lemma: change-measure}
Let $\nu, \mu$ be two mutually absolutely continuous finite measures on a measure space $(\XX, \TT)$.  Let
  $g=d\nu/d\mu$. Let $\kappa: \XX \times \XX \to \RRR_+$ be a kernel. Define
  another kernel $\kappa': \XX \times \XX \to \RRR_+$ as
  $$\kappa'(x,y)=\sqrt{\frac{g(x)}{g(y)}}\kappa(x,y), \; x,y \in \XX . $$
 Denote by $\KK$ [ resp. $\KK'$] the  integral operator associated with $\kappa$ [resp. $\kappa'$]
 on $L^2(\XX, \TT, \nu)$ [resp.  $L^2(\XX, \TT, \mu)$].  
  Then $\|\KK\|_{L^2(\nu)} = \|\KK'\|_{L^2(\mu)}$. \end{Lemma}

\textbf{Proof:} Note that the operator $\Acal: L^2(\XX,\TT,\nu) \to L^2(\XX, \TT,\mu)$ defined as
 $(\Acal f) =\sqrt{g} f$, $f \in L^2(\nu)$, is an isometry.  Also, for $f \in L^2(\mu)$
 $$( \Acal \KK \Acal ^{-1} f)(x)= \sqrt{g(x)} \int_\XX \kappa(x,y)\frac{1} {\sqrt{g(y)}} f(y) \mu(dy) = \int_\XX \kappa'(x,y)f(y)\mu(dy).$$
 Thus $\Acal \KK \Acal ^{-1} = \KK'$.  The result now follows on noting that 
 $\|\Acal \KK \Acal ^{-1}\|_{L^2(\mu)} = \|\KK\|_{L^2(\nu)}$.  \qed \\

For the rest of this subsection we will take $(\XX, \TT)$ to be the cluster space introduced in Section \ref{sec:model-irg} (see above 
\eqref{eqn:phi-def}).  Given rate functions $(a,b,c)$, $\mu_t(a,b,c), \kappa_t(a,b,c), \phi_t(a,b,c), \KK_t(a,b,c), \rho_t(a,b,c)$ are as
introduced above Lemma \ref{theo:error-cont-norm}.

\begin{Lemma}
\label{lemma:RN-derivative}
Let $(a_i, b_i, c_i)$, $i=1,2$ be two sets of rate functions.  
 Suppose that $a_1, c_1$ are strictly positive on $(0, T]$ and
$a_2 \le a_1, c_2 \le c_1$ on $[0, T]$.  Also suppose that for some $ \delta \in (0, e^{-T}/T)$, $c_1 \le c_2+\delta$ on $[0, T]$.
Fix $t \in [0, T]$. 
Let $\mu_i = \mu_t(a_i, b_i, c_i)$, $i=1,2$.  Then $\mu_2 \ll \mu_1$ and
\begin{equation*}
\frac{d\mu_2}{d\mu_1} (s,w)=\frac{a_2(s)}{a_1(s)} \times \exp  \{ -\int_{s}^{T} w(u)[ c_2(u)-c_1(u)]du \} \times \Pi_{i=1}^{w(T)-2} \left(\frac{c_2(\tau_i)}{c_1(\tau_i)} \right), \; (s,w) \in [0, t] \times \WW .
\end{equation*}
where $\tau_i=\tau_i(s, w)$ is ($\mu_1$ a.s.) the $i^{th}$ jump of $w$ after time $s$.
\end{Lemma}
\textbf{Proof:} 
Recall the probability measure $\nu_s$ on $\WW$ introduced below \eqref{eqn:phi-def}.  Write $\nu^i_s = \nu_s(a_i, b_i, c_i)$, $i=1,2$.
Note that (see \eqref{eqn:mu-def}), for $s \in [0, t]$,  $\mu_i(ds\, dw) = \nu^i_s(dw) a_i(s) ds$, $i=1,2$.  Thus to prove the result it
suffices to show that, for all $s \in [0, t]$, $\nu^2_s \ll \nu^1_s$ and
$\frac{d\nu^2_s}{d\nu^1_s} = L_s^T$, where, for $t \in [s, T]$,
\[
L_s^t(w)= \exp  \{ -\int_{s}^{t} w(u)[ c_2(u)-c_1(u)]du \} \times \Pi_{i \ge1}\left(\frac{c_2(\sigma_i)}{c_1(\sigma_i)}  {\bf1}_{\{\sigma_i \le t\}}\right),\]
and  $\sigma_i(w)$ is ($\nu_s^1$ a.s.) the $i^{th}$ jump of $w$.
For this, it suffices in turn to show that $\int_{\WW} L_s^t(w) \nu_s^1(dw) = 1$ for all $t \in [s, T]$ (see eg. Theorem T3, p.166 of \cite{Bremaud}).  The process $\{L_s^t\}_{ t \in [s, T]}$ on $\WW$ with the canonical filtration is a local martingale
under $\nu_s^1$ (see Theorem T2, p.166, \cite{Bremaud}) so to check the martingale property, it suffices to check (see (2.4) on page 166, and Theorem T8 on page 27 of \cite{Bremaud}), that
\[
\int_{\WW} [\int_s^T L_s^u(w) |c_1(u)-c_2(u)| du ]  \nu_s^1(dw)< + \infty.
\]
Note that $L_s^u(w)  \le \exp \{ T w(T) \delta\} $ since $c_1- \delta \le c_2 \le c_1$.   Thus
\begin{align*}
\int_{\WW} [\int_s^T L_s^u(w) |c_1(u)-c_2(u)| du ]  \nu_s^1(dw)
&\le 2T \int_{\WW} e^{ T w(T) \delta } \nu_0^1(dw).
\end{align*}
Note that, under $\nu_0^1$,   $w(T)$ is stochastically bounded by the sum of two independent Geom($e^{-T}$)  (see proof of  Lemma
 \ref{theo:error-wst}). Thus the last integral is bounded by
 $(\E e^{TZ\delta})^2$, where $Z$ is a Geom($e^{-T}$) random variable.  This expectation is finite since $T\delta < e^{-T}$.  The result follows. \qed \\

Proof of Lemma \ref{theo:error-cont-norm} will make use of parts (iv) and (v) of Lemma \ref{lemma:prop-kernel}.  In order to use (v) we will need the kernels
to be suitably bounded.  For that we use a truncation of the kernels, the associated error of which is estimated through the following lemma.\\

\begin{Lemma}
\label{lemma:kappa-trunc1}
  For $A \in (0, \infty)$, $t \in [0, T]$ and rate functions  $(a,b,c)$, define the kernel
$\kappa_{A,t}(a,b,c) \equiv \kappa_{A,t}: \XX \times \XX \to [0, \infty)$ as
$$\kappa_{A,t}(\bfx, \bfy) = \kappa_{t}(\bfx, \bfy) 1_{\{w(T) \le A, \tilde w(T) \le A\}},\;
\mbox{ where }
\kappa_t = \kappa_t(a,b,c), \bfx = (s,w), \bfy = (r, \tilde w).$$
Denote by $\KK_{A,t}$ the integral operator corresponding to $\kappa_{A,t}$ on $L^2(\XX, \TT, \mu_t)$ and $\rho(\kappa_{A,t})$ its norm, where $\mu_t = \mu_t(a,b,c)$.  Then there exist $A_0, C_3, C_4 \in (0, \infty)$ such that 
$$\rho(\kappa_t) - C_3e^{-C_4A} \le \rho(\kappa_{A,t}) \le \rho(\kappa_t),$$
for all rate functions $(a,b,c)$, $t \in [0, T]$ and $A \ge A_0$.
\end{Lemma}
\textbf{Proof:} We will suppress $t$ in the notation.  Since $\kappa_A \le \kappa$, from Lemma \ref{lemma:prop-kernel} (ii) and (iii) we have
 $$\rho(\kappa) - \rho (\kappa-\kappa_A) \le \rho(\kappa_A) \le \rho(\kappa).$$
 Consequently,
 $$\rho(\kappa) - \rho (\kappa_A) \le \rho (\kappa-\kappa_A) \le  \| \kappa-\kappa_A\|_2,$$
 where $\|\kappa\|_2$ denotes the $L^2(\mu \times \mu)$ norm of the kernel $\kappa$.  Note that for $\bfx, \bfy \in \XX$
\begin{equation}
\label{eqn:kappa-bound}
\kappa(\bfx,\bfy)=\int_0^t w(u)\tilde w(u) b(u)du \le  T w(T)\tilde w(T), \mbox{ for } \mu\times\mu \mbox{ a.e. } (\bfx, \bfy) = ((s,w), (r,\tilde w)).
\end{equation}
Let $\Lambda=\{ (s, w)\in \XX : w(T) > A\}$.  Then for fixed $\bfx=(s, w) \in \Lambda^c$
\begin{align*}
\int_\XX (\kappa(\bfx,\bfy)-\kappa_A(\bfx,\bfy))^2 \mu(d\bfy)
&\le \int_\XX (\kappa(\bfx,\bfy) {\bf1}_\Lambda (\bfy) )^2 \mu(d\bfy)\\
&\le\int_{\XX} [ T^2 w^2(T)\tilde w^2(T){\bf1}_{\{\tilde w(T) > A \}}] \mu(dr d\tilde w) )\\
&\le T^3 w^2(T) \E_0 [w^2(T){\bf1}_{\{w(T) > A \}}], 
\end{align*}
where in the second line we have used  \eqref{eqn:kappa-bound}  and $\E_0$ in the third line denotes the expectation 
corresponding to the probability measure $\nu_0$ on $\WW$.
Next, for fixed $\bfx=(s, w) \in \Lambda$,  we have in a similar manner
\begin{equation*}
\int_\XX (\kappa(\bfx,\bfy)-\kappa_A(\bfx,\bfy))^2 \mu(d\bfy) \le  T^3 w^2(T) \E_0 [w^2(T)].
\end{equation*}
Thus for any $\bfx=(s,w) \in \XX$,
\begin{equation*}
\int_\XX (\kappa(\bfx,\bfy)-\kappa_A(\bfx,\bfy))^2 \mu(d\bfy) \le T^3 w^2(T) \E_0 [w^2(T){\bf1}_{\{w(T) > A \}}] + {\bf1}_{\{w(T) > A \}} T^3 w^2(T) \E_0 [w^2(T)].
\end{equation*}
Integrating with respect to $\bfx \in \XX$, and noting that $\nu_s(w(T) \ge \alpha)  \le \nu_0(w(T) \ge \alpha)$, $\alpha \ge 0$,  we have
\begin{align*}
\|(\kappa-\kappa_A)\|_2^2
&=\int_\XX \int_\XX (\kappa(\bfx,\bfy)-\kappa_A(\bfx,\bfy))^2 \mu(d\bfy)  \mu(d\bfx)\\
&\le T^3 \E_0 [w^2(T){\bf1}_{\{w(T) > A \}}] \int_0^t a(s) \E_0 [ w^2(T) ]ds+ T^3 \E_0 [w^2(T)]\int_0^t a(s) \E_0[{\bf1}_{\{w(T) > A \}} w^2(T)] ds\\
&\le 2 T^4 \E_0 [w^2(T)] \E_0 [w^2(T){\bf1}_{\{w(T) > A \}}].
\end{align*}
As noted in the proof of  Lemma \ref{theo:error-wst}, under $\nu_0$, $w(T)$ is stochastically dominated by 
$ Z^*_1 +Z^*_2$
where $Z^*_1, Z^*_2$ are two independent copies of Geom($e^{-T}$). Therefore
$$ \E_0 [w_0^2(T)] \le \E (Z^*_1 +Z^*_2)^2 =2 e^{2T}(2-e^{-T})<4e^{2T}$$
and for a suitable $A_0 \in (0, \infty)$
\begin{equation*}
\E_0 [w_0^2(T){\bf1}_{\{w_0(T) > A \}}] \le \sum_{k>A} k^2 \times (k-1) e^{-2T}(1-e^{-T})^{k-2} \le (1-e^{-T})^{A/2}
\end{equation*}
for all $A \ge A_0$. Combining the above estimates, for all $A \ge A_0$
$$\|\kappa-\kappa_A\|_2 < 8 T^4 e^{2T}(1-e^{-T})^{A/2} .$$
The result follows. \qed \\

In the proof of Lemma \ref{theo:error-cont-norm} we will apply Lemma \ref{lemma: change-measure} to measures $\mu_1, \mu_2$ of the form in
Lemma \ref{lemma:RN-derivative}. This latter lemma shows that (under the conditions of the lemma) $\mu_2\ll\mu_1$.  However, to use 
Lemma  \ref{lemma: change-measure} we need the two measures to be mutually absolutely continuous.  To treat this difficulty we will use
an additional truncation  introduced in the lemma below and the elementary fact in Lemma \ref{elemz}.

\begin{Lemma}
\label{lemma:kappa-trunc2}
Given rate functions $(a,b,c)$ and $A \in (0, \infty)$ and $t \in [0, T]$, let $\kappa_t, \kappa_{A,t}$ be as in Lemma  \ref{lemma:kappa-trunc1}.
For $(s,w) \in \XX$, let $\tau(s, w) = \inf \{u > s: w(u)-w(u-) \neq 0\}$.  For $\delta > 0$, let
$\Lambda_{\delta} = \{(s,w) \in \XX: s \le \delta \mbox{ and } \tau(s, w) \le s+ \delta \}$. Define the kernel $\kappa_{A, \delta, t}$
as
$$
\kappa_{A, \delta, t}(\bfx, \bfy)=\kappa_{A,t}(\bfx,\bfy) {\bf 1}_{\Lambda_{\delta}^c}(\bfx){\bf 1}_{\Lambda_{\delta}^c}(\bfy).$$
Then there exists a $C_5 \in (0, \infty)$ such that for all rate functions $(a,b,c)$ and  $\delta, A \in (0, \infty)$ 
$$\rho(\kappa_{A,t}) -C_5 A^2 \delta \le \rho(\kappa_{A,\delta,t}) \le \rho(\kappa_{A,t}).$$
\end{Lemma}
\textbf{Proof:} Once more we will suppress $t$ from the notation.  As in the proof of  Lemma \ref{lemma:kappa-trunc1}, we have
$$\rho(\kappa_A) -\|\kappa_A-\kappa_{A,\delta}\|_2 \le \rho(\kappa_{A,\delta}) \le \rho(\kappa_A).$$
For $\bfx \in \Lambda^c$, 
\begin{align*}
\int_\XX ( \kappa_A(\bfx,\bfy)-\kappa_{A,\delta} (\bfx,\bfy))^2 \mu(d\bfy)
&\le \int_{\XX} \kappa^2_A(\bfx,(r,\tilde w)) {\bf1}_{\{\tau \le r+\delta \}} {\bf1}_{\{r \le \delta \}} \mu(dr d\tilde w)\\
&\le T^2A^4  \int_0^\delta \nu_r \{ \tau(r,w) \le r+ \delta \}  dr,
\end{align*}
where the above inequality uses the bound  $\kappa_A \le TA^2$.

Also, for fixed $\bfx \in \Lambda$,
$$\int_\XX ( \kappa_A(\bfx,\bfy)-\kappa_{A,\delta} (\bfx,\bfy))^2 \mu(d\bfy) \le T^2A^4 \cdot T$$
Thus for all $\bfx \in \XX$,
$$\int_\XX ( \kappa_A(\bfx,\bfy)-\kappa_{A,\delta} (\bfx,\bfy))^2 \mu(d\bfy) \le T^2A^4  \int_0^\delta \nu_r \{ \tau(r,w) \le r+ \delta \} dr 
+ {\bf1}_\Lambda(\bfx) T^3 A^4.$$
Finally
$$\|\kappa_A-\kappa_{A,\delta}\|_2^2 =
\int_ \XX \int_\XX ( \kappa_A(\bfx,\bfy)-\kappa_{A,\delta} (\bfx,\bfy))^2 \mu(d\bfy) \mu(d\bfx) \le 2T^3A^4 \int_0^\delta \nu_r \{ \tau(r,w) \le r+ \delta \} dr.$$
The result follows on observing that $\nu_r \{ \tau(r,w) \le r+ \delta \}=1-\exp\{-\int_r^{r+\delta} 2c(u)du \} \le 2 \delta$.
 \qed \\
 
 We will use the following elementary lemma.
 \begin{Lemma}
 \label{elemz}
Let $\gamma_0$, $\gamma_1$ be finite measures on a measure space $(\XX, \TT)$ such that $\gamma_0 \ll \gamma_1$. Let $G \in \TT$ be such that
$\{\bfx: d\gamma_0/ d\gamma_1 > 0 \} \supset G$.  For $i=0,1$, let $\gamma_i^G$ be restriction of $\gamma_i$   to $G$: $\gamma_i^G(\cdot)
= \gamma_i(\cdot \cap G)$. Then $\gamma_0^G$ and $\gamma_1^G$ are mutually absolutely continuous with
$\frac{d\gamma_0^G}{ d\gamma_1^G}(\bfx) = \frac{d\gamma_0}{ d\gamma_1}(\bfx) {\bf 1}_G(\bfx)$ and
$\frac{d\gamma_1^G}{ d\gamma_0^G}(\bfx) = (\frac{d\gamma_0}{d\gamma_1}(\bfx))^{-1} {\bf 1}_G(\bfx)$,
a.s. $\gamma_1^G$ and $\gamma_0^G$.
\end{Lemma} 

Lemma \ref{theo:error-cont-norm} requires establishing an estimate of the form $C_2(-\log \delta)^3 \delta^{1/2}$ for both
$|\rho_t- \rho_{+,t}|$ and $|\rho_t-\rho_{-,t}|$.  Proofs for the two cases are similar and so we only provide details for $|\rho_t-\rho_{-,t}|$ and leave the 
other case for the reader.  Given $\delta \in (0, \infty)$ and rate functions $(a,b,c)$ as in the statement of Lemma \ref{theo:error-cont-norm}, we denote
$(a_{\delta}, b_{\delta}, c_{\delta}) = ((a-\delta)^+, (b-\delta)^+, (c-\delta)^+)$.  Note that
\beqn
|\rho_t-\rho_{-,t}| &=& |\rho_t(a, b, c)-\rho_t(a_{\delta},b_{\delta},c_{\delta})|\nonumber \\
&
\le&  |\rho_t(a, b, c)-\rho_t(a_{\delta},b,c_{\delta})|+|\rho_t(a_{\delta}, b, c_{\delta})-\rho_t(a_{\delta},b_{\delta},c_{\delta})|.\label{eqn:only-use-once1}
\eeqn
We treat the first term on the right side in Lemma \ref{lemtria} while the second term is estimated
in Lemma \ref{lemtrib}.

\begin{Lemma}\label{lemtria}
	Let $(a,b,c)$ be rate functions as in the statement of Lemma \ref{theo:error-cont-norm}.
There exists $C_6, \delta_2 \in (0, \infty)$ such that for all $\delta \in (0,\delta_2)$ and $t \in [0, T]$
$$ |\rho_t(a,b,c)-\rho_t(a_{\delta},b,c_{\delta})| < C_6 (-\log \delta)^3 \delta^{1/2}.$$
\end{Lemma}
\textbf{Proof:} Since the kernel
$\kappa_t(a,b,c)$ does not depend on $a,c$, $\kappa_t(a,b,c) = \kappa_t(a_{\delta},b,c_{\delta})=\kappa_t$.  Henceforth suppress $t$ from the notation.
Let $\rho = \rho_t(a,b,c)$ and $\rho_{\delta} = \rho_t(a_{\delta},b,c_{\delta})$.  For $\eps > 0$, let $D \subset \XX$ be defined as
$$D\equiv D_{\eps} =\{(s,w): w(T)>A \mbox{ or } (\tau(s,w) \le s+\eps \mbox{ and } s \le \eps )\}^c$$
and define the kernel $\kappa_{D}$ as
$$\kappa_{D}(\bfx,\bfy)=\kappa(\bfx,\bfy) {\bf1}_{D}(\bfx){\bf1}_{D}(\bfy).$$
Using Lemmas \ref{lemma:kappa-trunc1} and Lemma \ref{lemma:kappa-trunc2},
$$ \rho-C_3 e^{-C_4 A}-C_5 A^2 \epsilon \le \rho(\kappa_D) \le \rho$$
and
$$ \rho_{\delta} -C_3 e^{-C_4 A}-C_5 A^2 \epsilon \le \rho_{\delta}(\kappa_D) \le \rho_{\delta}, $$
where $\rho(\kappa_D)$ [resp. $\rho_{\delta}(\kappa_D)$] is the norm of the corresponding integral operator on $L^2(\mu)$ [resp. $L^2(\mu_\delta)$],
where $\mu = \mu(a,b,c)$ and $\mu_{\delta} = \mu(a_{\delta},b,c_{\delta})$.
Thus we have
\begin{equation}
\label{eqn:only-use-once}
|\rho-\rho_\delta| < 2C_3 e^{-C_4 A} + 2C_5 A^2 \epsilon + |\rho(\kappa_D)-\rho_\delta(\kappa_D)|.
\end{equation}
We now estimate $|\rho(\kappa_D)-\rho_\delta(\kappa_D)|$. 
By Lemma \ref{lemma:RN-derivative}, $\mu_\delta \ll \mu$ and 
\begin{equation*}
g (s,w)=_{\sss def}\frac{d\mu_\delta}{d\mu}(s,w)=\frac{a_\delta(s)}{a_0(s)} \times \exp  \{ -\int_{s}^{T} w(u)[ c_\delta(u)-c(u)]du \} \times \Pi_{i=1}^{w(T)-2} \left(\frac{c_1(\tau_i)}{c_0(\tau_i)} \right),
\end{equation*}
$\mu$ a.s., where $\tau_i$ are as in the statement of Lemma \ref{lemma:RN-derivative}.
For $\mu$ a.e. $(s,w) \in D$ we have 
$$w(t) \le w(T) \le A  \mbox{ and } \tau_1(s,w) > \eps, $$ 
consequently $c_{\delta}(\tau_i) \ge (\theta \eps - \delta)^+$ for all $i$, where $\theta$ is as in the statement of Lemma
\ref{theo:error-cont-norm}.  Also, since $a$ is bounded away from $0$, we can find $d_1 \in (0, \infty)$ such that
$a_{\delta}(s) \ge (d_1- \delta)^+$.  Thus for $\mu$ a.e. $\bfx \in D$
\begin{equation}\label{ins2144}
 (d_1- \delta)^+  \left((\theta \eps - \delta)^+\right)^A< g(\bfx) < \exp\{ T A \delta \}.
\end{equation}
Denote by $\mu^D$ [resp. $\mu_{\delta}^D$] the restrictions of $\mu$ [resp. $\mu_{\delta}$] to $D$.  Then from Lemma \ref{elemz}, whenever $\delta
< \delta_0 = \min \{\theta \eps , d_1\}$, $\mu^D$ and $\mu_{\delta}^D$ are mutually absolutely continuous and
$$\frac{d\mu_\delta^D}{d\mu^D}(\bfx) = g(\bfx) {\bf 1}_D(\bfx), \; \mbox{ a.e. } d\mu^D .
$$
For the rest of the proof we consider only $\delta < \delta_0$.  Then,
\begin{equation}\label{ins1119}
 (1- \frac{\delta}{d_1}) \left(1 - \frac{\delta}{\theta \eps}\right)^A< g(\bfx) < \exp\{ T A \delta \}.
\end{equation}

Note that $\rho(\kappa_D, \mu) = \rho(\kappa_D, \mu^D)$ and $\rho(\kappa_D, \mu_{\delta}) = \rho(\kappa_D, \mu^D_{\delta})$.
Also by Lemma \ref{lemma: change-measure}, $\rho(\kappa_D, \mu_{\delta}^D) = \rho(\kappa'_D, \mu^D)$, where
$$
\kappa'_D(\bfx, \bfy) = \kappa_D(\bfx , \bfy) \sqrt{\frac{g(\bfx)}{g(\bfy)}}{\bf 1}_D(\bfx){\bf 1}_D(\bfy), \; \bfx, \bfy \in \XX .
$$
Thus using Lemma \ref{lemma:prop-kernel}
\beq |\rho(\kappa_D)-\rho_\delta(\kappa_D)| &=&|\rho(\kappa_D, \mu)-\rho(\kappa_D, \mu_\delta)|\\
&=& |\rho(\kappa_D, \mu^D)-\rho(\kappa_D, \mu_\delta^D)|\\
&=&  |\rho(\kappa_D, \mu^D)-\rho(\kappa'_D, \mu^D)|\\
&\le& \rho(| \sqrt{\frac{g(\bfx)}{g(\bfy)}}-1|\kappa_D(\bfx,\bfy), \mu^D)  \\
&\le& T\sup_{\bfx,\bfy \in \XX}\left(| \sqrt{\frac{g(\bfx)}{g(\bfy)}}-1|\kappa_D(\bfx,\bfy)\right).\eeq
From \eqref{ins1119} we see that whenever $\delta \le \delta_0$,
\begin{equation*}
\sup_{\bfx,\bfy \in D}\left( \sqrt{\frac{g(\bfx)}{g(\bfy)} }\vee \sqrt{\frac{g(\bfy)}{g(\bfx)}}  \right)
< (1+\frac{2}{d_1} \delta) \times \exp\{ T A \delta \} \times (1+ \frac{2\delta}{\theta \epsilon})^A \equiv d(\delta, \eps, A).
\end{equation*}
Noting that $\kappa_D \le TA^2$, we have
\[
\sup_{\bfx,\bfy \in \XX} \left( | \sqrt{\frac{g(\bfx)}{g(\bfy)}}-1| \kappa_D(\bfx,\bfy)\right)  \le TA^2 |d(\delta, \eps, A)-1|.
\]
Thus
$$
|\rho(\kappa_D) - \rho_{\delta}(\kappa_D)| \le T^2 A^2 |d(\delta, \eps, A)-1|.
$$
Combining this with \eqref{eqn:only-use-once}, we have
\be
\label{ins1144}
|\rho-\rho_\delta| < 2C_3e^{-C_4 A} + 2C_5A^2 \epsilon + TA^2 |d(\delta, \eps, A)-1|.
\ee
Take $A= -\frac{1}{C_4} \log \delta$ and $\epsilon = \delta^{1/2}$.  Note that when $\delta$ is sufficiently small, $\delta \le \frac{1}{2} \min\{\theta \eps , d_1\}$
and so the above inequality holds for such $\delta$.
Also, with this choice, we can find $\delta_1 , d_2 \in (0, \infty)$ such that the sum of the first two terms  on the right side in \eqref{ins1144} is bounded by
$$ 2C_3\delta + \frac{2C_5}{C_4^2} (-\log \delta)^2 \delta^{1/2} \le d_2(-\log \delta)^2\delta^{1/2}  \; \mbox{ for all } \delta \le \delta_1.$$
Also note that with the above choice of $\eps, A$, $d(\delta, \eps, A) \to 1$ as $\delta \to 0$.  Furthermore
$$
d(\delta, \eps, A)  = (1+ O(\delta))(1+O(\delta)) (1+ O(\delta^{1/2} (-\log \delta))).$$
Thus,we can find $d_3 , \delta_2 \in(0, \infty)$ such that whenever $\delta \le \delta_2$
$$
 TA^2 |d(\delta, \eps, A)-1| \le d_3(-\log \delta)^3\delta^{1/2}.$$
 The result follows on combining the above estimates. \qed \\
We now estimate the second term in \eqref{eqn:only-use-once1}.
\begin{Lemma}
\label{lemtrib}
There exists $C_7 \in (0, \infty)$ and $\delta_3 \in (0,1)$ such that for all $t \in [0, T]$
$$ |\rho_t(a_{\delta},b,c_{\delta})-\rho_t(a_{\delta},b_{\delta},c_{\delta})| < C_7 (-\log \delta)^2 \delta$$
\end{Lemma}
\textbf{Proof:}
Denote $\kappa_t(a_{\delta}, b_{\delta},c_{\delta}) = \kappa_{\delta,t}=\kappa_{\delta}$ and recall that $\kappa_t(a_{\delta},b_{\delta},c_{\delta})=\kappa_t=\kappa$. We suppress $t$ in rest of the proof. For $A > 0$, let $\kappa_{A}$ [resp. $\kappa_{\delta, A}$] be the
truncated kernels as in Lemma \ref{lemma:kappa-trunc1} associated with $\kappa$ [resp. $\kappa_{\delta}$]. 
Clearly
$$|\kappa_{A}-\kappa_{\delta, A}| \le TA^2\delta.$$
Using this bound and Lemma \ref{lemma:kappa-trunc1}, for all $A \ge A_0$
\beq
|\rho(a_{\delta},b,c_{\delta})-\rho(a_{\delta},b_{\delta},c_{\delta})| &=& |\rho(\kappa, \mu_{\delta}) - \rho(\kappa_{\delta}, \mu_{\delta})|\\
&\le & TA^2\delta + 
 |\rho(\kappa, \mu_{\delta}) - \rho(\kappa_{A}, \mu_{\delta})|+ |\rho(\kappa_{\delta}, \mu_{\delta}) - \rho(\kappa_{\delta, A}, \mu_{\delta})|\\
 & \le & TA^2\delta + 2C_3 \exp\{-C_4A\},
 \eeq
where $\mu_{\delta} = \mu(a_{\delta}, b_{\delta}, c_{\delta}) =  \mu(a_{\delta}, b, c_{\delta})$ is as in Lemma \ref {lemtria}.

Result follows on taking $A = (-\log \delta)/C_4$ and $\delta$ sufficiently small.  \qed\\ \ \\

{\bf Proof of Lemma \ref{theo:error-cont-norm}.}
The proof is immediate on using Lemmas \ref{lemtria} and \ref{lemtrib} in \eqref{eqn:only-use-once1}.
\qed
 \\

\subsection{Proof of Lemma \ref{theo:error-cont-int}}
\label{sec:measure-theoretic}
We begin with an elementary lemma which allows one to regard the operators $\KK_t(a,b,c)$, $t \in [0, T]$, to be defined on a common Hilbert space.
Recall that for a kernel $\kappa$ on $\XX \times \XX$ and a finite measure $\mu$ on $\XX$, we denote by
$\rho(\kappa, \mu)$ the norm of the integral operator associated with $\kappa$ on $L^2(\XX, \TT, \mu)= L^2(\mu)$.
\begin{Lemma}
\label{lemcom}
Let $(a,b,c)$ be rate functions.  Then for all $t \in [0, T]$
$$\rho(\kappa_t(a,b,c), \mu_t(a,b,c))=\rho(\kappa_t(a,b,c),\mu_T(a,b,c)).$$
\end{Lemma}
\textbf{Proof:} Write $\kappa_t(a,b,c) = \kappa_t$, $\mu_t(a,b,c)= \mu_t$.  Denote the integral operator corresponding to
$\kappa_t$ on $L^2(\mu_t)$ [resp. $L^2(\mu_T)$] by $\KK_t$ [resp. $\KK_t^T$].
Let $\XX_t=[0,t] \times \WW$, then $\mu_t$ is supported on $\XX_t$ and $\kappa_t$ is supported on $\XX_t \times \XX_t$. Thus for any $\psi \in L^2(\mu_T)$, $\KK_t^T \psi=\KK_t \psi$ is also supported on $\XX_t$, and this implies $\|\KK_t^T\| \le \|\KK_t\|$. On the other hand, for any $\psi \in L^2(\mu_t)$,
$\KK_t \psi = \KK_t^T (\psi{\bf 1}_{[0, t]\times \WW})$
$$\|\KK_t \psi\|_{L^2(\mu_t)} \le \|\KK_t^T\| \|\psi{\bf 1}_{[0, t]\times \WW}\|_{L^2(\mu_T)} = \|\KK_t^T\| \|\psi\|_{L^2(\mu_t)}.$$
Thus $\|\KK_t^T\| \ge \|\KK_t\| $. \qed \\

The following theorem concerning IRG models on a general type space is a corollary of Theorem 3.1 and Theorem 3.12 in \cite{bollobas-riordan-janson}. 
\begin{Theorem}\cite{bollobas-riordan-janson}
\label{theo:bollobas}
Let $(\XX, \TT, \mu)$ be a type space.
Consider the weight function  $\phi \equiv 1$ on this space. Let $\kappa_n, \kappa$ be symmetric kernels on $\XX \times \XX$
 such that
$$ x_n \to x \mbox{ and } y_n \to y \mbox{ implies } \kappa_n(x_n,y_n) \to \kappa(x,y).$$
Let $\CC^{(1)}_n(\kappa)$ denote  the size of the largest component in $\RG_n(\kappa, \mu, \phi)$. Then\\
(i) If $\rho(\kappa,\mu) \le 1$, then $\CC^{(1)}_n(\kappa_n)/n \convp 0$. Furthermore, if $\|\kappa\|_{\infty} < +\infty$, then $\CC^{(1)}_n(\kappa_n) = O (\log n)$.\\
(ii)If $\rho(\kappa,\mu) > 1$, then $\CC^{(1)}_n(\kappa_n) = \Theta (n)$.
\end{Theorem}

\textbf{Proof of Lemma \ref{theo:error-cont-int}:} 
Let $\kappa_t(a_0, b_0, c_0) = \kappa_t$ and $\mu_t(a_0, b_0, c_0) = \mu_t$.  Note that for $0 \le t_1 \le t_2 \le T$, since $\kappa_{t_1} \le
\kappa_{t_2}$, we have from Lemma \ref{lemcom} that 
$$\rho(t_1) = \rho(\kappa_{t_1}, \mu_{t_1}) = \rho(\kappa_{t_1}, \mu_{T}) \le \rho(\kappa_{t_2}, \mu_{T}) = \rho(\kappa_{t_2}, \mu_{t_2})=\rho(t_2).$$
Thus $\rho(t)$ is nondecreasing in $t$.  Next note that, since $w(u)\tilde w(u) b_0(u)$ is non decreasing in $u$ for $\mu_T \times \mu_T$ a.e.
$((s,w), (r, \tilde w))$, $\kappa_t(\bfx, \bfy)$ is convex in $t$ for a.e. $\bfx, \bfy$, i.e. for $\mu_T \times \mu_T$ a.e. $(\bfx , \bfy) \in \XX \times \XX$,
and all $t_1, t_2 \in [0, T]$, and $\alpha, \beta \in [0,1]$, $\alpha+\beta = 1$,
$$\kappa_{\alpha t_1+ \beta t_2}(\bfx, \bfy) \le \alpha \kappa_{t_1}(\bfx, \bfy) +\beta \kappa_{t_2}(\bfx, \bfy).$$
Thus
\begin{align*}
\rho(\alpha t_1 +\beta t_2)
&=\rho( \kappa_{\alpha t_1 +\beta t_2}, \mu_T) \\
&\le \rho(\alpha \kappa_{t_1}+\beta \kappa_{t_2}, \mu_T) \\
&\le \rho(\alpha \kappa_{t_1}, \mu_T)+\rho(\beta \kappa_{t_2}, \mu_T) \\
&= \alpha \rho(\kappa_{t_1}, \mu_T)+\beta \rho(\kappa_{t_2}, \mu_T)  \\
&=\alpha\rho( t_1) +\beta \rho(t_2),
\end{align*}
where lines 3 and 4 above use parts  (ii) and (iii) of Lemma \ref{lemma:prop-kernel} and line 2
uses the convexity of $\kappa_{\cdot}$.

Thus  $\rho$ is convex on $[0, T]$.  Also since $\rho(0) = 0$ and $\rho(t) > 0$ for $t > 0$, we have that
$\rho$ is strictly increasing on $[0, T]$ and has a strictly positive left derivative on $(0, T]$.  This proves parts (i) and (iii) of Lemma \ref{theo:error-cont-int}.

We now consider part (ii).  For $\delta > 0$, let
$$\rho^{\delta,+}(t)=\rho_t((a_0+\delta)\wedge 1,(b_0+\delta)\wedge 1,(c_0+\delta)\wedge 1), \;
\rho^{\delta,-}(t)=\rho_t((a_0-\delta)^+,(b_0-\delta)^+,(c_0-\delta)^+).$$
Similarly define $\mu_t^{\delta, -}$, $\kappa_t^{\delta, -}$,  $\mu_t^{\delta, +}$ and $\kappa_t^{\delta, +}$.
We will argue by contradiction.\\

Suppose first that  $\rho(t_c) > 1$, then by Lemma \ref{theo:error-cont-norm} and the continuiity of $\rho(t)$, there exist $\epsilon, \delta >0$ 
such that $\rho^{\delta,-}(t_c-\epsilon)>1$. 
Denote $\RG^{\delta, -}_t(\kappa_n) = \RG_n(\kappa^{\delta, -}_{n,t}, \mu_t^{\delta,-}, \phi_t)$,
where $\kappa^{\delta, -}_{n,t}$ is defined as in \eqref{eqn:kappan-def}, replacing $b$ there with $(b_0-\delta)^+$.
Since $\kappa^{\delta, -}_{n,t}$ converges to  $\kappa^{\delta, -}_{t}$ uniformly on compact subsets of $\XX \times \XX$, by Theorem \ref{theo:bollobas}, the size of the largest component of $\RG^{\delta, -}_t(\kappa_n)$, whp,  is $\Theta(n)$ and consequently the volume of 
the largest component of $\RG^{\delta, -}_t(\kappa_n)$ is, whp,  at least $\Theta(n)$.
 By Lemma \ref{lemma:couple-bfia}, and as $\barx_n(t) \to x(t)$, we have whp
$$\com_n(t_c-\epsilon)  \ge_d \IA_n((a_0-\delta)^+,(b_0-\delta)^+,(c_0-\delta)^+)_{t_c-\epsilon}$$
and from Lemma \ref{lemma:rgiva-irg-def} the largest component in
$\IA_n((a_0-\delta)^+,(b_0-\delta)^+,(c_0-\delta)^+)_{t_c-\epsilon}$ has the same distribution as that in 
 $\RG^{\delta, -}_{t_c-\epsilon}(\kappa_n)$.
However, by Theorem 1.1 of \cite{spencer2007birth} the largest component size of Bohman-Frieze model for $t < t_c$ is,
whp,   $\Theta( \log n)$ which contradicts the fact that the volume of 
the largest component of $\RG^{\delta, -}_{t_c-\epsilon}(\kappa_n)$ is, whp,  at least $\Theta(n)$.  Thus we have shown that $\rho(t_c) \le 1$.

Suppose now that $\rho(t_c) < 1$. then there exists $\epsilon, \delta >0$ such that $\rho^{\delta,+}(t_c+\epsilon) <1$.  Then a similar argument as above shows   that, whp,
$$\com_n(t_c+\epsilon) \le_d \IA_n((a_0+\delta)\wedge 1,(b_0+\delta)\wedge 1,(c_0+\delta)\wedge 1))_{t_c+\epsilon}=_{def}  \IA_{n,t_c+\epsilon} ^{\delta}.$$
 Lemma \ref{lemma:logn-bound}  implies that whp the largest component in $\IA_{n,t_c+\epsilon}^\delta$ is $O(\log^4 n)$. However from Theorem 1.1 of \cite{spencer2007birth}, for $t > t_c$, the largest component size of Bohman-Frieze model is whp $\Theta(n)$.
 This contradiction shows that $\rho(t_c) \ge 1$.  Combining the above arguments  we  have $\rho(t_c)=1$. \qed\\

\section{Proof of Proposition \ref{prop:main}}
\label{sec:analysis-s2s3}

We shall now study the asymptotics of $\calS_2$ and $\calS_3$, namely the sum of squares and cubes of the component sizes. We first analyze  $\calS_2$ since the asymptotics for this will be required in the analysis of $\calS_3$. 

\subsection{Analysis of $\bars_2(\cdot)$ near criticality}
\label{s2sec}
Let us start with the sum of squares.  Recall from  \eqref{ins2141} that $\calS_2(t)$ denotes the sum of squares of the component sizes  in $\BF(t)$ and $\bars_2(t) = \calS_2(t)/n$. Also recall the limiting functions $s_k(t)$, $k=2,3$, introduced in \eqref{ins-s2} and \eqref{ins-s3}. Note that these functions are non-decreasing 
and they blow up at the critical point $t_c$ (see \eqref{eqn:s2-scaling-crit} and \eqref{eqn:s3-scaling-crit}).
Let $y(t) = 1/s_2(t)$. Then this function satisfies the differential equation 
\begin{equation}
	\label{eqn:yt-diff}
	y^\prime(t) = -x^2(t) y^2(t) - (1-x^2(t)), \; y(0)= 1, \; t \in [0, t_c].
\end{equation}
  Note
 that $y$ is a monotonically decreasing function with $y(t)\to 0$ as $t\to t_c$ and  as shown in Theorem 3.2 of \cite{janson2010phase}, the scaling behavior near $t_c$ is 
\begin{equation}
	y(t) = \frac{1}{\alpha}(t_c - t)+  O(t_c-t)^2
	\label{eqn:yt-scaling}
\end{equation} 
as $t \uparrow t_c$, where $\alpha$ is as in \eqref{eqn:alpha-def}. Let $Y_n(t) = 1/\bars_2(t)$. To simplify notation, we  suppress the dependence of the process $Y$ on $n$
when convenient. Note that \eqref{eqn:s2-tn-n-alpha} is equivalent to showing 
\begin{equation}
	n^{1/3}\left|Y(t_n)-\frac{1}{\alpha n^\gampar}\right|\convp 0.
	\label{eqn:ynt-to-show}
\end{equation}
Here, and throughout Section \ref{sec:analysis-s2s3},  $t_n=t_c - 1/n^\gampar$ with $\gampar\in (1/6,1/5)$. From  \eqref{eqn:yt-scaling} 
\be \label{ins1607}\left|y(t_n) - \frac{1}{\alpha n^\gampar}\right| = O\left(\frac{1}{n^{2\gampar}}\right) =o\left(\frac{1}{n^{1/3}}\right).\ee
Thus to show \eqref{eqn:ynt-to-show} it is enough to prove the following
\begin{Proposition}
	\label{prop:s2-crit}
	As $n\to\infty$
	\[ n^{1/3}\sup_{s\leq t_n}\left|Y(s) - y(s)\right| \convp  0. \]

\end{Proposition} 

  We shall prove this via a sequence of lemmas. We begin with some notation. 
  Recall that $\bfC_n^{\sss BF}(t) \equiv (\CC_n^{\sss (i)}(t) : i \ge 1)\equiv (\CC_i(t) : i \ge 1)$ is the component size vector, $I_n(t)$  the size of the largest component, and
$X_n(t)$ the number of singletons, in  $\BF_n(t)$.   Let $\sum_i$ denote the summation over all components and $\sum_{i<j}$ denote the summation over all pairs of components $(i,j)$ with $i<j$.

The first Lemma identifies the semimartingale decomposition  of the process $Y(\cdot)$ as well as the predictable quadratic variation  $\langle M\rangle$ of the martingale $M$ in the decomposition.  Recall the natural filtration associated with the BF process introduced in Section \ref{sec:model-equiv}.

\begin{Lemma}
	\label{lemma:mart-s2}
	The process $Y_n(\cdot)$ can be decomposed as 
	\be Y_n(t) = 1+ \int_0^t A_n(s) ds + M_n(t), \; t \in [0, t_c] \label{ins1006}\ee
	where 
	\\(a)  $M_n$ is a RCLL martingale with respect to the natural filtration $\{\FF_t\}_{t\geq 0}$ of the BF process. 
	\\(b) The process $A_n = A_1^n+R_1^n$ where (suppressing $n$)
	\begin{equation}
		\label{eqn:gnu}
		A_1(u)  =  -Y^2(u)\barx^2(u) -(1-\barx^2(u)) + (1-\barx^2(u))\frac{Y^2(u)}{n^2}\sum_i \calC_i^4(u), \;\;u \le t_c
	\end{equation}
	and for some  $C_7 \in (0, \infty)$,
\begin{equation}
	|R_1^n(u)| \le C_7 \left( \frac{1}{n} + \frac{I_n^2(u)}{n} \right), \mbox{ for all } n \in \NNN \mbox{ and } u \le t_c.  \label{eqn:r1s}
\end{equation}
	\\(c) Predictable quadratic variation of $M_n$ is given as
	\[\langle M_n\rangle(t) = \int_0^t B_n(u) du \]
	where  $B_n$ is such that
	\begin{equation}
		B_n(u)\le \frac{4}{n} + \frac{4 Y_n^2(u) I_n^2(u)}{n}, \; u \le t_c.
		\label{eqn:yi-bound}
	\end{equation}
\end{Lemma} 
{\bf Proof:}  
We will suppress $n$ from the notation when convenient.  Note that 
$$Y(t) = 1 + \sum_{s \le t} \Delta Y(s), \mbox{ where } \Delta Y(s) = Y(s)- Y(s-).$$
We now analyze the possible jumps of $Y$.  Note than any jump in $Y$ corresponds to a jump of one of the Poisson processes $\PP_\bfe$,
$\bfe = (e_1,e_2) \in \EE^2$ (recall the notation from Section \ref{sec:cont-time-bf}).   A jump of $\PP_\bfe$ at a time instant $u$ could result in two different kinds of
jumps in $Y$.\\
%
%
%
%
%
%
%
%
%

(i) {\bf Merger caused by the first edge $e_1$:} In this case $\Delta \calS_2(u)= \calS_2(u) - \calS_2(u-) =2$ which implies
$$\Delta Y(u) \equiv \alpha_1(u-) = -\frac{2Y^2(u-)}{n}\left[ 1-O\left( \frac{2Y(u-)}{n}\right)\right] .$$
(ii) {\bf Merger caused by the second edge $e_2$:} In this case, suppose components $i$ and $j$ merge, then
 $\Delta \calS_2(u)=2\CC_i(u-) \CC_j(u-)$ and thus
\[\Delta Y_n(u) \equiv \alpha_2^{i,j}(u-)= -2\frac{\calC_i(u-)\calC_j(u-)}{n}Y_n^2(u-)\left[1- O\left(2\frac{\calC_i(u-)\calC_j(u-) Y_n(u-)}{n}\right)\right].\]
With these observations we can represent $Y$ in terms of stochastic integrals with respect to $\PP_\bfe$ as follows.  Define
$$
\clh_1(u) = \{\bfe = (e_1, e_2) \in \EE^2: e_1 = (v_1, v_2) \mbox{ where both } v_1, v_2 \mbox{ are singletons at time } u\},$$
 
 $$
\clh_2^{(i,j)}(u) = \{\bfe = (e_1, e_2) \in \EE^2\setminus \clh_1(u): e_2 = (v_1, v_2) \mbox{ where one vertex is in } \clc^i(u) \mbox{ while the other is in } \clc^j(u)  \}.$$
Also let
$$\clu_{\bfe}(u) = \alpha_1(u) {\bf 1}_{\clh_1(u)}(\bfe), \; \clu_{\bfe}^{i,j}(u) = \alpha_2^{i,j}(u){\bf 1}_{\clh_2^{(i,j)}(u)}(\bfe).$$
Then
\be \label{ins2023}
Y(t) = 1 + \sum_{\bfe \in \EE^2} \int_{(0, t]} \left (\clu_{\bfe}(s-) + \sum_{i< j} \clu_{\bfe}^{i,j}(s-)\right) \PP_{\bfe}(ds).\ee
Recalling that  $\PP_{\bfe}$ is a rate $2/n^3$ Poisson process, one can write $Y$ as
$$Y(t) = 1 + \int_{[0,t]} A(s) ds + M(t),$$
where
$$
A(s) = \frac{2}{n^3} \sum_{\bfe \in \EE^2}\left (\clu_{\bfe}(s) + \sum_{i< j} \clu_{\bfe}^{i,j}(s)\right).$$
Note that
\beqn
\sum_{\bfe \in \EE^2}{\bf 1}_{\clh_1(s)}(\bfe) &=& {X_n(s) \choose 2}{n \choose 2}=\frac{n^4}{4}\barx_n^2(s) \cdot (1+O(1/n)),\label{ins2029a}\\
\sum_{\bfe \in \EE^2}{\bf 1}_{\clh_2^{(i,j)}(s)}(\bfe) & = &
  \left[{n\choose 2}-{X_n(s)\choose 2}\right] \calC_i(s)\calC_j(s) = \frac{n^2}{2}(1-\barx_n^2(s)) \calC_i(s)\calC_j(s) \cdot (1+O(1/n)).\label{ins2029b}
  \eeqn
Thus we get 
$A(s) = A_1(s) + R_1(s)$, 
where $A_1$ represents the leading order terms:
\begin{align*}
	A_1(s) &= - \frac{n}{2}\barx_n^2(s) \cdot \frac{2Y^2(s)}{n} -\sum_{i<j} \frac{1}{n}(1-\barx_n^2(s))\CC_i(s) \CC_j (s)\cdot \frac{2Y^2(s) \CC_i(s) \CC_j(s)}{n}\\
	&= -\barx_n^2(s) Y^2(s)-(1-\barx_n^2(s)) Y^2(s)\cdot \frac{1}{n^2}\sum_{i<j} 2\CC_i^2(s)\CC_j^2(s)\\
	&= -\barx_n^2(s) Y^2(s)-(1-\barx_n^2(s)) Y^2(s)\cdot \frac{1}{n^2}[(\sum_{i} \CC_i^2(s))^2-\sum_i \CC_i^4(s)]\\
	&= -\barx_n^2(s) Y^2(s)-(1-\barx_n^2(s)) +(1-\barx_n^2(s)) Y^2(s)\cdot \frac{1}{n^2}\sum_i \CC_i^4(s)
\end{align*}
and the last line follows from the fact $Y=\frac{n}{\sum_i \CC_i^2}$. 
The term $R_1$ consists of the lower order terms and using the observations $Y \le 1$, $\barx_n \le 1$ and $|A_1| \le \barx_n^2 Y^2+(1-\barx_n^2) \le 2$, it can be
estimated as follows.
\begin{align*}
	|R_1(u)| 
	&\le |A_1(u)| \cdot \frac{d_1}{n} + \frac{n}{2}\barx_n^2(u) \cdot \frac{2Y^2(u)}{n} \cdot \frac{d_2 Y(u)}{n}\\
	&+ (1-\barx_n^2(u)) Y^2(u)\cdot \frac{1}{n^2}\sum_{i<j} 
	[2\CC_i^2(u)\CC_j^2(u) \cdot \frac{d_3 \CC_i(u)\CC_j(u) Y}{n}(u)]\\
	&\le \frac{2d_1}{n} + \frac{d_2}{n} + \frac{d_3 I^2(u)}{n} \cdot \frac{Y^2(u)}{n^2} \sum_{i<j} 2\CC_i^2(u)\CC_j^2(u)\\
	&\le \frac{2d_1+d_2}{n} + \frac{d_3 I^2(u)}{n}.
\end{align*}
Next, from \eqref{ins2023} and using independence of Poisson processes $\PP_{\bfe}$,
$$
\langle M\rangle(t) = \frac{2}{n^3} \sum_{\bfe \in \EE^2} \int_{(0, t]} \left ((\clu_{\bfe}(s))^2 + \sum_{i< j} (\clu_{\bfe}^{i,j}(s))^2 \right) ds
\equiv \int_{(0,t]} B(s) ds,$$
 where, using \eqref{ins2029a} and \eqref{ins2029b} once more, we can estimate $B$ as follows.
\begin{align*}
	B(u) 
	&\le \frac{n}{2}\barx_n^2(u) \cdot 2 \cdot \left(\frac{2 Y^2(u)}{n}\right)^2 + (1-\barx_n^2(u)) \cdot 2 \cdot \sum_{i<j} \frac{\CC_i(u)\CC_j(u)}{n} \left(\frac{2Y^2(u)\CC_i(u)\CC_j(u)}{n}\right)^2\\
	&\le \frac{4}{n} + \frac{4 Y^2(u) I^2(u)}{n} \cdot \frac{Y^2(u)}{n^2}\sum_{i<j} 2\CC_i^2(u)\CC_j^2(u)\\
	&\le \frac{4}{n} + \frac{4 Y^2(u) I^2(u)}{n}.
\end{align*}
This completes the proof of the lemma.
 
\qed

The following result  bounds the difference  $Y-y$ through an application of Gronwall's lemma.
\begin{Lemma}
	\label{lemma:gronwall}
	For every $n \in \NNN$,
	\[\sup_{s\leq t_n}|Y_n(s)-y(s)|\leq \eps_n e^{2t_c}\]
	where 
	\[\eps_n= 4t_c\sup_{s\leq t_c}|\barx(s)-x(s)|+\sup_{s\leq t_n}|M_n(t)|+C_8\int_0^{t_n} \frac{I^2_n(s)}{n}ds,\]
	$C_8 = (2C_7+1)$ and $C_7$, $M_n$ are as in Lemma \ref{lemma:mart-s2}.
\end{Lemma}
{\bf Proof:} Since $y$ solves \eqref{eqn:yt-diff} and $Y$ satisfies \eqref{ins1006}, using the fact that $\barx, x$ take values in $[0,1]$, we have for any fixed $t\leq t_n$
\begin{align}
	|Y(t)-y(t)| \leq\left|\int_0^t \left (Y^2(s)\barx^2(s)-y^2(s)x^2(s)\right) ds \right|&+ \int_0^t|\barx^2(s)-x^2(s)|ds +\int_0^t |R_1(s)|ds \nonumber \\
	&+\int_0^t \left(\frac{Y^2(s)}{n^2} \sum_{i}\calC_i^4(s) \right)ds +\sup_{s\leq t_n} |M(t_n)|. \label{ins1012}
\end{align}
Let us now analyze each term individually. Writing 
\[Y^2(s)\barx^2(s)-y^2(s)x^2(s) =\barx^2(s)(Y^2(s)-y^2(s)) + y^2(s)(\barx^2(s)-x^2(s))\]
and using the fact that $Y,y,x, \bar x$ take values in the interval $[0,1]$ we get 
\[\left|\int_0^t Y^2(s)\barx^2(s)-y^2(s)x^2(s) ds \right|\leq 2\int_0^t|Y(s)-y(s)|ds+ 2t_c\sup_{s\leq t_c}|\barx(s)-x(s)|\]
Similarly the second term in \eqref{ins1012} can be estimated as
\[\int_0^t|\barx^2(s)-x^2(s)|ds \leq 2t_c\sup_{s\leq t_c}|\barx(s)-x(s)|. \]
From \eqref{eqn:r1s} the integrand in the third term in \eqref{ins1012} can be bounded as 
$$
	|R_1(s)| \le C_7 \left( \frac{1}{n} + \frac{I_n^2(s)}{n} \right)  \le 2C_7\frac{I_n^2(s)}{n}.
$$
Similarly, the integrand in the fourth term in \eqref{ins1012} can be bounded as
\[\frac{Y^2(s)}{n^2} \sum_{i} \calC_i^4(s) \leq \frac{I^2_n(s)}{n} Y^2(s)\frac{1}{n} \sum_{i} \calC_i^2(s).\]
Noting that by definition $Y(s) \frac{1}{n} \sum_{i} \calC_i^2(s)=1$ and $Y(s)\leq 1$ we get 
\[\int_0^t \left(\frac{Y^2(s)}{n^2} \sum_{i} \calC_i^4(s) \right )ds \leq \int_0^{t_n} \frac{I^2_n(s)}{n}ds. \]
Combining we get, for all $t \le t_n$
\[|Y(t)-y(t)|\leq \eps_n+\int_0^t 2|Y(s)-y(s)|ds.\]
Thus Gronwall's lemma (see Theorem 5.1, Appendix in \cite{ethier-kurtz}) proves the result. 
\qed

One last ingredient in the proof of Proposition \ref{prop:s2-crit} is the following estimate on the martingale in Lemma \ref{lemma:mart-s2}.
\begin{Lemma}
\label{prop:n1-sup-M}
As $n \to \infty$, 
	\begin{equation}
		n^{1/3}\sup_{s\leq t_n} |M(s)| \convp 0. \notag
	\end{equation} 
\end{Lemma}
Proof of Lemma \ref{prop:n1-sup-M} is given in Section \ref{secins1049}.\\ \ \\

{\bf Proof of Proposition \ref{prop:s2-crit}: } In view of Lemma \ref{lemma:gronwall} it is enough to show $n^{1/3}\eps_n\convp 0$ as $n\to\infty$, where $\eps_n$ is as defined in the Lemma \ref{lemma:gronwall}. Let us analyze each of the terms in $\eps_n$. First note that by Lemma \ref{lemma:error-diff-eqn}, for any $\betapar<1/2$, and in particular for $\betapar=1/3$
\[\sup_{s\leq t_c} n^{\betapar}|\barx(s)-x(s)|\convp 0, \; \mbox{ as } n\to\infty.\]
Next, from Proposition \ref{prop:com-size}
$$ \{ I_n(t) \le m(n,t), \forall t < t_c-n^{-\gamma}\}  \mbox{ whp as } n \to \infty,$$
where
$ m(n,t)= B \frac{(\log n)^4}{(t_c-t)^2}$.
 Thus, recalling that $\gamma \in (1/6, 1/5)$, we have, whp, 
\begin{align}
	\int_0^{t_n} \frac{I^2_n(s)}{n}ds &\leq \frac{B^2(\log {n})^8}{n}\int_0^{t_n} \frac{1}{(t_c-s)^4} ds \nonumber \\
	&= \frac{B^2n^{3\gampar}(\log {n})^8}{n}=o\left(\frac{1}{n^{1/3}}\right). \label{ins1602}
\end{align}
The result now follows on combining these estimates with Lemma \ref{prop:n1-sup-M}.  \qed

\subsection{Proof of Lemma \ref{prop:n1-sup-M}} \label{secins1049}

We shall prove the lemma by first
showing that $n^{\betapar}\sup_{s\leq t_n}|M(s)| \convp 0$ for any $\betapar < 1/5$  and then sharpening the estimates to allow for any $\betapar \le 1/3$ . 
\begin{Lemma}
	\label{lemma:mart-bound-1}
	As $n \to \infty$, 
	\begin{equation}
	\prob\left(\sup_{s\leq t_n} |M_n(s)|> \frac{n^{3\gampar/2}(\log{n})^6}{\sqrt{n}} \right) \to 0	
	\label{eqn:mart-bound-1}
	\end{equation} 
	and
	\[\prob(Y_n(t)\leq 2 y(t) \; \forall  t< t_n)\to 1.\]
\end{Lemma} 
{\bf Proof:} 
Fix  $\gamma_1 \in (1/3, 1/2)$  and define  stopping times $\tau_i$, $i=1,2$  by
\[\tau_1 = \inf\{t: I_n^1(t) > m(n,t) \}, \;\; \tau_2= \inf\{t: |\bar{x}(t) - x(t)|> n^{-\gamma_1} \}.\]
From Proposition \ref{prop:com-size} and Lemma \ref{lemma:error-diff-eqn}  
\begin{equation}
	\prob(\tau_1\wedge \tau_2 > t_n) \to 1
	\label{eqn:tau-tn}
\end{equation}
as $n\to \infty$. Let $\tau^* = t_n\wedge \tau_1 \wedge \tau_2$. For \eqref{eqn:mart-bound-1}, in view of \eqref{eqn:tau-tn}, it suffices to prove the statement
with $t_n$ replaced by $\tau^*$.  
By Doob's maximal inequality we have 
\[\expec(\sup_{s\leq \tau^*}|M(s)|^2) \leq 4 \expec(\langle M\rangle(\tau^*)) = 4  \expec \int_0^{\tau^*} B(u) du .\]
Furthermore, from \eqref{eqn:yi-bound},  
\begin{align*}
	B(s) \leq \frac{4}{n}+  \frac{4 Y^2(s) I^2_n(s)}{n} \leq \frac{4}{n}+  \frac{4 I^2_n(s)}{n}. \notag
\end{align*}
Since for $t< \tau^*$ we have $I_n(t) \leq B(\log{n})^4/(t_c-t)^2 $ (see \eqref{ins1711}), we have 
$$
	\E [\langle M\rangle(\tau^*)] \leq \int_0^{t_n}\left ( \frac{4}{n}+ \frac{4B^2(\log{n})^8}{n(t_c-t)^4}\right) dt 
	  \leq d_1\frac{n^{3\gampar} (\log{n})^{10}}{n}.
$$
Combining the estimates, we have 
\[\expec(\sup_{s\leq \tau^*}|M(s)|) \leq d_1^{1/2}\frac{n^{3\gampar/2} (\log{n})^5}{\sqrt{n}}.\]
A simple application of Markov's inequality and \eqref{eqn:tau-tn} gives \eqref{eqn:mart-bound-1}. To get the final assertion in the Lemma, note that on the set
\[B_n = \set{\tau_1> t_n} \cap \set{\tau_2> t_n} \cap \set{\sup_{s\leq \tau^*} |M(s)| < \frac{n^{3\gampar/2}(\log{n})^6}{\sqrt{n}}}\]
$|\bar{x}(t) - x(t)|< n^{-\gamma_1}$ and (see \eqref{ins1602})
$$
	\int_0^{t_n} \frac{I^2_n(s)}{n}ds =o\left(\frac{1}{n^{1/3}}\right). $$
Therefore, since $\gampar <1/5$,
the error $\eps_n$ in Lemma \ref{lemma:gronwall} satisfies whp
\[\eps_n <\frac{n^{3\gampar/2}(\log{n})^7}{\sqrt{n}} = o\left(\frac{1}{n^\gampar}\right).\]
 Noting that $y(t)$ is a monotonically decreasing function,  the above along with \eqref{ins1607} implies that 
$Y(t) = (1+o_p(1))y(t)$ for $t\leq t_n$ on the set $B_n$. Since $\prob(B_n)\to 1$, the result follows. 
\qed \\ \ \\

{\bf Proof of Lemma \ref{prop:n1-sup-M}:} Along with the stopping times $\tau_1, \tau_2$ introduced in Lemma \ref{lemma:mart-bound-1}, consider the stopping time
\[\tau_3 = \inf\set{t: Y(t)> 2y(t)}\]
Then  Lemma \ref{lemma:mart-bound-1} and \eqref{eqn:tau-tn} imply that 
\[\prob(\tau_1\wedge \tau_2\wedge \tau_3 > t_n) \to 1\]
as $n\to\infty$. Thus to complete the proof, it is enough to show that for the stopping time $\tau = t_n\wedge \tau_1\wedge \tau_2\wedge \tau_3$,  $n^{1/3}\expec(\sup_{s\leq \tau} M(t)) \to 0$ as $n\to\infty$. Once again by Doob's maximal inequality it is enough to show that 
\be \label{ins1632b} n^{2/3}\expec(\langle M\rangle(\tau)) \to 0\qquad \mbox{as } n\to \infty. \ee
Now note that \eqref{eqn:yi-bound} implies that 
\begin{align*}
	\E[\langle M\rangle(\tau)]
	&\leq \E \int_0^{\tau} \left ( \frac{4}{n}+ \frac{4 Y^2(s) I^2_n(s) }{n} \right )ds \\
	&\leq \frac{4t_n}{n} + \frac{4}{n}\int_0^{t_n} 4y^2(s) \frac{2B^2(\log{n})^8}{(t_c-s)^4} ds \\
	&\le d_1 \left ( \frac{1}{n} + \frac{1}{n}\int_0^{t_n} \frac{B^2(\log{n})^8}{\alpha^2(t_c-s)^2} ds  + \frac{(\log{n})^8}{n}\right).
\end{align*}
In the second line of the above display we have used the fact that $I_n(t) \le m(n,t)$ for all $t \le \tau_1$ and in the last line we have used \eqref{eqn:yt-scaling}.
Thus
\[\E[n^{2/3}\langle M\rangle(\tau)] \leq d_2 \left (\frac{1}{n^{1/3}} + \frac{n^{2/3+\gampar}(\log{n})^8}{n} \right).\]
Since $\gampar <1/5$ we have \eqref{ins1632b} and this completes the proof.  \qed

\subsection{Analysis of $\bars_3$ near criticality}
We will now analyze the sum of cubes of component sizes near criticality and consequently prove \eqref{eqn:s3-s2}. Define
\[z(t) = \frac{s_3(t)}{s_2^3(t)}, \; t \in [0, t_c).\]
Then differential equations \eqref{ins-s2} \eqref{ins-s3} imply (see \cite{janson2010phase}) that $z$ solves the differential equation 
\begin{equation}
	\label{eqn:diff-z}
	z^\prime(t) = 3 x^2(t) y^3(t) - 3 x^2(t) y(t) z(t), \; z(0)= 1, \; t \in [0, t_c)
\end{equation}
and furthermore
$z(t)\to \beta$
as $t\to t_c$. Now consider the process 
\[Z_n(t)= \frac{\calS_3(t)/n}{(\calS_2(t)/n)^3}= Y^3_n(t) \frac{\calS_3(t)}{n}.\]
Then to show \eqref{eqn:s3-s2}, it is enough to show the following proposition:
\begin{Proposition}
	\label{prop:zt-convg}
Fix any $\gampar\in (1/6,1/5)$ and let as before $t_n = t_c-n^{-\gampar}$. Then 
\[|Z(t_n)-z(t_n)|\convp 0\]
as $n\to \infty$.
\end{Proposition}
 The analysis is similar to that for $\calS_2(\cdot)$  as carried out in the Section \ref{s2sec}. 
We begin by writing the semimartingale decomposition  for $Z_n$  and identifying the predictable quadratic variation  $\langle \ti M\rangle$ of the martingale $\ti M$ in the decomposition.  
\begin{Lemma}
	\label{lemma:mart-decomp-z}
	 The process $Z_n$ can be decomposed as 
	\be Z_n(t) = 1+ \int_0^t \ti A_n(s) ds + \ti M_n(t), \; t \in [0, t_c] \label{ins1006s3}\ee
	where
	\\(a)  $\ti M_n$ is a RCLL martingale with respect to the natural filtration $\{\FF_t\}_{t\geq 0}$ of the BF process. 
	\\(b) The process $\ti A_n = \ti A_1^n+\ti R_1^n + \ti R_2^n$ where (suppressing $n$), for $u \in [0, t_c]$,
	\[\ti A_1(u)= 3 \barx^2(u)Y^3(u) - 3Z(u)Y(u)\barx^2(u),\]
		 \be |\ti R_1(u)|=\left|  3 (1-\barx^2(u))\left[ Z(u)Y(u)\sum_{i}\frac{\calC_i^4(u)}{n^2}-Y^3(u)\sum_{i}\frac{\calC_i^5(u)}{n^2}\right]\right| \le \frac{6I^3(u)Y^2(u)}{n} \label{ins1110}
		 \ee
	and 
		\[ |\ti R_2(u)|\leq \frac{30Y^2(u)I^3_n(u)}{n}+ \frac{12 Y^3(u)I^5_n(u)}{n^2}. \]	
	\\(c) Predictable quadratic variation of $\ti M_n$ is given as
	\[\langle \ti M_n\rangle(t) = \int_0^t \ti B_n(u) du, \]
	and the process $\ti B_n$ can  be bounded as
	\begin{equation}
		\ti B(u)\leq C_9 \left(\frac{Y^4(u) I^4(u)}{n} + \frac{ Y^5(u)I^6(u)}{n^2}+ \frac{Y^6(u) I^8(u)}{n^3}+ \frac{Y^7(u) I^{10}(u)}{n^4} + \frac{Y^8(u) I^{12}(u)}{n^5}\right)
		\label{eqn:zr3-bound}\ee
		for some $C_9 \in (0, \infty)$.
	\end{Lemma} 
{\bf Proof:}  
We will proceed as in the proof of Lemma \ref{lemma:mart-s2}.  Note that 
$$Z(t) = 1 + \sum_{s \le t} \Delta Z(s), \mbox{ where } \Delta Z(s) = Z(s)- Z(s-).$$
Next note that if $\Delta \calS_2(u) = a$ and $\Delta \calS_3(u) = b$, then
$$
	Z(u) = \frac{\bars_3(u-)+b/n}{(\bars_2(u-)+ a/n)^3}
	= \frac{Y^3(u-)(\bars_3(u-)+b/n)}{(1+aY(u-)/n)^3}.$$
	Using the estimate $|(1+x)^{-3}- (1-3x)|\leq 6 x^2$, for $0<x<1$, we get
	$$Z(u) 
= Y^3(u-)\left(\bars_3(u-)+\frac{b}{n}\right)\left(1-\frac{3aY(u-)}{n}+\tilde{R}(a,u-)\right), $$
where 
\[|\tilde{R}(a,u-)|\leq \frac{6 a^2 Y^2(u-)}{n^2}.\]
Thus 
\[\Delta Z(u)\equiv \ti \zeta (a,b,u-)= -\frac{3aZ(u-)Y(u-) }{n}+ \frac{Y^3(u-)b}{n} + \ti R_{a,b}(u-),\]
where the remainder term 
\begin{align*}
	|\ti R_{a, b}(u-)|&\leq Z(u-)|\tilde{R}(a,u-)| + \frac{3abY^4(u-)}{n^2}+\frac{bY^3(u-)|\tilde{R}(a,u-)|}{n}\\
	&\leq \frac{6 a^2 Z(u-) Y^2(u-)}{n^2} + \frac{3abY^4(u-)}{n^2}+\frac{6a^2bY^5(u-)}{n^3}.
\end{align*} 
Any jump in $\calS_2, \calS_3$ or $Z$ corresponds to a jump of one of the Poisson processes $\PP_\bfe$,
$\bfe = (e_1,e_2) \in \EE^2$.   A jump of $\PP_\bfe$ at a time instant $u$ could result in the following different values for $a,b$.\\

(i) {\bf Merger caused by the first edge $e_1$:} In this case $a=2$ and $b=6$ and $\ti R_{a,b}$ can be estimated as
$$
	|\ti R_{2,6}(u-)| 
	\le \frac{24 Y^4(u-) I(u-)}{n^2} + \frac{36 Y^4(u-)}{n^2}+\frac{144 Y^5(u-)}{n^3}
	\le \frac{204 Y^4(u-) I(u-)}{n^2}.$$
(ii) {\bf Merger caused by the second edge $e_2$:} In this case, suppose components $i$ and $j$ merge, then
 $$a\equiv \theta_{i,j}(u-) = 2\calC_i(u-)\calC_j(u-) , \;  b \equiv \eta_{i,j}(u-)=3\calC_i^2(u-)\calC_j(u-)+ 3\calC_i(u-)\calC_j^2(u-)$$ and noting that
 $a \le 2 I^2$, $b \le 6 I^3$ and
\begin{equation}
	Z(u)=Y^3(u) \bars_3(u)\leq Y^2(u) I_n(u), \label{eqn:zu-bound}
\end{equation}
 we have
$$
|\ti R_{a, b}(u-)|
	\leq \frac{6 a^2 Z(u-) Y^2(u-)}{n^2} + \frac{3abY^4(u-)}{n^2}+\frac{6a^2bY^5(u-)}{n^3}
	\leq \frac{30 a Y^4(u-) I^3(u-)}{n^2}+\frac{72 a Y^5(u-) I^5(u-)}{n^3}.
$$
With these observations we can represent $Z$ in terms of stochastic integrals with respect to $\PP_\bfe$ as follows.  Recall
$
\clh_1,\clh_2^{(i,j)}$ introduced in the proof of Lemma \ref{lemma:mart-s2}.
Define $\ti \alpha_1(u) = \ti \zeta(2,6,u)$ and  $\ti \alpha_2^{i,j}(u) = \ti \zeta (\theta_{i,j}(u),\eta_{i,j}(u),u)$.  Also, let
$$\ti \clu_{\bfe}(u) = \ti \alpha_{1}(u) {\bf 1}_{\clh_1(u)}(\bfe), \; \ti \clu_{\bfe}^{i,j}(u) = \ti \alpha_2^{i,j}(u){\bf 1}_{\clh_2^{(i,j)}(u)}(\bfe).$$
Then
\be \label{ins2023s3}
Z(t) = 1 + \sum_{\bfe \in \EE^2} \int_{(0, t]} \left (\ti \clu_{\bfe}(s-) + \sum_{i< j} \ti \clu_{\bfe}^{i,j}(s-)\right) \PP_{\bfe}(ds).\ee
Recalling that  $\PP_{\bfe}$ is a rate $2/n^3$ Poisson process, one can write $Z$ as
$$Z(t) = 1 + \int_{[0,t]} \ti A(s) ds + \ti M(t),$$
where
$$
\ti A(s) = \frac{2}{n^3} \sum_{\bfe \in \EE^2}\left (\ti \clu_{\bfe}(s) + \sum_{i< j} \ti \clu_{\bfe}^{i,j}(s)\right).$$
Also, once again using independence of Poisson processes $\PP_{\bfe}$,
$$
\langle \ti M\rangle(t) = \frac{2}{n^3}\sum_{\bfe \in \EE^2} \int_{(0, t]} \left ((\ti \clu_{\bfe}(s))^2 + \sum_{i< j} (\ti \clu_{\bfe}^{i,j}(s))^2 \right) ds
\equiv \int_{(0,t]} \ti B(s) ds.$$
The proof is now completed upon using \eqref{ins2029a} and \eqref{ins2029b} as for Lemma  \ref{lemma:mart-s2}.
\qed\\

As is clear from the above lemma, a precise analysis of $Z$ will  involve considering several terms of the form $I^\thtpar Y^\betapar$.
The following lemma shows that such terms are asymptotically negligible for suitable $\thtpar, \betapar$. 
\begin{Lemma}
	\label{lemma:powers-bd}
	For any $\thtpar, \betapar\geq 0$ and $p > 0$ satisfying $\gamma(2\thtpar-\betapar -1)< p$
\[\int_0^{t_n}\frac{I^\thtpar_n(u) Y^\betapar_n(u)}{n^p} du \convp 0,\]
	as $n\to\infty$. 
	\end{Lemma}  
{\bf Proof:} From Proposition \ref{prop:com-size}, 
$$  \{I_n(t) \leq B(\log n)^4/(t_c-t)^{2} \mbox{ for all } t\leq t_n \} \mbox{ occurs whp}. 
$$
Also, by Proposition \ref{prop:s2-crit} and \eqref{ins1607} we know that 
\[\sup_{t\leq t_n}|Y(t)- y(t)|=o_p(y(t_n))\] 
and from \eqref{eqn:yt-scaling}, for $t$ near $t_c$, $y(t_c)\sim (t_c-t)/\alpha$. Using these bounds in the above integral proves the result. 
\qed \\ \ \\

{\bf Proof of Proposition \ref{prop:zt-convg}. }By the integral representation of the process $Z$ given in Lemma \ref{lemma:mart-decomp-z} and the fact that $z$ solves the differential equation \eqref{eqn:diff-z}, we have that 
\begin{align}
	|Z(t)-z(t)|\leq 3 \int_0^{t} |\barx^2(u)Y^3(u)&-x^2(u)y^3(u)|du+ 3\int_0^{t} \left|Z(u)Y(u)\barx^2(u)-z(u)y(u)x^2(u)\right|du \nonumber \\
	&+ \int_0^{t_n} |\ti R_1(u)|du + \int_0^{t_n} |\ti R_2(u)|du +\sup_{t\leq t_n} |\ti M(t)|. \label{ins1057}
\end{align}
Now the integrand in the first term can be bounded as 
\begin{align}
	3|\barx^2(u)Y^3(u)-x^2(u)y^3(u)|&\leq 3Y^3(u)|\barx^2(u)-x^2(u)|+ 3\barx^2(u)|Y^3(u)-y^3(u)| \notag\\
	&\leq 6|\barx(u)-x(u)|+ 9|Y(u)-y(u)|.\label{eqn:xy-bound}
\end{align}
The integrand in the second term in \eqref{ins1057} can be decomposed  as 
\begin{align*}
z(u)y(u)x^2(u)-	Z(u)Y(u)\barx^2(u) = x^2(u)&y(u)(z(u)-Z(u))\\
   & + y(u)Z(u)(x^2(u)-\barx^2(u))+ Z(u)\barx^2(u)(y(u)-Y(u)).
\end{align*}
Thus the second integral in \eqref{ins1057} can be bounded by 
\begin{equation}
	3 \int_0^t |(z(u)-Z(u))|du+ 6 \int_0^{t_n} Z(u)|x(u)-\barx(u)|du + 3 \int_0^{t_n}Z(u)|Y(u)-y(u)|du.
	\label{eqn:integ-z-bd}
\end{equation}
Combining \eqref{eqn:integ-z-bd} and \eqref{eqn:xy-bound} we get that 
\[|Z(t)-z(t)| \leq \eps_n+ 3 \int_0^t |(Z(u)-z(u))|du\]
where 
\begin{align*}
	\eps_n = 9t_c\sup_{s\leq t_n}\left(|\barx(s)-x(s)|+|Y(s)-y(s)|\right) &+ 6 \int_0^{t_n} Z(u)|x(u)-\barx(u)|du + 3 \int_0^{t_n}Z(u)|Y(u)-y(u)|du \\
	&+ \int_0^{t_n} |\ti R_1(u)|du + \int_0^{t_n} |\ti R_2(u)|du +\sup_{t\leq t_n} |\ti M(t)|\\
	&\hspace{-2in} = \eta_1+\eta_2+\eta_3+\eta_4+ \eta_5 + \eta_6.
\end{align*}
By Gronwall's lemma, it is enough to show that $\eps_n\to 0$ in probability as $n\to\infty$. Let us show each of the six constituents of $\eps_n$ satisfy this asymptotics. By Lemma \ref{lemma:error-diff-eqn} and Proposition \ref{prop:s2-crit}~ $\eta_1\to 0$ in probability. Again by Lemma \ref{lemma:error-diff-eqn}, for any $\betapar< 1/2$, whp, 
$$
	\eta_2 \leq \frac{6}{n^\betapar}\int_0^{t_n} Z(u) du 
	\leq \frac{6}{n^\betapar} \int_0^{t_n} Y^2(u) I_n(u) du.$$
	Using Lemma \ref{lemma:powers-bd}, the last term converges to $0$ in probability as $n \to \infty$.  Thus $\eta_2 \to 0$ in probability.
	
An identical argument, using Proposition \ref{prop:s2-crit} instead of Lemma  \ref{lemma:error-diff-eqn} shows that $\eta_3\to 0$ in probability. For $\eta_4$, 
note that from \eqref{ins1110},\[|R_1(u)|\leq \frac{6 I^3(n,u)Y^2(u)}{n}.\] 
Lemma \ref{lemma:powers-bd} now shows that $\eta_4\to 0$ in probability. A similar argument, using the bounds in Lemma \ref{lemma:mart-decomp-z} on $R_2(u)$ establishes that $\eta_5\to 0$ in probability. 
For $\eta_6$ note that for an arbitrary stopping time $\tau$
\[ \E \sup_{t \le t_n\wedge \tau} |\ti M_t|^2 \le 4 \E[\langle M\rangle(t\wedge \tau)] = 4\E \int_0^{t\wedge \tau} B(u) du.\]
The bound on $\ti B(u)$ in \eqref{eqn:zr3-bound} along with Lemma \ref{lemma:powers-bd} and a localization argument similar
to the one used in the proof of Lemma \ref{prop:n1-sup-M} now shows that $\eta_6$  converges to $0$ in probability.

\subsection{Proof of Proposition \ref{prop:main}}
We now complete the proof of Proposition \ref{prop:main}.
Proof of \eqref{eqn:s2-tn-n-alpha} follows from Proposition \ref{prop:s2-crit} and the discussion immediately above the proposition.
Proof of \eqref{eqn:s3-s2} is immediate from Proposition \ref{prop:zt-convg}.  Finally we consider \eqref{eqn:max-s2}.
From Proposition \ref{prop:s2-crit} and \eqref{ins1607} 
\be \frac{\calS_2(t_n)}{\alpha n^{1+\gampar}} \convp  1\label{ins2038} \ee
Also, from Proposition \ref{prop:com-size} 
\[\prob\left(\frac{\calC_n^{\sss(1)}(t_n)}{ n^{2\gampar}\log^4{n}}\leq B\right)\to 1\]
as $n\to\infty$.  Combining, and recalling that $\gampar\in (1/6, 1/5)$, we have 
\[\frac{n^{2/3}\calC_n^{\sss (1)}(t_n)}{\calS_2(t_n)} \convp 0.\]
This completes the proof of Propostion \ref{prop:main}. \qed

\section{Proof of Theorem \ref{theo:main}}
\label{sec:proof-main}
We  will now complete the proof of  Theorem \ref{theo:main}. As always, we write the component sizes as
$$\bfC_n^{\sss BF}(t)\equiv (\CC_n^{\sss (i)}(t) : i \ge 1)\equiv  (\CC_i(t): i\ge1 );$$
and write the scaled component sizes as
\be
\label{ins1407}
\bar \bfC_n^{\sss BF}(\lambda)\equiv  \left(\frac{\beta^{1/3}}{n^{2/3}}\calC^{\sss (i)}_n \left(t_c+ \beta^{2/3} \alpha\frac{\lambda}{n^{1/3}} \right) : i \ge 1 \right) \equiv \left(\bar{\calC}_i(\lambda) : i\ge 1 \right) \ee
Then Proposition \ref{prop:main} proves that  with
\[\lambda_n = -\frac{n^{-\gamma+1/3}}{\alpha \beta^{2/3}}\]
and  $\gamma \in (1/6,1/5)$ we have, as $n \to \infty$,
\be \label{ins1409}	\frac{\sum_i \left( \bar \CC_i(\lambda_n)\right)^3}{\left[\sum_i \left( \bar \CC_i (\lambda_n)\right)^2\right]^3} \convp 1, \;
	\frac{1}{\sum_i \left( \bar \CC_i (\lambda_n)\right)^2} +\lambda_n \convp 0, \; 
	\frac{ \bar \CC_1(\lambda_n)}{\sum_i \left( \bar \CC_i(\lambda_n)\right)^2} \convp 0. \ee

We shall now give an idea of the proof of the main result, and postpone precise arguments to the next two sections. 
The first step is to observe that the asymptotics in \eqref{ins1409} imply that the  $\bar {\bfC}_{\sss bf}$ process  at time $\lambda_n$ satisfies the  regularity conditions of Proposition 4 of \cite{aldous1997brownian}. 
The second key observation is that the scaled components merge in the critical window at a rate close to that for the
multiplicative coalescent.  Indeed,  note that for any given time $t$  components $i<j\in \BF(t)$  merge in a small time interval $[t, t+dt)$ at rate
\[\frac{1}{n}(1-\barx^2(t)) \calC_i(t)\calC_j(t).\] 
Thus letting $\lambda = (t-t_c)n^{1/3} / (\alpha \beta^{2/3})$ be the scaled time parameter, 
 in the time interval $[\lambda, \lambda+d\lambda)$, these two components merge at rate 
\begin{align*}
	\gamma_{ij}(\lambda)&=\frac{(1-\barx^2(t_c+ \beta^{2/3} \alpha\frac{\lambda}{n^{1/3}}))}{n} \frac{\beta^{2/3} \alpha}{n^{1/3}} \CC_i\left(t_c+ \frac{ \beta^{2/3} \alpha\lambda}{n^{1/3}}\right)\CC_j \left(t_c+ \frac{ \beta^{2/3} \alpha\lambda}{n^{1/3}}\right) \\
	&= \alpha \left(1-\barx^2\left(t_c+ \beta^{2/3} \alpha\frac{\lambda}{n^{1/3}}\right)\right) \bar \CC_i(\lambda)  \bar \CC_j(\lambda).
\end{align*} 
Now since, for large $n$,
\[\barx^2\left(t_c+ \beta^{2/3} \alpha\frac{\lambda}{n^{1/3}}\right)\approx x^2(t_c)\]
and from \cite{janson2010phase}, $\alpha(1-x^2(t_c)) =1$ (see \eqref{eqn:alpha-def}) we get 
\[\gamma_{ij}(\lambda) \approx \bar \CC_i(\lambda)  \bar \CC_j(\lambda)\] 
which is exactly the rate of merger for the multiplicative coalescent. The above two facts  allow us to complete the proof using ideas similar to those in \cite{aldous2000random}.
Let us now make these statements precise. 

As before, throughout this section  $t_n=t_c-n^{-\gamma}=t_c+\beta^{2/3} \alpha\frac{\lambda_n}{n^{1/3}}$, where $\gamma$ is fixed in $(1/6, 1/5)$.
We will first show that $\bar \bfC_n^{\sss BF}(\lambda) \convd \bfX(\lambda)$ in $\ldown$ for each $\lambda \in \RRR$ and at the end of the section show that, in fact,
$\bar \bfC_n^{\sss BF} \convd \bfX$ in $\DD((-\infty, \infty): \ldown)$.  Now fix $\lambda \in \RRR$.  By choosing $n$ large enough we can ensure that
$\lambda \ge \lambda_n$.  Henceforth consider only such $n$.  Recall that 
$\com_n(t)$ denotes the subgraph of $\BF_n(t)$ obtained by deleting all the singletons. Let $ \sum_{i \in \com}$ denote the summation over all components in $\com_n$, and $\sum_i$ denote the summation over all components in $\BF_n$.
Since
\be
\sum_i \left (  \bar \CC_i(\lambda) \right)^2 - \sum_{i \in \com} \left (  \bar \CC_i(\lambda) \right)^2 \le \frac{d_1}{n^{4/3}} \sum_{i=1}^{X_n(t)} 1 = O(1/n^{1/3}),
\label{ins1503}\ee
it suffices to prove Theorem \ref{theo:main} and verify Proposition \ref{prop:main} with $\BF_n(t)$ replaced by $\com_n(t).$ We write $\sum_i$ instead of $\sum_{i \in \com}$ for simplicity of the notation from now on. We begin in Section \ref{ins2122} by constructing a coupling
of $\{\com_n(t)\}_{t \ge t_n}$ with two other random graph processes, sandwiching our process between these two processed,  and proving statements analogous to those in Theorem  \ref{theo:main} for scaled component vectors
associated with these processes. Proof of Theorem  \ref{theo:main} will then be completed in Section \ref{ins2121}.\\ \ \\

\subsection{Coupling with the multiplicative coalescent} \label{ins2122}


\textbf{Lower bound coupling:} Let, for $t\ge t_n$, $\com_n^-(t)$ be a modification of $\com_n(t)$ such that $\com_n^-(t_n)=\com_n(t_n)$,  and when $t > t_n$, we change the dynamics of the random graph to the Erd\H{o}s-R\'{e}nyi type. More precisely, recall from Section \ref{sec:model-equiv}
that a jump in $\BF_n(t)$ can be produced by three different kinds of events.  These are described in items (i), (ii) and (iii) in Section \ref{sec:model-equiv}.
$\com^-_n(t)$, $t \ge t_n$ is constructed from $\com^-_n(t_n)$ by erasing events of type (i) and (ii) (i.e. immigrating doubletons and attaching singletons)
and changing the probability of edge formation between two non-singletons (from that given in \eqref{ins1516}) to the fixed value $b_n^*(t_n)/n$. Since
$b_n^*(t)$ is nondecreasing in $t$, we have that $\com_n(t_n + \cdot) \ge_d \com_n^-(t_n + \cdot)$.
Denote by $\bar{\bfC}_n^-(\lambda) = \left( \bar \calC_i^-(\lambda): i\geq 1\right)$ the scaled (as in \eqref{ins1407}) 
component size vector for $\com_n^-(t)$. 
From Proposition 4 of \cite{aldous1997brownian}, it follows that for any $\lambda \in \Rbold$,
\be \label{ins1320} \bar {\bfC}_n^-(\lambda) \convd \bfX(\lambda) \ee
in  $\ldown$.
Indeed, note that the first and third convergence statements in \eqref{ins1409} hold with $\bar \CC_i$ replaced with  $\bar \calC_i^-$
since the contributions made by singletons to the scaled sum of squares is $O(n^{-1/3})$ (see \eqref{ins1503}) and to the sum of cubes is even smaller.  This
shows that the first and third requirements in Proposition 4 of  \cite{aldous1997brownian} (see equations (8), (10) therein) are met.
To show the second requirement in Proposition 4 of  \cite{aldous1997brownian}, using the second convergence in \eqref{ins1409},
\beqn
&&\lim_{n\to \infty} \left ( \left ( n^{2/3}\beta^{-1/3}\right)^2 \frac{b_n^*(t_n)}{n}  \frac{\beta^{2/3}\alpha (\lambda - \lambda_n)}{n^{1/3}} -
\frac{1}{\sum_i \left( \bar \calC_i^-(\lambda_n)\right)^2}\right)\label{ins1608}\\
&=& \lim_{n\to \infty} \alpha b_n^*(t_n)\lambda - \lambda_n (\alpha b_n^*(t_n) - 1)\nonumber\\
&=& \lambda - \lim_{n\to \infty} \lambda_n (\alpha b_n^*(t_n) - 1),\nonumber
\eeqn
where the last equality follows on observing that, as $n \to \infty$,
$b_n^*(t_n) \convp 1-x^2(t_c)$ and $\alpha(1-x^2(t_c))=1$.  Also,
\beq
 \lim_{n\to \infty} \lambda_n |\alpha b_n^*(t_n) - 1|  &=& \lim_{n\to \infty}\frac{n^{-\gamma+1/3}}{\beta^{2/3}} |b_n^*(t_n) - \alpha^{-1}|\\
 &=&  \lim_{n\to \infty}\frac{n^{-\gamma+1/3}}{\beta^{2/3}} |b_0(\bar x(t_n)) - b_0(x(t_c))| \\
 &\le & d_1  \lim_{n\to \infty}n^{-\gamma+1/3} |\bar x(t_n)- x(t_c)|\\
 &\le &  \lim_{n\to \infty} d_2 \left (n^{-\gamma+1/3} |\bar x(t_n) - x(t_n)|  + n^{-\gamma+1/3}|t_n - t_c| \right),
 \eeq
 where the second equality follows from \eqref{eqn:b-def1}.    The first term on the last line converges to $0$ using Lemma \ref{lemma:error-diff-eqn}.
 For the second term note that $n^{-\gamma+1/3}|t_n-t_c| = n^{-\gamma+1/3}n^{-\gamma}$ which converges to $0$ since $\gamma > 1/6$.
Thus we have shown that the expression in \eqref{ins1608} converges to $\lambda$ as $n \to \infty$ and therefore the second requirement
in Proposition 4 of  \cite{aldous1997brownian}  (see equation 9 therein) is met as well.  This proves that $\bar{\bfC}_n^-(\lambda) \convd \bfX(\lambda)$
in  $\ldown$, for  every $\lambda \in \RRR$.  Although Proposition 4 of \cite{aldous1997brownian} only proves convergence at any fixed point $\lambda$,
from the Feller property of the multiplicative coalescent process proved in Proposition 6 of the same paper it now follows that, in fact,
$\bar{\bfC}_n^- \convd \bfX$ in $\DD((-\infty, \infty): \ldown)$.\\ \ \\

\textbf{Upper bound coupling:} Let us construct $\{\com_n^+(t):  t \ge t_n \}$ in the following way. Let $t_n^+ = t_c + n^{-\gamma}$ and let
$$\lambda_n^+ = (t_n^+ - t_c) n^{1/3}/(\alpha \beta^{2/3}) = n^{1/3 - \gamma}/(\alpha \beta^{2/3}).$$
Let $\com_n^+(t_n)$ be the graph obtained by including all immigrating doubleton and attachments during time $t \in [t_n, t_n^+]$ to the graph of $\com_n(t)$, along with all the attachment edges. Namely, we construct
$\com_n^+(t_n)$ by including in $\com_n(t_n)$ all events of type (i) and (ii) of Section \ref{sec:model-equiv} that occur over $[t_n, t_n^+]$.
For $t > t_n$ the graph evolves in the Erd\H{o}s-R\'{e}nyi way such that edges are added between each pair of vertices in the fixed rate $b_n^*(t_n^+)/n$. The coupling between $\com_n^+(\cdot + t_n)$ and $\com_n(\cdot + t_n)$ can be achieved as follows:
 Construct a realization of $\{\com_n(t): t_n \le t \le t_n^+ \}$ first, then use $b_n^*(t_n^+)-b_n^*(t)$ to make up for all the additional edges in $\com_n^+(t)$ for $t_n \le t \le t_n^+$. Note that $\com_n(t_n+ \cdot) \le_d \com_n^+(t_n+\cdot)$ over $[0, t_n^+-t_n]$. \\
Let $\bar \bfC_n^+(\lambda) = \left(\bar \calC_i^+(\lambda): i\geq 1\right)$ be the scaled (as in \eqref{ins1407}) 
component size vector for $\com_n^+$.  We will once more apply Proposition 4 of \cite{aldous1997brownian}. We first show that  the three convergence statements in \eqref{ins1409} hold with $\bar \CC_i$ replaced with  $\bar \calC_i^+$.  For this it will be convenient to consider processes
under the original time scale.  Write $\calC_{n}^{\sss(i)}(t_n) \equiv \cle_i$.  Also denote by $\{\cle_i^+\}$ the component vector obtained by adding
all events of type (ii) only, to $\com_n(t_n)$ (i.e. attachment of singletons to components in $\com_n(t_n)$), over $[t_n, t_n^+]$. Since
$c^*$ is bounded by $1$,  $\cle_i^+$ is stochastically dominated by the sum of $\cle_i$ independent copies of Geometric($p$), with $p=e^{t_n-t_n^+}=e^{-2n^{-\gamma}}$. Thus
\[
u_i \stackrel{\sss def}{=} \cle_i^+-\cle_i \le_d \text{Negative-binomial}(r,p)\text{ with } r=\cle_i, p=e^{-2n^{-\gamma}}.
\]
The random graph $\com_n^+(t_n)$ contains components other than $\{\cle_i^+\}$.  These additional components correspond to the ones
obtained from doubletons immigrating over $[t_n, t_n^+]$.  Since there are at most $n$ vertices, the number $N$ of such doubletons is bounded by $n/2$. Denote by $\{\tilde \cle_i^+\}_{i=1}^N$ the components corresponding to such
doubletons. Once again using the fact that $c^* \le 1$, we have that
\[
\ti \cle_i^+ \le_d  2+ \text{Negative-binomial}(2,p)\text{ with }  p=e^{-2n^{-\gamma}}.\]
 Write 
 $$\calS_k=\sum_i (\cle_i)^k,\, \calS_k^+=\sum_i (\cle_i^+)^k+ \sum_{i=1}^N (\tilde \cle_i^+)^k \mbox{ for } k=2, 3 \mbox{ and }
  I = \max_i \cle_i, I^+=\max\{ \max_i \cle_i^+, \max_i \tilde \cle_i^+\}.$$
  The following proposition shows that Propostion \ref{prop:main} holds with $(\calS_2(t_n), \calS_3(t_n), \clc_n^{\sss (1)}(t_n))$ replaced with
  $(\calS_2^+(t_n), \calS_3^+(t_n), I^+(t_n))$.
 \begin{Proposition}
As $n \to \infty$,
\begin{align*}
I^+ &= \Theta(I)\\
\frac{\calS_2^+}{\calS_2} &\convp 1\\
\frac{\calS_3^+}{\calS_3} &\convp 1\\
n^{4/3}\left( \frac{1}{\calS_2}-\frac{1}{\calS_2^+} \right) &\convp 0.
\end{align*}
\end{Proposition}
\textbf{Proof:} 
An elementary calculation shows that if $U$ is $\text{Negative-binomial}(r,e^{-2n^{-\gamma}})$ then for some $d_1 \in (0, \infty)$
$$\prob (U \ge 3\gamma^{-1} r) \le \frac{d_1}{n^3}$$
and thus, as $n \to \infty$,
$$\prob(\max_i \cle_i^+ \ge (1+ 3\gamma^{-1})I) \le \prob(u_i \ge 3\gamma^{-1}\cle_i \mbox{ for some } i = 1, \cdots n) \to 0.$$
A similar calculation shows that, for some $d_2 \in (0, \infty)$, as $n \to \infty$.
$$\prob(\max_{i= 1, \cdots N} \tilde \cle_i^+  \ge d_2) \to 0.$$
 The first statement in the proposition now follows on combining the above two displays.

Next, note that for Negative-binomial($r,p$), the first, second and third moments are
\begin{align*}
M_1 &=\frac{1}{p} r(1-p)\\
M_2 &=\frac{1}{p^2}[r^2(1-p)^2+r(1-p)]\\
M_3 &=\frac{1}{p^3}[r^3(1-p)^3+3r^2(1-p)^2+r(4-9p+7p^2-2p^3)].
\end{align*}
From \eqref{ins2038} and \eqref{eqn:s3-s2} it follows that
 $\calS_2 = \Theta(n^{1+\gamma})$ and $\calS_3=\Theta(n^{1+3\gamma})$.  Also, clearly,   $\sum_i \CC_i =O(n)$.

Write $D_2 \stackrel{\sss def}{=} \calS_2^+-\calS_2=\sum_{i=1}^N (\tilde \cle_i^+)^2  + \sum_i(2\CC_i u_i + u_i^2)$, then 
$$\E[D_2 |\{ \CC_i\}_i] \le d_2\left( n\cdot n^{-\gamma}+\sum_i[ (\CC_i)^2 n^{-\gamma} +(\CC_i)^2 n^{-2\gamma}+\CC_i n^{-\gamma} ] \right)=O(n)$$
thus $D_2/\calS_2 \convp 0$ and consequently $\calS_2^+/\calS_2 \convp 1$.

Write $D_3 \stackrel{\sss def}{=} \calS_3^+ -\calS_3=\sum_{i=1}^N (\tilde \cle_i^+)^3 + \sum_i[3(\CC_i)^2 u_i + 3 \CC_i u_i^2 + u_i^3 ]$.
One can similarly show that 
$$ \E[D_3 | \{\CC_i\}_i] =O(n^{1+2\gamma})$$ 
thus $D_3/\calS_3 \convp 0$ and  so $\calS_3^+/\calS_3 \convp 1$.

To prove the third convergence, it suffices to prove
\be \label{ins2046} \frac{n^{4/3} D_2}{(\calS_2)^2} \convp 0. \ee
By the asymptotics shown above, we have
$$\frac{n^{4/3} D_2}{(\calS_2)^2} = O( n^{4/3+1-2(1+\gamma) })=O(n^{1/3-2\gamma})$$
As $\gamma > 1/6$, \eqref{ins2046} follows and thus the proof is completed. \qed\\

For scaled
component size vector  of $\com_n^+$, the above proposition shows that the  statements in \eqref{ins1409} hold with $\bar \CC_i$ replaced with  $\bar \calC_i^+$.
In particular, the first and third requirements in Proposition 4 of  \cite{aldous1997brownian}  are met by $\{\bar \calC_i^+\}$
Also, using the second convergence in \eqref{ins1409}, a calculation similar to that for \eqref{ins1608} shows that
$$
\lim_{n\to \infty} \left ( \left ( n^{2/3}\beta^{-1/3}\right)^2 \frac{b_n^*(t_n^+)}{n}  \frac{\beta^{2/3}\alpha (\lambda - \lambda_n)}{n^{1/3}} -
\frac{1}{\sum_i \left( \bar \calC_i^+(\lambda_n)\right)^2}\right) \to \lambda .$$

 Therefore the second requirement
in Proposition 4 of  \cite{aldous1997brownian}   is satisfied.  This proves that 
\be \label{ins1318} \bar \bfC_n^+(\lambda) \convd \bfX(\lambda). \ee
in  $\ldown$, for  every $\lambda \in \RRR$.  Using  Proposition 6 of \cite{aldous1997brownian}  once again it now follows that
$\bar \bfC_n^+ \convd \bfX$ in $\DD((-\infty, \infty): \ldown)$.

\subsection{Completing the proof of Theorem \ref{theo:main}}\label{ins2121}

By \cite{aldous1998entrance, aldous2000random}, there is a natural partial order $\preceq$ on $\ldown$. Informally, 
interpreting an element of $\ldown$ as a sequence of cluster sizes, $\bfx,\bfy \in \ldown$, $\bfx \preceq \bfy$ if $\bfy$ can be obtained from $\bfx$ by adding new clusters and coalescing together clusters. The coupling constructed in Section \ref{ins2122} gives that, for every, $\lambda \in (\lambda_n, \lambda_n^+)$
$$ \bar{\bfC}_n^-(\lambda) \preceq \bar\bfC_n^{\sss BF}(\lambda) \preceq \bar \bfC_n^+(\lambda).$$
Since, as $n \to \infty$,  $\lambda_n \to - \infty$ and $\lambda_n^+ \to +\infty$, \eqref{ins1320}, \eqref{ins1318} along with Lemma 15 of \cite{aldous2000random}
yield that
$$\bar\bfC_n^{\sss BF}(\lambda) \convd \bfX(\lambda)$$
for all $\lambda\in \RRR$.

Finally we argue convergence in $\DD((-\infty, \infty):\ldown)$.
For $\bfx, \bfy \in \ldown$, let $\bfd^2(\bfx, \bfy) = \sum_{i=1}^{\infty} (x_i-y_i)^2$, $\bfx = \{x_i\}$, $\bfy = \{y_i\}$. Then $\bfd^2(\bfx,\bfy) < \sum_i y_i^2 -\sum_i x_i^2$
whenever $\bfx \preceq \bfy$. To prove that $\bar\bfC_n^{\sss BF} \to \bfX$ in $\DD((-\infty, \infty):\ldown)$ it suffices to prove that 

\be \label{ins1338}
\sup_{\lambda \in [\lambda_1, \lambda_2]} \bfd(\bar \bfC_n^{\sss BF}, \bar{\bfC}_n^-) \convp 0, \mbox{ for all} -\infty < \lambda_1 < \lambda_2 < \infty .\ee
Fix $\lambda_1, \lambda_2$ as above.  Then
\be
\label{ab1339}
\sup_{\lambda \in [\lambda_1, \lambda_2]} \bfd(\bar \bfC_n^{\sss BF}, \bar{\bfC}_n^-) \le \sup_{\lambda \in [\lambda_1, \lambda_2]}
[\sum_i (\bar \CC_i^+(\lambda))^2-\sum_i (\bar \CC_i^-(\lambda))^2].
\ee
Let, for $\lambda \in \RRR$, 
$$\clu_+(\lambda) = \sum_i (\bar \CC_i^+(\lambda))^2, \; \clu_-(\lambda) = \sum_i (\bar \CC_i^-(\lambda))^2 \mbox{ and } \clv(\lambda) = \clu_+(\lambda)- \clu_-(\lambda).$$
From Lemma 15 of
\cite{aldous2000random}, $\clv(\lambda) \convp 0$ for every $\lambda \in \RRR$.  Thus it suffices to show that $\clv$ is tight in $\DD((-\infty, \infty): \RRR_+)$.
Note that both $\clu_+$ and $\clu_-$ are tight in $\DD((-\infty, \infty): \RRR_+)$.  Although, in general difference of relatively compact sequences in the $\DD$-space need not be relatively compact, in the current
setting due to properties of the multiplicative coalescent this difficulty does not arise.  Indeed, if $\{\bfX^\bfx(t), t \ge 0\}$ denotes the multiplicative coalescent on the positive real line with initial condition
$\bfx \in \ldown$ then, for $\delta$ sufficiently small 
$$\sup_{\tau \in \clt(\delta)}\E \left(\bfd^2(\bfX^\bfx(\tau) , \bfx)\wedge 1\right)\le \E \left[\sum_i(X^\bfx_i(\delta))^2 -\sum_i x_i^2\right] \le 2\sum_{i<j} \delta x_ix_j \cdot 2 x_i x_j \le 2\delta||\bfx||^4,$$
where, $||\bfx||= (\sum x_i^2)^{1/2}$, $\clt(\delta)$ is the family of all stopping times (with the natural filtration) bounded by $\delta$.
Using the above property, the Markov property of the coalescent process and the tightness of $\sup_{\lambda \in [\lambda_1, \lambda_2]}\clu_+(\lambda)$, $\sup_{\lambda \in [\lambda_1, \lambda_2]}\clu_-(\lambda)$ one can verify Aldous' tightness criteria (see Theorem VI.4.5 in  \cite{jacod-shiryaev}) for $\clv$ thus proving the desired tightness.
\qed


\vspace{1in}

{\bf Acknowledgements}
AB and XW have been supported in part by the National Science Foundation (DMS-1004418), the Army Research Office
(W911NF-0-1-0080, W911NF-10-1-0158) and the US-Israel Binational Science Foundation (2008466). SB's research has been supported in part by the UNC research council and a UNC Junior Faculty development award. SB would like to thank SAMSI for many interesting discussions with the participants of the complex networks program held at SAMSI in the year 2010-2011. 

\bibliographystyle{plain}

\end{document}